%% file: 2019_10_07.tex
\documentclass[10pt,twoside,english,reqno,a4paper]{amsart}

\usepackage{geometry}
\usepackage{datetime}
\usepackage{fancyhdr}
\usepackage{amsmath}
\newlength{\defbaselineskip}
\input{style.tex}
\input{macro.tex}

\input{todo.tex}

\title{\bf Regularizing effects concerning elliptic equations\\ with a superlinear gradient term}

\author{Marta Latorre, Martina Magliocca and Sergio Segura de Le\'on}

\address{M. Latorre: Departamento de Matem\'atica Aplicada, Ciencia y Tecnolog\'ia de los Materiales y Tecnolog\'ia Electr\'onica,
Universidad Rey Juan Carlos,
C/Tulipán s/n 28933, M\'ostoles, Spain.
{\it E-mail address:} {\tt marta.latorre@urjc.es }}

\address{M. Magliocca: Dipartimento di Matematica, Sapienza Università di Roma, Piazzale Aldo Moro 5, 00185 Roma, Italy.
{\it E-mail address:} {\tt  magliocca@mat.uniroma1.it }}

\address{S. Segura de Le\'on: Departament d'An\`{a}lisi Matem\`atica,
Universitat de Val\`encia,
Dr. Moliner 50, 46100 Burjassot, Spain.
{\it E-mail address:}  {\tt sergio.segura@uv.es }}


\begin{document}

 \maketitle

\begin{center}
\end{center}

\begin{abstract}
We consider the homogeneous Dirichlet problem for an elliptic equation driven by a linear operator with discontinuous coefficients and having a subquadratic gradient term. This gradient term behaves as $g(u)|\nabla u|^q$, where $1<q<2$ and $g(s)$ is a continuous function. Data belong to $L^m(\Omega)$ with $1\le m <\frac N2$ as well as measure data instead of $L^1$-data, so that unbounded solutions are expected. Our aim is, given $1\le m<\frac N2$ and $1<q<2$, to find the suitable behaviour of $g$ close to infinity which leads to existence for our problem. We show that the presence of $g$ has a regularizing effect in the existence and summability of the solution. Moreover, our results adjust with continuity with known results when either $g(s)$ is constant or $q=2$.
\end{abstract}

\section{Introduction}

This paper is concerned to an elliptic problem, in an open bounded set $\Omega\subset\R^N$, whose model is:
\begin{equation}\label{problema}
\left\{
\begin{array}{ll}
-\Delta u = g(u)\,|\nabla u |^q + f(x) & \mbox{ in } \Omega\,,\\[2mm]
u=0 & \mbox{ on } \partial\Omega\,,
\end{array}\right.
\end{equation}
where
\begin{enumerate}
\item[H1.]  $g\,:\,\R\to(0,+\infty)$ is a continuous positive function;
\item[H2.] $1<q<2$;
\item[H3.] $f \in L^m(\Omega)$ such that $m\ge 1$. Eventually, we will also consider measures instead of $f\in L^1(\Omega)$.
\end{enumerate}
Our aim is, given $q$ and $m$, to find the suitable behaviour of $g$ close to infinity which leads to existence for problem \eqref{problema}.
We may measure the behaviour of $g$ through the exponent $\alpha$ such that the limit $\lim_{|s|\to\infty}|s|^\alpha g(s)$ is positive and finite. For the sake of simplicity, we will assume that $g$ is continuous and satisfies $g(s)\le \frac{\gamma}{|s|^\alpha}\;$, with $\;\alpha, \gamma >0$. So, we look for the possible exponents $\alpha$ for which we can obtain existence of solution to this problem.

Solutions to \eqref{problema} are considered in a weak sense (i.e., having finite energy) when $m\ge \frac{2N}{N+2}$. Nevertheless, this notion has no meaning when $m$ is closer to 1. In these cases the notion of weak solution must be replaced with the notion of entropy solution or that of renormalized one. Entropy solutions were introduced in \cite{B-V} for $L^1$--data and in \cite{BGO} for measure data which are absolutely continuous with respect to the capacity. On the other hand, renormalized solutions were handled in \cite{M,DMOP}. Since both notions are equivalent, in the present paper renormalized solution is the chosen notion and only the renormalized formulation will be used in what follows.

\subsection{Background}

Problems related to \eqref{problema} have been widely studied in recent years. Recall that, when $0\le q<1$ and data belong to the dual space $H^{-1}(\Omega)$, it can be solved applying the theory of pseudomonotone operators (see, for instance, \cite{LL}). Likewise, this theory also applies when $q=1$ and the norm of the source $f$ is small enough to get coerciveness. Without assuming any smallness condition, an existence result holds true as proved by Bottaro and Marina in \cite{BMa} for linear equations, and by Del Vecchio and Porzio in \cite{DVP} in a nonlinear framework.

The other growth limit, $q=2$, deserves some remarks. With additional hypotheses, equations having gradient terms with quadratic growth have been studied in a series of papers in the 80's, mainly by Boccardo, Murat and Puel. For gradient terms satisfying a sign condition, we refer to \cite{BMP1, BMP2, BBM}, while for existence of bounded solutions, to \cite{BMP3, BMP4}. The first attempt to study equations with a gradient term having natural growth (without the sign condition or an additional zero order term), was carried out by Ferone and Murat in \cite{FM} (see also \cite{FM1, FM2} for extensions).
They consider the case $g(s)$ constant and prove a sharp smallness condition on $f\in L^{\frac N2}(\Omega)$ which leads to an existence result.
More precisely, it was proved that if $\|f\|_{L^{\frac N2}(\Omega)}$ is small, then there exists a solution $u$ which also satisfies the further regularity $\pare{e^{\delta |u|}-1}\in H_0^1(\Omega)$, for $\delta$ less than a constant which only depends on $\|f\|_{L^{\frac N2}(\Omega)}$, the coerciveness of the principal term and the best constant in Sobolev's inequality. More general data in this quadratic growth were considered in \cite{S, P1}, under the assumption $g\in L^1(\R)$: It is studied existence for all $L^1$--data (in \cite{S}) and for all Radon measures (in \cite{P1}). The assumption $g\in L^1(\R)$ turns out to be optimal (see \cite[Proposition 5.1]{S}). The exhaustive analysis of the necessary growth condition on $g$ to obtain a solution for every datum $f\in L^m(\Omega)$, where $m>1$, was made by Porretta and Segura de Le\'on in \cite{PS}. The main result of \cite{PS}, when it is applied to an equation governed by a linear operator, states:
\begin{enumerate}
  \item[Q1.] Given any $f\in L^m(\Omega)$ with $m\ge\frac N2$: there exists a solution to problem \eqref{problema}, with $q=2$, under the assumption $\displaystyle\lim_{|s|\to\infty}g(s)=0$; this solution is bounded when $m>\frac N2$.
  \item[Q2.] Given any $f\in L^m(\Omega)$ with $\frac{2N}{N+2}\le m<\frac N2$: if  $ g(s)\leq \frac{\gamma}{|s|}$, with $\gamma<\frac{N(m-1)} {N-2m}$, then there exists a solution to problem \eqref{problema}, with $q=2$, which belongs to $H_0^1(\Omega)\cap L^{\frac{Nm}{N-2m}}(\Omega)$.
  \item[Q3.] Given any $f\in L^m(\Omega)$ with $1<m<\frac{2N}{N+2}$:  if $g(s) \le  \frac\gamma {|s|}$, with $\gamma<\frac{N(m-1)}{N-2m}$,  then there exists an entropy solution to problem \eqref{problema}, with $q=2$, which belongs to $W_0^{1,m^*}(\Omega)$.
\end{enumerate}
One of the main objectives of the present paper is to extend these results to the case $1<q<2$.
Some consequences of the quadratic case for our problem are:
\begin{enumerate}
  \item[C1.]  If  $\displaystyle  \lim_{s\to\pm\infty} g(s)=0$, then there exists a solution for every $f\in L^m(\Omega)$, with $m \ge\frac N2$ (see Proposition \ref{prop1} below).
  \item[C2.] If $g\in L^{\frac 2q}(\R)$, then there exists a solution for every $f\in L^1(\Omega)$ (see Remark \ref{obs2} below).
\end{enumerate}

In the subquadratic setting (i.e., $1<q<2$) and for $g(s)$ constant, the general theory was developed by Grenon, Murat and Porretta in \cite{GMP0} (for the range $1+\frac2N\le q<2$) and \cite{GMP} (with full generality). Their aim is to find the ``optimal" exponent $m$, which depends on $q$, such that there exists a solution for data in $L^m(\Omega)$ satisfying a smallness condition. Moreover,  Alvino, Ferone and Mercaldo showed in \cite{AFM} the sharp condition on datum $f$ which guarantee the existence of solution. It was proved in \cite{GMP0} that:
 \begin{enumerate}
   \item[S1.] If $1+\frac2N\le q<2$, $m\ge\frac{N(q-1)}{q}$  and $\|f\|_m$ is small enough, then there exists a solution $u\in H_0^1(\Omega)$ to problem \eqref{problema}, with $g(s)$ constant, which satisfies the further regularity $|u|^{\sigma}\in H_0^1(\Omega)$, with $\sigma=\frac{m(N-2)}{2(N-2m)}$.
        \end{enumerate}
       The extension studied in \cite{GMP} leads to:
 \begin{enumerate}
   \item[S2.] If $\frac N{N-1}<q<1+\frac2N$, $m\ge \frac{N(q-1)}{q}$ and $\|f\|_m$ is small enough, then there exists a renormalized solution to problem \eqref{problema}, with $g(s)$ constant, which satisfies the regularity $(1+|u|)^{\sigma-1}u\in H_0^1(\Omega)$, with $\sigma=\frac{m(N-2)}{2(N-2m)}$, and $|\nabla u|\in L^{N(q-1)}(\Omega)$.
   \item[S3.] If $q=\frac{N}{N-1}$ and $\|f\|_{L^m(\Omega)}$ is small enough for certain $m>1$, then there exists a renormalized solution to problem \eqref{problema}, with $g(s)$ constant, satisfying the regularity
   $(1+|u|)^{\sigma-1}u\in H_0^1(\Omega)$, with $\sigma=\frac{m(N-2)}{2(N-2m)}$, and $|\nabla u|\in L^{m^*}(\Omega)$.
      \item[S4.] If $1<q<\frac N{N-1}$, $m\ge1$  and $\|f\|_1$ is small enough, then there exists a renormalized solution to problem \eqref{problema}, with $g(s)$ constant, which satisfies $|u|\in M^{\frac N{N-2}}(\Omega)$ and $|\nabla u|\in M^{\frac N{N-1}}(\Omega)$. Here, $M^q(\Omega)$ stands for the Marcinkiewicz space (see Subsection \ref{ele1} below). Actually, sources more general than $L^1$--functions are handled, namely, finite Radon measures.
 \end{enumerate}

 The final picture looks as follows:

  
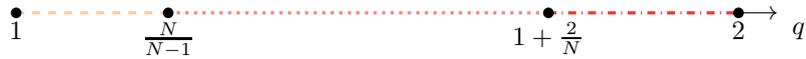
\begin{figure}[H]
\centering
\begin{tikzpicture}
\draw [->,thin] (9.5,0) -- (10,0);
\draw [dashed,very thick, color=orange!40!white] (0,0) to (2,0);
\draw [dotted,very thick, color=red!50!white] (2,0) to (7,0);
\draw [dash dot,very thick, color=red!80!white] (7,0) to (9.5,0);
\fill (0,0) circle (2pt) node[below] {$1$};
\fill (2,0) circle (2pt) node[below] {$\frac{N}{N-1}$};
\fill (7,0) circle (2pt) node[below] {$1+\frac{2}{N}$};
\fill (9.5,0) circle (2pt) node[below] {$2$};
\node [below left] at (10.5,0) { $q$};
\end{tikzpicture}
\end{figure}

The right zone indicates solutions of finite energy, the left zone shows the points $q$ where measure data can be considered while the central zone is where they obtain renormalized solutions with $L^m$-data.

We point out that, in \cite{GMP}, the authors obtain existence for every datum with a zero order term, which has a regularizing effect. In some sense, the singular term $\frac{|\N u|^q}{|u|^\al}$, where $\alpha>0$, induces a similar effect, so that we expect better estimates than those in the case $\alpha=0$ (see Remark \ref{remark1}). Note, nevertheless, that the term $\frac{|\N u|^q}{|u|^\al}$ behaves in a superlinear way with respect to the gradient power when $\alpha<q-1$ (see Subsection \ref{threshold}).

Finally, we remark that a problem similar to \eqref{problema} has recently been studied in \cite[Proposition $3.4$]{salva}.

\subsection{Our results}

As we have mentioned, given $q$ and $m$, our goal is to look for the best exponent $\alpha$ to get existence for our problem. The identity we find is
\begin{equation}\label{ident}
  \alpha=\frac{N(q-1)-mq}{N-2m}\quad\hbox{for } m>1,\,\, 1<q<2\,.
\end{equation}
This value of $\alpha$ is intuitively deduced in Subsection \ref{threshold}.  We remark that, when $\alpha=0$, it yields $m=\frac{N(q-1)}{q}$ recovering the threshold occurring in S1. and S2.  above. In general, $m=m(q,\alpha)$ is given by
\begin{equation}\label{ident-m}
m=\frac{N(q-1-\alpha)}{q-2\alpha} \,.
\end{equation}

According to the value of $m$ and the connection between $\alpha$ and $q$, there are two different types of results:
\begin{enumerate}
  \item As in Q2. or S1., if $m\ge \frac{2N}{N+2}$, then we get finite energy solutions. Otherwise, renormalized ones are obtained. We also note that there are points $(q,\alpha) \in [1,2]\times[0,1]$ satisfying
 \begin{equation*}
 \frac{N(q-1-\alpha)}{q-2\alpha}<1\,.
 \end{equation*}
 In this area, talking about Lebesgue spaces looses sense. This means that measure data are allowed.
  \item In full agreement with the above picture, we have to deal with three zones. If $q>1+\alpha$, we are within the superlinear framework. In this setting, we may only expect existence of solutions for sources satisfying a smallness condition. The limit case $q=1+\alpha$ corresponds to a linear gradient term in which we get existence when this term is small enough. Furthermore, the sublinear case $q<1+\alpha$ guaranties existence of solutions for all data and all gradient terms. The informal deduction of this classification will be shown in Subsection \ref{threshold}.
\end{enumerate}

We state our main results in Theorem \ref{main-p} and Theorem \ref{main-a} below, where all possible situations are considered.

Roughly speaking, we may  illustrate the relation between $\al,\,q,\,m$ in the following picture.

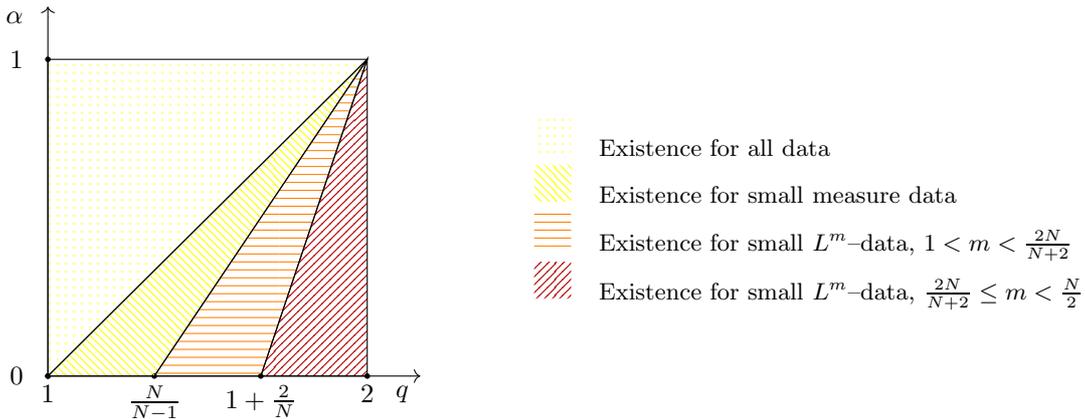
\begin{tabular}{ll}
\begin{minipage}{.43\textwidth}

\begin{figure}[H]

					\begin{tikzpicture}[scale=.7]

\fill (0,6.8) circle (0pt) node[left] {$\al\hspace*{2mm}$};
\fill (0,6) circle (1.5pt) node[left] {$1\hspace*{2mm}$};
\fill (0,0) circle (1.5pt) node[left] {$0\hspace*{2mm}$};

\fill (6,0) circle (1.5pt) node[below] {$2$};
\fill (0,0) circle (1.5pt) node[below] {$1$};
\fill (2,0) circle (1.5pt) node[below] {$\frac{N}{N-1}$};
\fill (4,0) circle (1.5pt) node[below] {$1+\frac{2}{N}$};
\fill (6.8,0) circle (0pt) node[below] {$q\hspace*{2mm}$};

\path[draw,pattern= dots
,pattern color=yellow!40!white] (0,0) -- (0,6) -- (6,6) -- cycle;

\path[draw,pattern=north west lines,pattern color=yellow] (0,0) -- (2,0) -- (6,6) -- cycle;

\path[draw,pattern=horizontal lines, pattern color=orange] (2,0) -- (4,0) -- (6,6) -- cycle;
			
\path[draw,pattern=north east lines,pattern color=red!70!black] (4,0) -- (6,0) -- (6,6) -- cycle;

	\draw[->] (0,0) -- (0,7) node  {};
	\draw[->] (0,0) -- (7,0) node {};
					\end{tikzpicture}

			\end{figure}
\end{minipage}
&
\begin{minipage}{.20\textwidth}
\begin{figure}[H]

\begin{center}
\small
					\parbox[b]{\textwidth}{
						\begin{tabular}{ll}
							\begin{tikzpicture}[scale=.5]
							\path[draw,white, pattern= dots,pattern color=yellow!40!white] (0,0) -- (0,1) -- (1,1) -- (1,0) -- cycle;
							\end{tikzpicture}
							& Existence for all data \\
							\begin{tikzpicture}[scale=.5]
							\path[draw,white, pattern=north west lines,pattern color=yellow] (0,0) -- (0,1) -- (1,1) -- (1,0) -- cycle;
							\end{tikzpicture}
							& Existence for small measure data \\
							\begin{tikzpicture}[scale=.5]
							\path[draw,white, pattern=horizontal lines,pattern color=orange] (0,0) -- (0,1) -- (1,1) -- (1,0) -- cycle;
							\end{tikzpicture}
							&  Existence for small $L^m$--data, $1<m<\frac{2N}{N+2}$ \\
							\begin{tikzpicture}[scale=.5]
							\path[draw,white, pattern=north east lines,pattern color=red!70!black] (0,0) -- (0,1) -- (1,1) -- (1,0) -- cycle;
							\end{tikzpicture}
							&  Existence for small $L^m$--data, $\frac{2N}{N+2}\le m<\frac{N}{2}$ \\
						\end{tabular}
					}
				\end{center}
			\end{figure}
\end{minipage}
\end{tabular}

We explicitly point out that, as $\alpha$ increases, the different zones drift to the right, so that the function $g$ induces a regularizing effect.   Moreover, the sublinear zone $0<\alpha<q-1$ appears, which entails existence for every datum. The bigger value of $\alpha$, the wider is this new zone.

This scheme adapts perfectly to what is expected since there is continuity with respect to known results. In fact, the $q$--axis coincides with the picture of results in \cite{GMP}, while the line $q=2$ depicts the results in \cite{PS}. Furthermore, the bound $\frac{N(m-1)}{N-2m}$, occurring in Q2. above, is the limit as $q\to2$ of the related bounds obtained for $q<2$ (see Remark \ref{obs7} below); a similar observation applies to Q3.

In order to achieve these bounds in the renormalized framework, we need to fine-tune our estimates as much as possible. So, we have to introduce a special way of applying known inequalities (see Lemma \ref{sob}). In this way we managed not to lose information when making our estimates.

To prove our results we use approximation techniques based on
\begin{itemize}
\item[1.] estimates of a suitable power of the solutions;
\item[2.] the strong convergence in $L^1$ of the gradient term.
\end{itemize}
Our estimates are obtained with variants of the method introduced by Grenon, Murat and Porretta in which certain powers of $G_k(u)=(|u|-k)^+\sg(u)$ are taken as test functions. The greatest difficulty arises when studying problem \eqref{problema} with a measure datum $\mu$. In this case one should take $T_j(G_k(u))$ as a test function getting
\begin{equation*}
\int_\Omega |\nabla T_j(G_k(u))|^2\,dx \le j \left[ \int_{\{ |u|>k\}} g(u) |\nabla u|^q\,dx +\|\mu\|_{\M}\right]\,.
\end{equation*}
An appeal to a lemma on Marcinkiewicz spaces (see \cite[Lemma 4.2]{BBM}) leads to
\begin{equation*}
\Big[|\nabla G_k(u)|\Big]_{\frac{N}{N-1}} \le C \left[ \int_{\{ |u|>k\}} g(u) |\nabla u|^q\,dx +\|\mu\|_{\M}\right]\,.
\end{equation*}
Having in mind that values $q>\frac{N}{N-1}$ are allowed, an estimate cannot be expected from this inequality. To overcome this trouble, we take $T_j(|G_k(u)|^\theta)\sg(u)$ as a test function, with the power $\theta$ close to $0$, and then, a suitable generalization of the above lemma (see Lemma \ref{marc} below) will be applied.

\subsection{Plan of the paper}

The plan of the paper is the following. Section 2 is devoted to introduce the assumptions and state the main results (see Theorem \ref{main-p} and Theorem \ref{main-a} below). We also include here our starting point (see Proposition \ref{prop1}), which is a simple consequence of the results of \cite{PS}. Section 3 deals with a priori estimates, while the convergence of approximate solutions is proved in Section 4. We point out that not only the superlinear case is seen, since we also deal with the linear case (in which existence for each data is achieved under a hypothesis of smallness on the gradient term) and even some sublinear cases that, as far as we know, have not already been handled. In this Section 4 the limit line $q=\frac{N+\alpha(N-2)}{N-1}$, which does not fit into the general scheme, is also studied. In Section 5 we end up analyzing what happens when data enjoy more summability than that strictly necessary to obtain existence.

\section{Assumptions and Statement of results}

\subsection{Notation}

Throughout this paper, $\Omega$ stands for a bounded open subset of $\R^N$, with $N>2$. The Lebesgue measure of a set $E\subset \Omega$ will be denoted by $|E|$. The symbols $L^s(\Omega)$ denote the usual Lebesgue spaces and $W^{1,s}_{0}(\Omega)$ the usual Sobolev spaces of measurable functions having weak gradient in $L^{s}(\Omega; \R^N)$ and zero trace on $\partial \Omega$. We will also use the notation $H_0^1(\Omega)$ instead of $W^{1,2}_{0}(\Omega)$. Let $1\le p<N$, in the sequel  $p^*=\frac{Np}{N-p}$, $p_*=\frac{Np}{Np-N+p}$ and $S_p$ stands for the constant in Sobolev inequality in $W_0^{1,p}(\Omega)$, that is,
 \[
 \left[\int_\Omega |u|^{p^*}dx\right]^{1/p^*}\le S_p \left[\int_\Omega |\nabla u|^{p}dx\right]^{1/p},\quad\hbox{for all }u\in W_0^{1,p}(\Omega)\,.
 \]
 We recall that this constant just depend on $N$ and $p$, and this dependence is continuous on $p$.
 On the other hand, $C^{PF}_p$ stands for the constant in the Poincar\'e--Friedrichs inequality in $W_0^{1,p}(\Omega)$, so that
 \[
 \left[\int_\Omega |u|^{p}dx\right]^{1/p}\le C^{PF}_p \left[\int_\Omega |\nabla u|^{p}dx\right]^{1/p},\quad\hbox{for all }u\in W_0^{1,p}(\Omega)\,.
 \]
 This constant depends on $\Omega$ and $p$.

 Two auxiliary real functions will be used throughout this paper. For every $k>0$, we define $T_k\>:\>\R\to\R$ and $G_k\>:\>\R\to\R$ as
 \begin{equation*}
 T_k(s)=\left\{\begin{array}{ll}
 s &\hbox{ if }|s|\le k\,,\\[3mm]
 k \frac s{|s|} &\hbox{ if }|s|> k\,;
 \end{array}\right.
 \end{equation*}
 \begin{equation*}
   G_k(s)=s-T_k(s)=(|s|-k)^+\sg(s)\,.
 \end{equation*}

\subsection{Assumptions}

We will deal with the problem
\begin{equation}\label{problem}
\left\{
\begin{array}{ll}
-\Div[A(x)\cdot\nabla u] = H(x, u, \nabla u) + f(x) & \mbox{ in } \Omega\,,\\[2mm]
u=0 & \mbox{ on } \partial\Omega\,,
\end{array}\right.
\end{equation}
and we assume the following statements.
\begin{enumerate}
  \item $A(x)$ is an $N\times N$ symmetric matrix which satisfies
  \begin{equation}\label{Hy1}
    \lambda |\xi|^2\le [A(x)\cdot\xi]\cdot\xi\le \Lambda |\xi|^2
  \end{equation}
  for almost all $x\in\Omega$ and $\xi\in\R^N$, and certain positive constants $\Lambda$ and $\lambda$.
  \item There exist a positive continuous function $g\>:\>\R\to(0,+\infty)$ and $1<q<2$ such that
   \begin{equation}\label{Hy2}
   |H(x, t, \xi)|\le g(t) |\xi |^q
   \end{equation}
  for almost all $x\in\Omega$, $t\in \R$ and $\xi\in\R^N$.
  \item The datum $f$ belongs to $L^m(\Omega)$, with $1\le m<\frac N2$. When $m=1$, instead of considering an $L^1$--function, we will choose a Radon measure (see problem \eqref{problem_measure}).
\end{enumerate}
As far as the function $g$ is concerned, we assume that there exist constants $\gamma, \alpha>0$   satisfying
  \begin{equation}\label{Hy2.5}
  g(t)\le \frac\gamma{|t|^\alpha}
  \end{equation}
  for all $t\in \R$.

  \begin{Remark}\rm
  \begin{enumerate}[(i)]
  \item Throughout this paper, the linearity of the principal part plays no role. We point out that our results also hold for equations driven by more general operators such as those of the Leray--Lions type with linear growth.
  \item We stress that condition \eqref{Hy2.5} have to be seen as a condition at infinity, that is, what is essential is
  $$
  \lim_{|s|\to\infty}g(s)|s|^\alpha=\ga\,.
  $$
  We are assuming \eqref{Hy2.5} for the sake of simplicity.
  \end{enumerate}
  \end{Remark}

    In what follows, we also consider the parameter
 \begin{equation}\label{mm}
 \sigma=\frac{m(N-2)}{N-2m}\,,\quad\hbox{i.e.}\quad m=\frac{N\sigma}{N+2\sigma-2}\,.
 \end{equation}
Observe that it implies
\begin{equation}\label{clave}
  (\sigma-1)m'=\frac\sigma 2 2^*\,.
\end{equation}

It is straightforward that
\begin{enumerate}
  \item $1\le\sigma<\infty$.
  \item $\sigma=1$ if and only if $m=1$.
  \item $1<\sigma<2$ if and only if $1<m<\frac{2N}{N+2}$.
  \item $\sigma\ge2$ if and only if $\frac{2N}{N+2}\le m <\frac N2$.
\end{enumerate}


\subsection{Notions of solution}

  According to the summability of the datum, we will find solutions to problem \eqref{problem} with finite energy or renormalized solutions. Definitions follow.

\begin{Definition}\label{fin-ener}
We will say that a function  $u\in H_0^{1}(\Omega)$  is a weak solution to problem \eqref{problem} if $H(x, u, \nabla u ) \in L^1 (\Omega )$ and
\begin{equation}\label{finenform}
 \int_{\Omega} [A(x)\cdot\nabla u] \cdot \nabla \varphi\, dx =
\int_{\Omega} H (x, u, \nabla u) \varphi\, dx+ \int_{\Omega} f(x) \varphi\, dx
\end{equation}
holds for all $\varphi \in H^{1}_0 (\Omega ) \cap L^{\infty} (\Omega )$.
 \end{Definition}

  We remark that, as a consequence of Sobolev's inequality, formulation \eqref{finenform} has sense only when $m\ge\frac{2N}{N+2}$, that is, $\sigma\ge2$.
  When $m<\frac{2N}{N+2}$ (so that $\sigma<2$), a different formulation must be required.
  The functional setting for the renormalized formulation lies on the space $\mathcal T_0^{1,2}(\Omega)$ of almost everywhere finite functions such that $T_k(u)\in H_0^{1}(\Omega)$ for all $k>0$.
  Functions in this space have a generalized gradient which (grosso modo) is defined by
  \[\nabla T_k(u)=(\nabla u)\chi_{\{|u|<k\}}\qquad\hbox{for all }k>0\,,\]
  (see \cite{B-V} or \cite{DMOP}).

    \begin{Definition}\label{renor}
    A function $u\>:\>\Omega\to\R$  is a
renormalized solution to problem \eqref{problem} having datum $f\in L^m(\Omega)$, with $1<m<\frac{2N}{N+2}$,  if it satisfies
\begin{enumerate}
  \item $u\in \mathcal T_0^{1,2}(\Omega)$;
  \item $\nabla u\in L^1(\Omega;\R^N)$;
  \item $H(x,u,\nabla u)\in L^1(\Omega)$;
\end{enumerate}
and
\begin{multline}\label{renor-f}
  \int_\Omega S'(u)\varphi [A(x)\cdot\nabla u]\cdot\nabla u \, dx+\int_\Omega S(u) [A(x)\cdot\nabla u]\cdot\nabla \varphi \, dx\\
  =  \int_\Omega H(x,u,\nabla u)S(u)\varphi\, dx+\int_\Omega f(x) S(u)\varphi\, dx
\end{multline}
holds for any Lipschitz function $S\>:\>\R\to\R$ with compact support
and for any $\varphi\in H^{1}(\Omega)\cap L^\infty(\Omega)$  such that $S(u)\varphi\in H_0^{1}(\Omega)$.
  \end{Definition}

\begin{Remark}\rm
Note that Definition \ref{renor} does not require any asymptotic condition on the energy term such as
\begin{equation}\label{as}
\lim_{n\to\infty}\frac{1}{n}\int_{\{n\le|u|\le 2n\}}|\N u|^2\,dx=0.
\end{equation}
Indeed, we will prove (in different steps) that
\begin{enumerate}
\item[R1.] If $1<\si<2$ (i.e. $1<m<\frac{2N}{N+2}$), then solutions enjoy a certain Sobolev regularity which implies \eqref{as}.
\item[R2.] If $\si=1$ (i.e. $m=1$), then condition \eqref{as} must be required to solutions.
\end{enumerate}
It is not difficult to check R1 if we assume condition
$\bra{(1+|u|)^{\frac{\sigma}{2}}-1}\in H_0^1(\Omega)$ with $1<\sigma<2$. Then
\[\frac{|\nabla u|^2}{(1+|u|)^{2-\sigma}}\in L^1(\Omega)\,.\]
Thus,
\begin{equation*}
  \frac1n\int_{\{n\le|u|\le 2n\}}|\nabla u|^2dx
  \le \frac1{n^{\sigma-1}}\left[2^{2-\sigma}\int_{\{n\le|u|\le 2n\}}\frac{|\nabla u|^2}{(1+|u|)^{2-\sigma}}dx\right]\le\frac C{n^{\sigma-1}}{\buildrel {n\to\infty} \over \longrightarrow}0
\end{equation*}
and condition \eqref{as} holds.

\end{Remark}

Up to now, we have taken $m=\frac{N(q-1-\al)}{q-2\al}$ with $m>1$. Nonetheless, this ratio can be strictly smaller than 1.
Then we take measure data and so consider problem
\begin{equation}\label{problem_measure}
\left\{
\begin{array}{ll}
-\Div[A(x)\cdot\nabla u] = H(x, u, \nabla u) + \mu & \mbox{ in } \Omega\,,\\[2mm]
u=0 & \mbox{ on } \partial\Omega\,,
\end{array}\right.
\end{equation}
being $\mu\in\M(\Omega)$ a bounded Radon measure, instead of problem \eqref{problem}.\\
As far as bounded Radon measures are concerned, we recall that every $\mu\in \M(\Omega)$ can be decomposed, in a unique way, as the sum $\mu=\mu_0+\mu_s$, where $\mu_0\in L^1(\Omega)+H^{-1}(\Omega)$ is the absolutely continuous (with respect to the capacity) part and $\mu_s$ is the singular one and it is concentrated on a set of null capacity. Further comments on measures data and the notion of capacity can be found in \cite[Section $2$]{BGO}, \cite[Section $2$]{DMOP}.

\begin{Definition}\label{measure}
    A function $u\>:\>\Omega\to\R$  is a
renormalized solution to problem \eqref{problem_measure} if it satisfies
\begin{enumerate}
  \item $u\in \mathcal T_0^{1,2}(\Omega)$;
  \item $\nabla u\in L^1(\Omega;\R^N)$;
  \item $H(x,u,\nabla u)\in L^1(\Omega)$;
\end{enumerate}
and
\begin{multline}\label{renor-mu}
  \int_\Omega S'(u)\varphi [A(x)\cdot\nabla u]\cdot\nabla u \, dx+\int_\Omega S(u) [A(x)\cdot\nabla u]\cdot\nabla \varphi \, dx\\
  =  \int_\Omega H(x,u,\nabla u)S(u)\varphi\, dx+\int_\Omega  S(u)\varphi\, d\mu_0
\end{multline}
holds for any Lipschitz function $S\>:\>\R\to\R$ with compact support
and for any $\varphi\in H^{1}(\Omega)\cap L^\infty(\Omega)$  such that $S(u)\varphi\in H_0^{1}(\Omega)$, and
\begin{align}
\lim_{n\to\infty}\frac{1}{n}\int_{\{n\le u\le 2n\}}\vp[A(x)\cdot\nabla u]\cdot\N u\,dx&=\integrale\vp\,d\mu_s^+,\label{as+}\\
\lim_{n\to\infty}\frac{1}{n}\int_{\{-2n\le u\le -n\}}\vp[A(x)\cdot\nabla u]\cdot\N u\,dx&=\integrale\vp\,d\mu_s^-,\label{as-}
\end{align}
for every $\vp\in C_b(\Omega)$, i.e. $\vp$ continuous and bounded in $\Omega$, and being $\mu_s^+$ and $\mu_s^-$ the positive and negative parts of $\mu_s$, respectively.
\end{Definition}

 \begin{Remark}\label{test}\rm
 Both in Definition \ref{renor} and Definition \ref{measure}, we will need to use test functions for which function $S$ has not compact support although $S'$ has. Most of them can be considered by a standard argument in the renormalized setting (see \cite{GMP0} for more details). This procedure consists of two steps, which we next apply to the main example $S(s)=T_k(s)$ in the case $f\in L^m(\Omega)$ with $m>1$.
 \begin{enumerate}[(1)]
   \item Take $S(s)=T_k(s)\vartheta_h(s)$, $\vp=1$ in \eqref{renor-f}  with
\begin{equation*}
\vartheta_h(s)=
\begin{cases}
\begin{array}{ll}
1 &\q |s|\le h\,,\\
\ds \frac{2h-|s|}{h} &\q h<|s|\le 2h\,,\\
0 & \q |s|>2h\,.
\end{array}
\end{cases}
\end{equation*}
Observe that $\vartheta_h(\cdot)$ is compactly supported and converges to $1$ as $h\to\infty$. In this way, \eqref{renor-f} becomes
\begin{multline*}
  \int_{\{|u|<k\}} \vartheta_h(u) [A(x)\cdot\nabla u]\cdot\nabla T_k(u) \, dx-\frac 1h\int_{\{h<|u|<2h\}} |T_k(u)|\, [A(x)\cdot\nabla u]\cdot\nabla u  \, dx\\
  =  \int_\Omega H(x,u,\nabla u)T_k(u)\vartheta_h(u)\, dx+\int_\Omega f T_k(u)\vartheta_h(u)\, dx\,.
\end{multline*}
   \item Check that letting $h$ go to $\infty$ is allowed in each term (in the second term, where $\vartheta_h'(u)$ appears, just apply condition \eqref{as}. It turns out that
       \begin{equation*}
        \int_{\Omega}  [A(x)\cdot\nabla u]\cdot\nabla T_k(u) \, dx
  =  \int_\Omega H(x,u,\nabla u)T_k(u)\, dx+\int_\Omega f T_k(u)\, dx\,.
       \end{equation*}
 \end{enumerate}
Analogous comments can be done when measure data are considered, taking the same test functions in \eqref{renor-mu} and using \eqref{as+}--\eqref{as-} instead of \eqref{as}, we get
       \begin{equation*}
        \int_{\Omega}  [A(x)\cdot\nabla u]\cdot\nabla T_k(u) \, dx
  =  \int_\Omega H(x,u,\nabla u)T_k(u)\, dx+\int_\Omega f T_k(u)\, d\mu_0+k\mu_s(\Omega)\,.
       \end{equation*}

 Throughout this paper, we will consider such general test functions without further comments.
 \end{Remark}

 \subsection{Main results}

 As we have seen in the Introduction, we get two different types of results:
 one in the superlinear setting and the other in the linear case. To justify this classification, we refer to the next Subsection.
 On the other hand, in both situations we should have in mind that, depending on the data, we will obtain finite energy solutions or renormalized ones.

 The results of this paper can be summarized in the following statements.

 \begin{Theorem}[Existence results in the superlinear case]\label{main-p}
 Using the above notation, assume that $\|f\|_{L^m(\Omega)}$ is small enough.
 \begin{enumerate}
   \item If $\frac{2N}{N+2}\le m<\frac N2$ and $\frac{N(q-1)-mq}{N-2m}\le \alpha<q-1$, then there exists a weak solution to problem \eqref{problem} satisfying $H(x,u,\nabla u)u\in L^1(\Omega)$, $H(x,u,\nabla u)\in L^{\frac 2q}(\Omega)$ and  the further regularity $|u|^{\frac{\sigma}{2}}\in H_0^1(\Omega)$.
   \item If $1<m<\frac{2N}{N+2}$ and $\frac{N(q-1)-mq}{N-2m}\le \alpha<q-1$, then there exists a renormalized solution to \eqref{problem} satisfying $(1+|u|)^{\frac{\sigma}{2}-1}u\in H_0^1(\Omega)$
 and $H(x,u,\nabla u)\in L^m(\Omega)$.
   \item If $m>1$ and $\alpha = \frac{N(q-1)-q}{N-2}$, then there exists a renormalized solution to \eqref{problem} satisfying $(1+|u|)^{\frac{\sigma}{2}-1}u\in H_0^1(\Omega)$.
 \end{enumerate}
Furthermore, assuming a source $\mu\in\M(\Omega)$ with $\|\mu\|_{\M}$ small enough, if $\frac{N(q-1)-q}{N-2}< \alpha<q-1$, then there exists a renormalized solution to \eqref{problem_measure}.
 \end{Theorem}

\medskip

 \begin{Theorem}[Existence results in the linear case]\label{main-a}
 With the same notation as above, assume that $\alpha=q-1$ and $\gamma$ is small enough.
 \begin{enumerate}
   \item If $\frac{2N}{N+2}\le m<\frac N2$, then there exists a weak solution to problem \eqref{problem}  which also satisfies $H(x,u,\nabla u)u\in L^1(\Omega)$, $H(x,u,\nabla u)\in L^{\frac 2q}(\Omega)$ and  the further regularity $|u|^{\sigma/2}\in H_0^1(\Omega)$.
   \item If $1<m<\frac{2N}{N+2}$, then there exists a renormalized solution to \eqref{problem} satisfying the regularity $(1+|u|)^{\frac{\sigma}{2}-1}u\in H_0^1(\Omega)$ and $H(x,u,\nabla u)\in L^m(\Omega)$.
 \end{enumerate}
 Furthermore, for every $\mu\in\M(\Omega)$ there exists a renormalized solution to \eqref{problem_measure}.
  \end{Theorem}

\subsection{Connection among parameters}\label{threshold}

The aim of this Subsection is to show the connection among all parameters of our problem which lead to existence of solution. The key argument is to find the best power $\sigma$ such that $u^\sigma$ can be taken as a test function.

We begin estimating the gradient term and seeing the connection between $q$ and $\alpha$. In order to simplify the incoming explanation, we consider the problem
\begin{equation}\label{pb}
\begin{cases}
\begin{split}
&\ds -\D u=\frac{F}{(1+u)^\al} & \mbox{ in }\Omega\,,\\
&\ds u=0 & \mbox{ on } \partial\Omega\,,
\end{split}
\end{cases}
\end{equation}
where $\al>0$ and $0<F \in L^\tau(\Omega)\mbox{ with } \tau<\frac{N}{2}$.
Note that $u>0$ by the classical maximum principle.\\
Basically, our aim is to prove a gradient estimate of the type
\[
\||\N u|^b\|_{L^1(\Omega)}\le c \|F\|_{L^\tau(\Omega)}^{\zeta}\,,
\]
for certain values $b\le 2$, $\zeta>0$. Once this step is concluded, we set $F=|\N u|^q$, i.e.
\[
\||\N u|^b\|_{L^1(\Omega)}\le c \||\N u|^b\|_{L^{\frac{\tau\,q}{b}}(\Omega)}^{\frac{q\zeta}{b}}\,,
\]
so we will deduce that
\begin{itemize}
\item we close the estimate choosing $\ds\frac{\tau\,q}{b}= 1$;
\item we are within the superlinear setting if and only if $\ds\frac{q\zeta}{b}>1$.
\end{itemize}
We take $(1+u)^{\si-1}-1$ as test function in \eqref{pb} for some $\sigma$. Then, defining $v=(1+u)^{\frac{\si}{2}}$, we obtain
\begin{equation}\label{dis1}
\integrale |\N v|^2\,dx\le c\integrale F \,v^{\frac{2}{\si}(\si-1-\al)}\,dx\,.
\end{equation}
Since  H\"older's inequality with $(\tau,\tau')$ implies
\begin{equation}\label{dis3}
\integrale F \,v^{\frac{2}{\si}(\si-1-\al)}\,dx\le \|F\|_{L^\tau(\Omega)}\|v\|_{L^{\frac{2}{\si}\tau'(\si-1-\al)}(\Omega)}^{\frac{2}{\si}(\si-1-\al)}\,,
\end{equation}
 we require
\begin{equation}\label{m(si,al)}
\tau=\tau(\si,\al)=\frac{N\si}{(N-2)(\al+1)+2\si},\q \t{i.e.}\q \frac{2}{\si}\tau'(\si-1-\al)=2^*\,.
\end{equation}
We estimate \eqref{dis3}  by applying Young's inequality with $\pare{\frac{2\tau'}{2^*},\frac{2\tau'}{2\tau'-2^*}}=\pare{\frac{2\tau'}{2^*},\frac{\tau(N-2)}{N-2\tau}}$. Then, invoking Sobolev's embedding too, we obtain
\[
\begin{split}
\integrale F \,v^{\frac{2}{\si}(\si-1-\al)}\,dx\le  \frac{1}{2}\|v\|_{H_0^1(\Omega)}^2+c\|F\|_{L^\tau(\Omega)}^{\frac{\tau(N-2)}{N-2\tau}}\,.
\end{split}
\]
Note that, having $\tau<\frac{N}{2}$,  this step makes sense. We gather \eqref{dis1}--\eqref{dis3} and deduce
\begin{equation}\label{disen}
\|v\|_{H_0^1(\Omega)}^2\le c\|F\|_{L^\tau(\Omega)}^{\frac{\tau(N-2)}{N-2\tau}}\,.
\end{equation}
Now, let $0<b<2$ and take into account $\||\N u|^b\|_{L^1(\Omega)}$. We omit the case $b=2$ since it can be dealt in the same way without passing to the change of variable $v=(1+u)^{\frac{\si}{2}}$. Then, by H\"older's inequality with $\pare{\frac{2}{b},\frac{2}{2-b}}$, we have that
\begin{equation}\label{dis2}
\begin{split}
\integrale |\N u|^b\,dx= \int_\Omega |\nabla v^\frac{2}{\sigma}|^b\,dx \le C\left(\integrale  |\N v|^2\,dx\right)^{\frac{b}{2}}\left( \integrale v^{\frac{2}{\si}\frac{b}{2-b}(2-\si)}\,dx\right)^{\frac{2-b}{2}}\,.
\end{split}
\end{equation}
We thus require
\begin{equation}\label{b(si)}
b=b(\si)=\frac{N\si}{N-2+\si},\q \t{i.e.}\q\frac{2}{\si}\frac{b}{2-b}(2-\si)=2^*\,.
\end{equation}
Note that $b<2$ if $\si<2$.
Thanks to Sobolev's embedding, the inequality in \eqref{dis2} becomes
\[
\||\N u|^b\|_{L^1(\Omega)}\le c \|v\|_{H_0^1(\Omega)}^b\|v\|_{L^{2^*}(\Omega)}^{b\frac{2-\si}{\si}}\le c S_2 \|v\|_{H_0^1(\Omega)}^{\frac{2b}{\si}}\,.
\]
Recalling \eqref{disen} too and taking $F=|\N u|^q$, we finally get
\[
\||\N u|^b\|_{L^1(\Omega)}
\le c\|F\|_{L^\tau(\Omega)}^{\frac{\tau(N-2)}{N-2\tau}\frac{b}{\si}}
\le c\| F\|_{L^{\tau q}(\Omega)}^{\frac{\tau(N-2)}{N-2\tau}\frac{bq}{\si}}
\le c\||\N u|^b\|_{L^{\tau\frac{q}{b}}(\Omega)}^{\frac{\tau(N-2)}{N-2\tau}\frac{q}{\si}}\,.
\]
Then
\begin{itemize}
\item we close the estimate taking $\ds 1= \frac{\tau q}{b}$. Thus, $\frac{b}{q}=\tau=\frac{N\sigma}{(N-2)(\alpha+1)-2\sigma}$ (by \eqref{m(si,al)}) from where we deduce taking into account \eqref{b(si)}
\begin{equation}\label{si}
\sigma = \frac{(N-2)(q-1-\alpha)}{2-q}
\end{equation}
and so
\begin{equation*}
\tau =\frac{b}{q}= \frac{N(q-1-\al)}{q(1-\al)}\,.
\end{equation*}
\item We are within the superlinear setting if and only if 
\[
\frac{\tau(N-2)}{N-2\tau}\frac{q}{\si}=\frac{q}{1+\al}>1\q \t{i.e.}\q q>1+\al\,.
\]
\end{itemize}
Since we want to keep us in a superlinear but still subquadratic setting, we will consider
\[
1+\al<q<2\,,
\]
which implies that $0<\al<1$. In other words, if $\al\ge 1$, then we are no longer in a superlinear gradient setting.
The linear one appears when $q=1+\alpha$, while we are in the sublinear setting when $q<1+\alpha$.


We now want to determine the relation between the gradient growth parameter $q$, the power growth parameter $\al$ and the data assumptions $f\in L^m(\Omega)$.
To this end, we focus on the source term $f$ and consider the simple problem
\begin{equation}\label{pb2}
\begin{cases}
\begin{split}
&\ds -\D u=f(x) & \mbox{ in }\Omega\,,\\
&\ds u=0 & \mbox{ on } \partial\Omega\,,
\end{split}
\end{cases}
\end{equation}
where $0<f \in L^m(\Omega)\mbox{ with } m<\frac{N}{2}$.
Now, if we  take again $\pare{(1+u)^{\si-1}-1}$, with $\si<2$ defined in \eqref{si}, as test function in \eqref{pb2} and reason as before, then we find that we need
\begin{equation*}
m(\si)=\frac{N\si}{N+2\si-2}\q \t{i.e.}\q m'(\si-1)= 2^*\frac{\si}{2}
\end{equation*}
Gathering this identity with \eqref{si}, we have
\[
m=\frac{N(q-1-\al)}{q-2\al},\q \t{i.e.}\q \al = \frac{N(q-1)-mq}{N-2m}\,.
\]
Therefore, we have informally deduced the need of conditions \eqref{mm}, \eqref{ident-m} and \eqref{ident}, respectively.

 \subsection{Our starting point}

 We begin with the following result which provides us of solution to approximating problems in Section 4. It follows from the results of \cite{PS}.

 \begin{Proposition}\label{prop1}
Consider two continuous functions $g_1, g_2 : \R\to(0,+\infty)$ such that $g_1\in L^{2/q}(\R)$ and $\ds\lim_{s\to\pm \infty} g_2(s) = 0$, and set $g=g_1+g_2$.
\begin{enumerate}
  \item If $f\in L^m(\Omega)$, with $m>N/2$, then there exists a solution to \eqref{problem} belonging to $H_0^1(\Omega)\cap L^\infty(\Omega)$.
  \item If $f\in L^{N/2}(\Omega)$, then there exists a solution to \eqref{problem} which belongs to $H_0^1(\Omega)\cap L^r(\Omega)$
  for all $1\le r<\infty$.
\end{enumerate}
\end{Proposition}

\begin{pf}
Note that the expression $G(u)=\int_0^ug_1(s)^{2/q}\, ds$ defines a real bounded function (due to $g_1\in L^{2/q}(\R)$). Now consider $\varphi$ a Lipschitz--continuous and increasing real function such that $\varphi(0)=0$.
Taking $e^{\frac{|G(u)|}\lambda}\varphi(u)$ as test function, $\lm$ as in \eqref{Hy1}, in \eqref{problem}, it follows from \eqref{Hy1}, \eqref{Hy2} and Young's inequality that
\begin{align*}
\lambda\int_\Omega e^{\frac{|G(u)|}\lambda}& \varphi^\prime(u) |\nabla u|^2\, dx+\int_\Omega |\varphi(u)| g_1(u)^{2/q}e^{\frac{|G(u)|}\lambda} |\nabla u|^2\, dx
\\
& \le\int_\Omega H(x, u,\nabla u)e^{\frac{|G(u)|}\lambda}\varphi(u)\, dx+\int_\Omega |f(x)|e^{\frac{|G(u)|}\lambda}\varphi(u)\, dx\\
& \le\int_\Omega  g_1(u) \varphi(u) e^{\frac{|G(u)|}\lambda} |\nabla u|^q\, dx+\int_\Omega  g_2(u) \varphi(u) e^{\frac{|G(u)|}\lambda} |\nabla u|^q\, dx+\int_\Omega |f(x)|e^{\frac{|G(u)|}\lambda}\varphi(u)\, dx
\\ &  \le \frac q2 \int_\Omega  g_1(u)^{2/q} |\varphi(u)| e^{\frac{|G(u)|}\lambda} |\nabla u|^2\, dx+\frac q2 \int_\Omega  g_2(u)^{2/q} |\varphi(u)| e^{\frac{|G(u)|}\lambda} |\nabla u|^2\, dx\\
&\qquad +(2-q)\int_\Omega |\varphi(u)| e^{\frac{|G(u)|}\lambda}\, dx+\int_\Omega |f(x)|e^{\frac{|G(u)|}\lambda}|\varphi(u)|\, dx\,.
\end{align*}
Simplifying, we deduce
\[
\lambda\int_\Omega e^{\frac{|G(u)|}\lambda}\varphi^\prime(u) |\nabla u|^2\, dx\le \frac q2\int_\Omega  g_2(u)^{2/q} |\varphi(u)| e^{\frac{|G(u)|}\lambda} |\nabla u|^2\, dx+\int_\Omega \left[|f(x)|+(2-q)\right]e^{\frac{|G(u)|}\lambda}|\varphi(u)|\, dx\,.
\]
Denoting $h(x)=C\left[|f(x)|+(2-q)\right]$, being $C$ an upper bound of $e^{\frac{|G(u)|}\lambda}$ (recall that $G$ is a bounded function), we have $h\in L^m(\Omega)$, here $m\ge N/2$. Then
\[
\lambda\int_\Omega \varphi^\prime(u) |\nabla u|^2\, dx\le  \frac q2 C\int_\Omega g_2(u)^{2/q} |\varphi(u)|  |\nabla u|^2\, dx+\int_\Omega h(x)|\varphi(u)|\, dx\,,
\]
for every Lipschitz--continuous and increasing real function $\varphi$ such that $\varphi(0)=0$.
Having in mind that $\lim_{s\to\pm \infty} g_2(s)^{2/q} = 0$, an appeal to the proofs of \cite[Theorem 2.1 and Theorem 2.2]{PS} shows that this estimate
 leads to existence for any $h\in L^m(\Omega)$ and consequently for every $f\in L^m(\Omega)$.
\end{pf}

\begin{Remark}\label{obs1}\rm
A straightforward consequence of the previous result is the existence of solutions for every $\alpha>0$ when $m\ge N/2$.
This is the reason to assuming $m<\frac N2$.
\end{Remark}

\begin{Remark}\label{obs2}\rm
The argument of the above result can also be applied to $L^1$--functions deducing existence of solution for any $f\in L^1(\Omega)$ when $g\in L^{\frac2q}(\R)$ (see \cite{S}, and \cite{P1} for its extension to measure data)
 and as consequence it is satisfied if $\;\alpha>q/2\;$. Nevertheless, this bound is not optimal since we will see that this fact holds for every $\;\alpha> q-1\;$ (note that $q-1<q/2$ if $1<q<2$). This gap will be studied in Theorem \ref{sub-lin1} below.
\end{Remark}

\section{A priori estimates}

Following \cite{GMP0, GMP}, the basic idea to get a priori estimates is to choose $|G_k(u)|^{\sigma-1}\sg(u)$ as test function in problem \eqref{problem}.
 Hence, we will consider three cases according to the value of the exponent $\sigma-1$. Roughly speaking, the easiest case is when $\frac{2N}{N+2} \le m < \frac{N}{2}$ (that is $\sigma \ge 2$) since then $|G_k(u)|^{\sigma-1}\sg(u)$ can be directly taken as test function. In the case $1<m<\frac{2N}{N+2}$ (that is $1<\sigma<2$), we have to replace it with $\bra{(\e+|G_k(u)|)^{\sigma-1}-\e^{\sigma-1}}\sg(u)$ since now the exponent does not define a Lipschitz--continuous function of $G_k(u)$. (Actually, we cannot take this function in the renormalized formulation, however we may follow the steps of Remark \ref{test} to approximate $\bra{(\e+|G_k(u)|)^{\sigma-1}-\e^{\sigma-1}}\sg(u)$ and lead to a similar estimate.) The last case is $m=1$ when the exponent vanishes and the test function must be bounded.

\subsection{Finite energy solutions}

\begin{Proposition}\label{teo1}
Let $f \in L^m(\Omega)$ with $\frac{2N}{N+2} \le m < \frac{N}{2}$ and let $ \alpha=\frac{N(q-1)-mq}{N-2m}$. Assume  \eqref{Hy1}, \eqref{Hy2}, \eqref{Hy2.5} and  that $u$ is a solution to problem \eqref{problem} in the sense of Definition \ref{fin-ener} such that $|u|^{\frac{\sigma}{2}}\in H_0^1(\Omega)$. (Observe that it yields $\sigma = \frac{(N-2)(q-1-\alpha)}{2-q}$ and $m=\frac{N(q-1-\alpha)}{q-2\alpha}$.)\\
Then, if $\|f\|_{L^m(\Omega)}$ is small enough, every such solution $u$  satisfies the following estimate:
\begin{equation*}
\int_\Omega |u|^{\sigma-2}|\nabla u|^2 \,dx \le M \,,
\end{equation*}
where $M$ is a positive constant which only depends on $N$, $q$, $m$, $\lm$, $\gamma$  and $\|f\|_{L^m(\Omega)}$.
\end{Proposition}

\begin{pf}
Let $k>0$. We start taking the test function $|G_k(u)|^{\sigma-1} \sg(u)$ in problem \eqref{problem}.
Then, by \eqref{Hy2}, we obtain
\begin{multline}\label{eq10}
\int_\Omega [A(x)\cdot\nabla u]\cdot \nabla(|G_k(u)|^{\sigma-1} \sg(u)) \\
 \le \int_\Omega g(u) \,|G_k(u)|^{\sigma-1} \sg(u) |\nabla u|^q \,dx + \int_\Omega  f\,|G_k(u)|^{\sigma-1} \sg(u) \,dx\,.
\end{multline}
On the left hand side, thanks to \eqref{Hy1}, we get
\begin{align*}
\int_\Omega [A(x)\cdot\nabla u]\cdot \nabla (|G_k(u)|^{\sigma-1} \sg(u)) &  \ge \lm \int_\Omega \nabla u \cdot \nabla \big[|G_k(u)|^{\sigma-1} \sg(u)\big] \,dx  \\&= \lm(\sigma-1)\int_\Omega |G_k(u)|^{\sigma-2}\,|\nabla G_k(u) |^2 \,dx
 = \lm\frac{4(\sigma-1)}{\sigma^2} \int_\Omega|\nabla |G_k(u)|^{\frac{\sigma}{2}}|^2\,dx \,.
\end{align*}
Recalling also \eqref{Hy2.5}, inequality \eqref{eq10} becomes
\begin{align}
\lm\frac{4(\sigma-1)}{\sigma^2} \int_\Omega|\nabla |G_k(u)|^{\frac{\sigma}{2}}|^2  & \le \int_\Omega g(u) \,|G_k(u)|^{\sigma-1} \sg(u) |\nabla u|^q \,dx + \int_\Omega  f\,|G_k(u)|^{\sigma-1} \sg(u) \,dx
 \nonumber \\ & \le \gamma \int_\Omega \frac{|G_k(u)|^{\sigma-1} }{|u|^\alpha}|\nabla u|^q\,dx + \int_\Omega |f| \,|G_k(u)|^{\sigma-1}\,dx =I_1+I_2\,.\label{ecu-2}
\end{align}
We start by performing some simple computations on the gradient term $I_1$.
\begin{align*}
I_1  = \gamma\,\int_\Omega\frac{|G_k(u)|^{\sigma-1} }{|u|^\alpha}|\nabla u|^q\, dx & \le \gamma\int_{\{|u|>k\}} \frac{|G_k(u)|^{\sigma-1} }{|G_k(u)|^\alpha} |\nabla u|^q \,dx
\\ & = \gamma\int_{\{|u|>k\}} |G_k(u)|^{\sigma-1-\alpha}\frac{|G_k(u)|^{(\frac{\sigma}{2}-1)q}}{|G_k(u)|^{(\frac{\sigma}{2}-1)q}} |\nabla G_k(u)|^q \,dx
\\ & =  \gamma\frac{2^q}{\sigma^q}\int_{\Omega} |G_k(u)|^{\sigma-1-\alpha-(\frac{\sigma}{2}-1)q}|\nabla |G_k(u)|^{\frac{\sigma}{2}}|^q \,dx\,,
\end{align*}
where we have used that $0<\alpha$ and that $|G_k(u)|\le |u|$ hold; we remark that no singularity appears since we are integrating on the set $\{|u|>k\}$. Then, applying H\"older's inequality, we deduce
\begin{equation}\label{ecu-3}
I_1 \le  \gamma \frac{2^q}{\sigma^q} \left(\int_{\Omega} |G_k(u)|^{[\sigma-1-\alpha-(\frac{\sigma}{2}-1)q]\frac{2}{2-q}} \,dx \right)^{\frac{2-q}{2}} \left( \int_{\Omega} |\nabla |G_k(u)|^{\frac{\sigma}{2}}|^2 \,dx \right)^{\frac{q}{2}} \,.
\end{equation}
\\[2mm]
Now, we will apply Sobolev's inequality. Indeed, since $\sigma = \frac{(N-2)(q-1-\alpha)}{2-q}$, the power of $G_k(u)$ in the first factor in \eqref{ecu-3} changes to
\begin{equation}\label{exp}
\left[\sigma-1-\alpha - \left(\frac{\sigma}{2}-1\right) q \right]\frac{2}{2-q} = \frac{\sigma}{2}\, 2^* \,.
\end{equation}
Therefore, it follows that
\begin{align*}
I_1 & \le \gamma \frac{2^q}{\sigma^q} \left(\int_{\Omega} |G_k(u)|^{\frac{\sigma}{2}2^*} \,dx \right)^{\frac{2-q}{2}} \||G_k(u)|^{\frac{\sigma}{2}} \|^{q}_{H_0^1(\Omega)}
\\ & \le \gamma \frac{2^q}{\sigma^q} \left(S_2^2\int_{\Omega} |\nabla |G_k(u)|^{\frac{\sigma}{2}}|^2 \,dx \right)^{\frac{2^*}{2}\frac{2-q}{2}} \||G_k(u)|^{\frac{\sigma}{2}} \|^{q}_{H_0^1(\Omega)}
\\ & =  \gamma \frac{2^q}{\sigma^q}S_2^{2^*\frac{2-q}{2}}  \||G_k(u)|^{\frac{\sigma}{2}} \|^{q+2^*\frac{2-q}{2}}_{H_0^1(\Omega)} = \gamma \frac{2^q}{\sigma^q}S_2^{2^*\frac{2-q}{2}}  \||G_k(u)|^{\frac{\sigma}{2}} \|^{2\frac{N-q}{N-2}}_{H_0^1(\Omega)} \,.
\end{align*}
\\[2mm]
On the other hand, we use H\"older's inequality on $I_2$ to get
\begin{equation*}
I_2 = \int_\Omega |f| \,|G_k(u)|^{\sigma-1} \,dx \le \left(\int_\Omega |f|^m \,dx\right)^{\frac{1}{m}} \left(\int_\Omega |G_k(u)|^{(\sigma-1)\,m'} \,dx \right)^\frac{1}{m'} \,.
\end{equation*}
Therefore, having in mind \eqref{clave},
\begin{align*}
I_2 \le \|f\|_{L^m(\Omega)}\left(\int_\Omega|G_k(u)|^{\frac{\sigma}{2}2^*}\,dx \right)^\frac{1}{m'} &  \le \|f\|_{L^m(\Omega)} \left(S_2^2\int_\Omega|\nabla|G_k(u)|^{\frac{\sigma}{2}}|^2\,dx \right)^{\frac{2^*}{2}\frac{1}{m'}}
\\&  =\|f\|_{L^m(\Omega)} S_2^{2\frac{\sigma-1}{\sigma}}\| |G_k(u)|^{\frac{\sigma}{2}}\|^{2\frac{\sigma-1}{\sigma}}_{H_0^1(\Omega)}
\,.
\end{align*}
Thus, inequality \eqref{ecu-2} becomes
\begin{equation*}
C_1\| |G_k(u)|^{\frac{\sigma}{2}}\|^2_{H_0^1(\Omega)}\le \gamma C_2  \||G_k(u)|^{\frac{\sigma}{2}} \|^{2\frac{N-q}{N-2}}_{H_0^1(\Omega)} + C_3\|f\|_{L^m(\Omega)} \| |G_k(u)|^{\frac{\sigma}{2}}\|^{2\frac{\sigma-1}{\sigma}} _{H_0^1(\Omega)}
\,,
\end{equation*}
for some positive constants $C_i$ only depending on $N$, $q$, $\lm$ and $m$ (this one through $\si$, by \eqref{mm}).
This is equivalent to
\begin{equation*}
C_1\| |G_k(u)|^{\frac{\sigma}{2}}\|^{\frac{2}{\sigma}}_{H_0^1(\Omega)}-\gamma C_2  \||G_k(u)|^{\frac{\sigma}{2}} \|^{2\frac{N-q}{N-2}-2\frac{\sigma-1}{\sigma}}_{H_0^1(\Omega)} \le C_3\|f\|_{L^m(\Omega)}
\,.
\end{equation*}
\\[2mm]
If we denote $Y_k=\| |G_k(u)|^{\frac{\sigma}{2}}\|^2_{H_0^1(\Omega)}$ and define the function
\begin{equation*}
F(y)=C_1y^{\frac{1}{\sigma}}-\gamma C_2  y^{\frac{N-q}{N-2}-\frac{\sigma-1}{\sigma}}\,,\quad y\ge 0\,,
\end{equation*}
we have obtained
\begin{equation}\label{FM}
F(Y_k) \le C_3 \|f\|_{L^m(\Omega)}\qq\forall k>0\,.
\end{equation}
Note that the continuous function $F(y)$ satisfies $F(0)=0$, $\lim_{y\to+\infty}F(y)=-\infty$, it is increasing until reaching certain $y^*$ and then it is decreasing,  so that it has a maximun $M^*$ at $y^*$, i.e., $M^* =F(y^*)= \max_{y} F(y)$. We explicitly remark that $M^*$ depends on $\gamma$ as well as on $q$, $N$, $\lambda$ and $m$.
Choosing constant
\begin{equation*}
K= \frac{M^*}{C_3} \,,
\end{equation*}
if we require $\|f\|_{L^m(\Omega)} <K$, then the equation $F(y)=C_3\|f\|_{L^m( \Omega)} < M^*$ has two roots:
\begin{equation*}
Y^- \;\mbox{ and }\; Y^+\,, \quad\mbox{ with }\quad Y^- < y^* < Y^+ \,.
\end{equation*}
It follows from $|u|^{\frac{\sigma}{2}}\in H_0^1(\Omega)$, that function $k\mapsto Y_k$ is continuous and goes to $0$ when $k \to \infty$. This fact implies $Y_k\le Y^-$ for all $k >0$, and so
\begin{equation*}
\frac{\sigma^2}{4}\int_\Omega |G_k(u)|^{\sigma-2}|\nabla G_k(u)|^2 \,dx =\int_\Omega|\nabla|G_k(u)|^{\frac{\sigma}{2}}|^2\,dx \le  Y^- \,,
\end{equation*}
for all $k>0$. Therefore, $\displaystyle\int_\Omega |u|^{\sigma-2}|\nabla u|^2 \,dx \le \frac{4}{\sigma^2} Y^- $.
\end{pf}

\begin{Remark}\label{remark1}\rm
We explicitly point out that our choice $m=\frac{N(q-1-\alpha)}{q-2\alpha}$ and our assumption $\alpha>0$ implies $\displaystyle m<\frac{N(q-1)}q$, so that the range for parameter $m$ is actually
$\frac{2N}{N+2} \le m< \frac{N(q-1)}q$. A simple consequence is that then $ q\ge \al\frac{N-2}{N}+1+\frac2N$, which, in particular, yields $\displaystyle q>1+\frac2N$.\\
If $\displaystyle q<\al\frac{N-2}{N}+ 1+\frac2N$, then we are allowed to consider data with a lower summability with respect to the case $\al=0$.
\end{Remark}

\begin{Remark}\label{remark2}\rm
The proof of Proposition \ref{teo1} for the case $\frac{N(q-1)-mq}{N-2m}<\alpha<q-1$ is similar to that of the limit case. The only differences begin in \eqref{exp} since now
\begin{equation*}
\beta:=\left[\sigma-1-\alpha - \left(\frac{\sigma}{2}-1\right) q \right]\frac{2}{2-q} < \frac{\sigma}{2}\, 2^* \,.
\end{equation*}
Therefore, H\"older's inequality must be applied once again in \eqref{ecu-3}:
\begin{align*}
I_1 & \le \gamma \frac{2^q}{\sigma^q} \left(\int_{\Omega} |G_k(u)|^{\beta} \,dx \right)^{\frac{2-q}{2}} \||G_k(u)|^{\frac{\sigma}{2}} \|^{q}_{H_0^1(\Omega)}
\\ & \le \gamma \frac{2^q}{\sigma^q}|\Omega|^{(2-q)(\frac1{2}-\frac\beta{2^*\sigma})} \left(\int_{\Omega} |G_k(u)|^{\frac{\sigma}{2}2^*} \,dx \right)^{\frac{\beta(2-q)}{2^*\sigma}} \||G_k(u)|^{\frac{\sigma}{2}} \|^{q}_{H_0^1(\Omega)}
\\ & \le \gamma \frac{2^q}{\sigma^q} |\Omega|^{(2-q)(\frac1{2}-\frac\beta{2^*\sigma})} \left(S_2^2\int_{\Omega} |\nabla |G_k(u)|^{\frac{\sigma}{2}}|^2 \,dx \right)^{\frac{\beta(2-q)}{2\sigma}} \||G_k(u)|^{\frac{\sigma}{2}} \|^{q}_{H_0^1(\Omega)}
\\ & =  \gamma \frac{2^q}{\sigma^q}S_2^{\frac{\beta(2-q)}{\sigma}} |\Omega|^{(2-q)(\frac1{2}-\frac\beta{2^*\sigma})} \||G_k(u)|^{\frac{\sigma}{2}} \|^{q+\frac{\beta(2-q)}{\sigma}}_{H_0^1(\Omega)} \,.
\end{align*}

From this point on, we can follow the same proof, we just note that now the constants also depend on  $\alpha$ and $|\Omega|$.
\end{Remark}

\begin{Remark}\rm
A relevant case occurs when $\alpha$ attains its limit value $\alpha=q-1$. Then $\beta=\sigma$ and so we have
\[
I_1\le \gamma\frac{2^q}{\sigma^q}S_2^{2-q}|\Omega|^{\frac{2-q}{N}}\|\,|G_k(u)|^{\frac\sigma 2}\|_{H_0^1(\Omega)}^2\,.
\]
A more accurate estimate follows from the Poincar\'e--Friedrichs inequality. It yields
\[
I_1\le \gamma\frac{2^q}{\sigma^q}(C^{PF}_2)^{2-q}\|\,|G_k(u)|^{\frac\sigma 2}\|_{H_0^1(\Omega)}^2\,.
\]
As a consequence, inequality \eqref{ecu-2} becomes
\[
\lm\frac{4(\sigma-1)}{\sigma^2} \|\,|G_k(u)|^{\frac\sigma 2}\|_{H_0^1(\Omega)}^2  \le \gamma\frac{2^q}{\sigma^q}(C^{PF}_2)^{2-q}\|\,|G_k(u)|^{\frac\sigma 2}\|_{H_0^1(\Omega)}^2 + \|f\|_{L^m(\Omega)} S_2^{2\frac{\sigma-1}{\sigma}}\|\, |G_k(u)|^{\frac{\sigma}{2}}\|^{2\frac{\sigma-1}{\sigma}}_{H_0^1(\Omega)}
\]
and an estimate for every $f\in L^m(\Omega)$ holds if
\[
\gamma\frac{2^q}{\sigma^q}(C^{PF}_2)^{2-q}<\lm\frac{4(\sigma-1)}{\sigma^2}\,.
\]
\end{Remark}

Hence, we have arrived at the following result.
\begin{Proposition}\label{prop2}
Let $f \in L^m(\Omega)$ with $\frac{2N}{N+2} \le m < \frac{N}{2}$ and let $ \alpha=q-1$. Assume  \eqref{Hy1}, \eqref{Hy2}, \eqref{Hy2.5} and that $u$ is a solution to problem \eqref{problem} in the sense of Definition \ref{fin-ener} such that $|u|^{\frac{\sigma}{2}}\in H_0^1(\Omega)$. \\
If
\[
\gamma<\lm\frac{2^{2-q}}{\sigma^{2-q}(C^{PF}_2)^{2-q}}(\sigma-1)\,,
\]
then such a solution $u$  satisfies the following estimate:
\begin{equation*}
\int_\Omega |u|^{\sigma-2}|\nabla u|^2 \,dx \le M \,,
\end{equation*}
for every $\|f\|_{L^m(\Omega)}$,
where $M$ is a positive constant which only depends on $N$, $q$, $m$, $\Omega$, $\lm$ and $\gamma$.
\end{Proposition}

\begin{Remark}\label{obs7}\rm
Noting that
$\displaystyle \sigma-1=\frac{N(m-1)}{N-2m}$,
it follows that the condition we have found in Proposition \ref{prop2} can be written as
\[
\gamma<\lm C^{2-q}\frac{N(m-1)}{N-2m}\q\t{for}\q C=\frac{2}{\si C_2^{PF}}\,.
\]
We point out that letting $q$ go to $2$, we obtain the same critical value appearing in \cite{PS}.
\end{Remark}

\subsection{Renormalized solutions with $L^m(\Omega)$ data}

In order to show that the parameters involved in all the cases are adjusted with continuity, the following result is necessary, it allows us to estimate sharply.

\begin{Lemma}\label{sob}
Let $v$ be a nonnegative function belonging to $W_0^{1,p}(\Omega)$   and consider $\varphi_\e(v)=(\e+v)^\nu$ for $\e>0$ and $0<\nu<1$. Then
\[
\Big(\int_{\Omega}\varphi_\e^{p^*}\, dx\Big)^{1/p^*}\le S_p\Big(A(\e)+\int_\Omega|\nabla \varphi_\e|^p\, dx\Big)^{1/p},
\]
where $A$ is a positive real function such that $\lim_{\e\to0}A(\e)=0$.
\end{Lemma}

\begin{pf} First note that $\varphi_\e\in W^{1,p}(\Omega)$ since $\varphi_\e$ is defined through a Lipschitz--continuous real function.
Now, extend $\varphi_\e$ to be $\e^\nu$ in $\R^N\backslash \Omega$.
We denote by $B_r$ the ball centered at the origin with radius $r$.
Fix $0<r<R$ in such a way that $\Omega\subset B_r$ and consider the cut--off function $\eta\in W^{1,\infty}(\R^N)$ with $0\le \eta\le1$ defined as
\[\begin{array}{ll}
  \eta(x)=1 & x\in B_r\,; \\[2mm]
  \eta(x)=0 & x\notin B_R\,; \\[2mm]
  |\nabla \eta(x)|=\frac 1{R-r} & x\in B_R\backslash B_r\,.
\end{array}\]
It follows from $\varphi_\e\eta\in W^{1,p}(\R^N)$, that
\[
\Big(\int_{\R^N}(\varphi_\e \eta)^{p^*}\, dx\Big)^{1/p^*}\le S_p \Big(\int_{\R^N}|\nabla (\varphi_\e \eta)|^{p}\, dx\Big)^{1/p}.
\]
As a consequence,
\begin{align*}
  \Big(\int_{\Omega}\varphi_\e^{p^*}\, dx\Big)^{1/p^*} & =\Big(\int_{\Omega}(\varphi_\e \eta)^{p^*}\, dx\Big)^{1/p^*}\le \Big(\int_{\R^N}(\varphi_\e \eta)^{p^*}\, dx\Big)^{1/p^*} \\
&   \le S_p \Big(\int_{\R^N}|\nabla (\varphi_\e \eta)|^{p}\, dx\Big)^{1/p} = S_p\Big(\int_{B_R\backslash B_r} \varphi_\e^p |\nabla \eta|^p\, dx+\int_\Omega\eta^p|\nabla \varphi_\e|^p\, dx\Big)^{1/p}\\
&   =S_p\Big(\frac{\e^{p\nu}}{(R-r)^p}|B_R\backslash B_r|+\int_\Omega |\nabla \varphi_\e|^p\, dx\Big)^{1/p},
\end{align*}
as desired.
\end{pf}

\begin{Remark}\label{remark2.5}\rm
We explicitly point out that a similar result holds for the Poincar\'e--Friedrichs inequality with $C_p^{PF}$ instead of $S_p$.
\end{Remark}

\begin{Proposition}\label{teo2}
Let $f \in L^m(\Omega)$ with $1 < m< \frac{2N}{N+2}$ and let $ \alpha=\frac{N(q-1)-mq}{N-2m}$. Assume \eqref{Hy1}, \eqref{Hy2}, \eqref{Hy2.5} and that $u$ is a renormalized solution to problem \eqref{problem} in the sense of Definition \ref{renor} such that $(1+|u|)^{\frac{\sigma}{2}-1}u \in H_0^1(\Omega)$. (Observe that then  $\sigma = \frac{(N-2)(q-1-\alpha)}{2-q}$ and $m=\frac{N(q-1-\alpha)}{q-2\alpha}$.)\\
Then, if $\|f\|_{L^m(\Omega)}$ is small enough, every such a solution $u$  satisfies the following estimate:
\begin{equation*}
\int_\Omega (1+|u|)^{\sigma-2}|\nabla u|^2 \,dx \le M \,,
\end{equation*}
where $M$ is a positive constant which only depends on $N$, $q$, $m$, $\gamma$, $\lm$ and $\|f\|_{L^m(\Omega)}$.
\end{Proposition}

\begin{pf}
Let $k>0$ and
fix $\e$ such that $0<\e <\min\{1, k\}$. We recall Remark \ref{test} and take the test function $S_{n,k}(u)\vp$, with
\begin{equation*}
S_{n,k}(u)=\int_0^{T_n(G_k(u))} (\e + |t|)^{\sigma-2} \,dt = \frac{1}{\sigma -1} \Big[ (\e + |T_n(G_k(u))| )^{\sigma-1}  -\e^{\sigma-1}\Big]\sg(u),\quad\hbox{ and }\quad \vp=1\,,
\end{equation*}
and so, by the growth condition \eqref{Hy2},
\begin{multline*}
  \int_\Omega [A(x)\cdot \nabla u] \cdot \nabla S_{n,k}(u) \,dx \le \int_\Omega g(u)\,|S_{n,k}(u)|\, |\nabla u|^q  \,dx+ \int_\Omega |f|\, |S_{n,k}(u)|\,dx\\
  = \int_{\{|u|>k\}} g(u)\,|S_{n,k}(u)|\, |\nabla u|^q  \,dx+ \int_\Omega |f|\, |S_{n,k}(u)|\,dx\,,
\end{multline*}
since $S_{n,k}(u)$ vanishes in the set $\{|u|\le k\}$.
On the left hand side we get
\begin{align*}
\int_\Omega [A(x)\cdot\nabla u] \cdot \nabla S_{n,k}(u) \,dx & = \int_\Omega (\e + |T_n(G_k(u))|)^{\sigma-2}[A(x)\cdot\nabla u] \cdot \nabla T_n(G_k(u)) \,dx
\\ & \ge \lm\frac{4}{\sigma^2} \int_\Omega  \big|\nabla (\e + |T_n(G_k(u))|)^{\frac{\sigma}{2}} \big|^2 \,dx \,
\end{align*}
by \eqref{Hy1}.\\
Therefore, we obtain
\begin{equation*}
\lm\frac{4}{\sigma^2} \int_\Omega  \big|\nabla (\e + |T_n(G_k(u))|)^{\frac{\sigma}{2}}\big|^2 \,dx
 \le \int_{\{|u|>k\}} g(u)|S_{n,k}(u)||\nabla u|^q\,dx + \int_\Omega |f|\,|S_{n,k}(u)|\,dx \,,
\end{equation*}
and letting $n \to \infty$ (which is licit thanks to the $\frac{\si}{2}$--power regularity), we have
\begin{equation}\label{eq4}
\lm \frac{4}{\si^2}\int_\Omega  |\nabla \varphi_\e^k(u)|^2 \,dx \le \int_{\{|u|>k\}} g(u)|S_k(u)||\nabla u|^q\,dx + \int_\Omega |f| \,|S_k(u)|\,dx  = I_1+ I_2\,,
\end{equation}
where we have denoted
\begin{equation*}
S_k(u) = \frac{1}{\sigma -1}  (\e + |G_k(u)| )^{\sigma-1} \,,
\end{equation*}
and
\begin{equation*}
\varphi^k_\e(u)= (\e + |G_k(u)| )^{\frac{\sigma}{2}} \,.
\end{equation*}
We start making some computations on $I_1$.
\begin{align*}
I_1 & \le \gamma \int_{\{|u|> k\}} \frac{|S_k(u)|}{|u|^\alpha}|\nabla u|^q\,dx \\
& =  \frac{\gamma}{\sigma -1} \int_{\{|u|>k\}} \frac{(\e + |G_k(u)| )^{\sigma-1}}{|u|^\alpha} |\nabla u|^q \,dx
\\
& \le \frac{\gamma}{\sigma -1}\int_{\{|u|>k\}}  (\e + |G_k(u)| )^{\sigma-1-\alpha}   |\nabla G_k(u)|^q \,dx
\end{align*}
owed to $\alpha>0$ and the fact that the inequality $\e + |G_k(u)| \le |u|$ holds in $\{ |u|>k \}$. Thus,
\begin{align*}
 I_1 & \le \frac{\gamma}{\sigma -1}\frac{2^q}{\sigma^q}\int_{\Omega}  (\e + |G_k(u)|)^{\sigma-1-\alpha} \frac{1}{(\e + |G_k(u)|)^{(\frac{\sigma}{2}-1)q}} |\nabla(\e + |G_k(u)|)^{\frac{\sigma}{2}}|^q \,dx
\\ & = \frac{\gamma}{\sigma -1} \frac{2^q}{\sigma^q} \int_{\Omega}  (\e + |G_k(u)|)^{\sigma-1-\alpha - (\frac{\sigma}{2}-1)q}|\nabla \varphi_\e^k(u)|^q \,dx \,.
\end{align*}
Moreover, applying H\"older's inequality we arrive at
\begin{equation}\label{eq5}
I_1 \le \frac{\gamma}{\sigma -1} \frac{2^q}{\sigma^q} \left( \int_{\Omega}  (\e + |G_k(u)|)^{\left[\sigma-1-\alpha - (\frac{\sigma}{2}-1)q \right]\frac{2}{2-q}}\,dx\right)^{\frac{2-q}{2}} \left( \int_{\Omega} |\nabla \varphi_\e^k(u)|^2 \,dx \right)^{\frac{q}{2}} \,.
\end{equation}
Since we have choosen $\sigma = \frac{(N-2)(q-1-\alpha)}{2-q}$, the power of $(\e+|G_k(u)|)$ in the first integrand is actually
\begin{equation*}
\left[\sigma-1-\alpha - \left(\frac{\sigma}{2}-1\right) q \right]\frac{2}{2-q} = \frac{\sigma}{2}\, 2^* \,.
\end{equation*}
Therefore, inequality \eqref{eq5} becomes
\[
I_1 \le \frac{\gamma}{\sigma -1}\frac{2^q}{\sigma^q} \left( \int_{\Omega}  |\varphi^k_\e(u)|^{2^*}\,dx\right)^{\frac{2-q}{2}} \|\nabla \varphi_\e^k(u)\|_{L^2(\Omega)}^{q}\,.
\]
Thanks to Lemma \ref{sob}, we may perform the following manipulations:
\begin{align*}
I_1 & \le \frac{\gamma}{\sigma -1}\frac{2^q}{\sigma^q}  S_2^{2^*\frac{2-q}{2}} \left(  \int_{\Omega}  |\nabla \varphi_\e^k(u)|^{2}\,dx + A(\e)\right)^{\frac{2^*}{2}\frac{2-q}{2}} \|\nabla \varphi_\e^k(u)\|_{L^2(\Omega)}^{q}
\\ & = \frac{\gamma}{\sigma -1}\frac{2^q}{\sigma^q}  S_2^{2^*\frac{2-q}{2}} \left(\|\nabla\varphi_\e^k(u)\|_{L^2(\Omega)}^2 +A(\e)\right)^{\frac{2^*}2\frac{2-q}2} \|\nabla\varphi_\e^k(u)\|_{L^2(\Omega)}^{q}
\\ & \le \frac{\gamma}{\sigma -1}\frac{2^q}{\sigma^q}  S_2^{2^*\frac{2-q}{2}} \left(\|\nabla\varphi_\e^k(u)\|_{L^2(\Omega)}^2 +A(\e)\right)^{\frac{2^*}2\frac{2-q}2+\frac q2} \,.
\end{align*}
\\[2mm]
On the other hand, we use H\"older's inequality in $I_2$ to get
\begin{equation*}
I_2 = \frac{1}{\sigma-1}\int_\Omega |f|\, (\e+|G_k(u)|)^{\sigma-1}  \,dx \le \frac{1}{\sigma-1} \left(\int_\Omega |f|^m \,dx\right)^{\frac{1}{m}} \left(\int_\Omega(\e+|G_k(u)|)^{(\sigma-1)\,m'}\,dx \right)^\frac{1}{m'} \,.
\end{equation*}
Therefore, on account of \eqref{clave} and applying Lemma \ref{sob} again,
\begin{align*}
I_2 & \le \frac{1}{\sigma-1}\|f\|_{L^m(\Omega)}\left(\int_\Omega(\e+|G_k(u)|)^{\frac{\sigma}{2}2^*}\,dx \right)^\frac{1}{m'}
 = \frac{1}{\sigma-1}\|f\|_{L^m(\Omega)}\left(\int_\Omega\varphi_\e^k(u)^{2^*}\,dx \right)^\frac{1}{m'}
\\ & \le \frac{1}{\sigma-1}\|f\|_{L^m(\Omega)}  S_2^{\frac{2^*}{m'}} \left(\int_\Omega |\nabla \varphi_\e^k(u)|^2\,dx +A(\e) \right)^{\frac{2^*}{2}\frac{1}{m'}}\,.
\end{align*}
Thus, inequality \eqref{eq4} becomes
\begin{align}\label{eq8}
\lm\frac 4{\sigma^2}\|\nabla\varphi_\e^k(u)\|_{L^2(\Omega)}^2 \le I_1+I_2 & \le \frac{\gamma}{\sigma -1}\frac{2^q}{\sigma^q}  S_2^{2^*\frac{2-q}{2}} \left(\|\nabla\varphi_\e^k(u)\|_{L^2(\Omega)}^2+A(\e)\right)^{\frac{2^*}2 \frac{2-q}{2}+\frac q2}
\\ \nonumber & + \frac{1}{\sigma-1}\|f\|_{L^m(\Omega)} S_2^{\frac{2^*}{m'}} \left(\|\nabla \varphi_\e^k(u)\|_{L^2(\Omega)}^2+A(\e)\right)^{\frac{2^*}{2m'}} \,.
\end{align}
If $k$ satisfies $G_k(u)=0$, then $\|\varphi_\e^k(u)\|_{H_0^1(\Omega)}=0$ and we are done. So, we will assume that $G_k(u)\ne 0$ and consequently $\lim_{\e\to0}\|\varphi_\e^k(u)\|_{H_0^1(\Omega)}\ne 0$. Then, we rearrange the terms of \eqref{eq8}, obtaining
\begin{multline}\label{eq9}
\lm\frac 4{\sigma^2}\|\nabla\varphi_\e^k(u)\|_{L^2(\Omega)}^{2-\frac{2^*}{m'}} \le \frac{\gamma}{\sigma -1}\frac{2^q}{\sigma^q}  S_2^{2^*\frac{2-q}{2}}\|\nabla\varphi_\e^k(u)\|_{L^2(\Omega)}^{2^*\frac{2-q}2+q-\frac{2^*}{m'}}\\
+B(\e)+
\frac{1}{\sigma-1}\|f\|_{L^m(\Omega)} S_2^{\frac{2^*}{m'}} \left(1+\frac{A(\e)}{\| \nabla\varphi_\e^k(u)\|_{L^2(\Omega)}^2}\right)^{\frac{2^*}{2m'}}\,,
\end{multline}
where
\[
B(\e)=\frac{\gamma}{\sigma -1}\frac{2^q}{\sigma^q}  S_2^{2^*\frac{2-q}{2}}\frac{\Big(\|\nabla\varphi_\e^k(u)\|_{L^2(\Omega)}^2+A(\e)\Big)^{\frac{2^*}2\frac{2-q}2+\frac q2}-\|\nabla\varphi_\e^k(u)\|_{L^2(\Omega)}^{{2^*}\frac{2-q}2+ q}}{\|\nabla\varphi_\e^k(u)\|_{L^2(\Omega)}^{\frac{2^*}{m'}}}
\]
defines a positive function which satisfies $\lim_{\e\to0}B(\e)=0$.
Denoting $Y_{k} = \||G_k(u)|^{\frac{\si}{2}}\|^2_{\ensp}$ and $Y_{k,\e} = \|\nabla\varphi_\e^k(u)\|^2_{L^2(\Omega)}$, inequality \eqref{eq9} changes to
\begin{equation*}
C_1Y_{k,\e}^{1-\frac{2^*}{2m'}} - \gamma C_2 Y_{k,\e}^{\frac{2^*}{2}\frac{2-q}{2}+\frac{q}{2}-\frac{2^*}{2m'}} \le B(\e) + C_3 \|f\|_{L^m(\Omega)} \left(1+\frac{A(\e)}{\| \nabla\varphi_\e^k(u)\|_{L^2(\Omega)}^2}\right)^{\frac{2^*}{2m'}} \,,
\end{equation*}
where each $C_i$ denotes a positive constant depending on $q$, $N$, $\lm$ and $m$.
\\[2mm]
Now, we consider again the function
\begin{equation*}
F(y)=C_1y^{1-\frac{2^*}{2m'}} -\gamma C_2 y^{\frac{2^*}{2}\frac{2-q}{2}+\frac{q}{2}-\frac{2^*}{2m'}} \,,
\end{equation*}
(note that $1-\frac{2^*}{2m'}=\frac{1}{\si}$ and $\frac{2^*}{2}\frac{2-q}{q}+\frac{q}{2}-\frac{2^*}{2m'}=\frac{N-q}{N-2}\frac{\si-1}{\si}$)
which has a maximun $M^*$ achieved at certain $y^*$, i.e., $M^* =F(y^*)= \max_{y} F(y)$.
Choosing constant
\begin{equation*}
K= \frac{M^*}{C_3} \,,
\end{equation*}
and requiring $\|f\|_{L^m(\Omega)}<K$, there exists $\e_0 \in(0,1)$ such that
\begin{equation*}
F(Y_{k, \e})\le M_{\e}=B(\e) + C_3 \|f\|_{L^m(\Omega)} \left(1+\frac{A(\e)}{\|\nabla \varphi_\e^k(u)\|_{L^2(\Omega)}^2}\right)^{\frac{2^*}{2m'}} < M^*
\end{equation*}
 for all $0<\e<\e_0$, and so the equation $F(y)=M_\e$ has two roots:
\begin{equation*}
Y_\e^- \;\mbox{ and }\; Y_\e^+\,, \quad\mbox{ with }\quad Y_\e^- < y^* < Y_\e^+ \,.
\end{equation*}
Observe that the continuity of $F$ leads to the continuity of the function $\e\mapsto Y_\e^-$.
\\[2mm]
From our hypothesis  $(1+|u|)^{\frac{\sigma}{2}-1}u\in H_0^1(\Omega)$, we have that function $k\mapsto Y_{k,\e}$ is continuous and goes to $0$ when $k \to \infty$.
Hence, $F(Y_{k, \e}) < M^*$ implies $Y_{k, \e}\le Y_\e^-$ for all $k >\e$ and, as a consequence,
\begin{align*}
\frac{\sigma^2}{4} \int_\Omega (1+|G_k(u)|)^{\sigma-2}|\nabla G_k(u)|^2 \,dx & \le \frac{\sigma^2}{4}\int_\Omega (\e +|G_k(u)|)^{\sigma-2}|\nabla G_k(u)|^2 \,dx
\\ & = \int_\Omega|\nabla(\e+|G_k(u)|)^{\frac{\sigma}{2}}|^2\,dx \le (Y_{\e}^-)^2 \,,
\end{align*}
for all $k>\e$. We point out that equation $F(y)=C_3\|f\|_{L^m(\Omega)}$ has two roots which will be denoted by $Y^-$ and $Y^+$, with $Y^-<Y_\e^- < y^* < Y_\e^+<Y^+$. Due to the continuity of function $F$ and since $\lim_{\e\to0}M_\e=C_3\|f\|_{L^m(\Omega)}$, it follows that $\lim_{\e\to0}Y_\e^-=Y^-$.
Hence,
\[ \int_\Omega(1+|G_k(u)|)^{\sigma-2}|\nabla G_k(u)|^{2} \,dx\le \frac{4}{\sigma^2} (Y^-)^2\]
 for all $k>0$ from where the desired estimate follows.
\end{pf}

\begin{Remark}\label{al-beta}\rm
As in Remark \ref{remark2}, we may extend the above result to the range $\frac{N(q-1)-mq}{N-2m}<\alpha<q-1$ with a constant depending also on $\alpha$ and $|\Omega|$.
\end{Remark}

In the same spirit than Proposition \ref{prop2}, a consequence of Proposition \ref{teo2} in the limit case $\alpha=q-1$ can be obtained. We also point out that when $q$ tends to $2$, it yields the same critical value found in \cite{PS}.
\begin{Proposition}\label{prop3}
Let $f \in L^m(\Omega)$ with $1<m<\frac{2N}{N+2}$ and let $ \alpha=q-1$. Assume \eqref{Hy1}, \eqref{Hy2}, \eqref{Hy2.5} and that $u$ is a renormalized solution to problem \eqref{problem} in the sense of Definition \ref{renor} such that $(1+|u|)^{\frac{\sigma}{2}-1}u\in H_0^1(\Omega)$. \\
If
\[
\gamma<\lm\frac{2^{2-q}}{\sigma^{2-q}(C^{PF}_2)^{2-q}}(\sigma-1)\,,
\]
then such solution $u$  satisfies the following estimate:
\begin{equation*}
\int_\Omega (1+|u|)^{\sigma-2}|\nabla u|^2 \,dx \le M \,,
\end{equation*}
where $M$ is a positive constant which only depends on $N$, $q$, $m$, $\Omega$, $\lm$, $\|f\|_{L^m(\Omega)}$ and $\gamma$.
\end{Proposition}
\begin{pf}
We may follow the same argument of the proof of Proposition \ref{teo2} until we reach the inequality \eqref{eq5}, which now is
\begin{equation*}
I_1 \le \frac{\gamma}{\sigma -1} \frac{2^q}{\sigma^q} \left( \int_{\Omega}  (\e + |G_k(u)|)^{\frac\sigma22}\,dx\right)^{\frac{2-q}{2}} \left( \int_{\Omega} |\nabla \varphi_\e^k(u)|^2 \,dx \right)^{\frac{q}{2}} \,.
\end{equation*}
Thus, the Poincar\'e--Friedrichs inequality yields
\begin{equation*}
I_1 \le \frac{\gamma}{\sigma -1} \frac{2^q}{\sigma^q}(C_2^{PF})^{2-q} \left( \int_{\Omega}  |\nabla \varphi_\e^k(u)|^2\,dx+A(\e)\right)^{\frac{2-q}{2}} \left( \int_{\Omega} |\nabla \varphi_\e^k(u)|^2 \,dx \right)^{\frac{q}{2}}
\end{equation*}
and so \eqref{eq8} becomes
\begin{multline*}
\lm\frac 4{\sigma^2}\|\nabla\varphi_\e^k(u)\|_{L^2(\Omega)}^2  \le \frac{\gamma}{\sigma -1}\frac{2^q}{\sigma^q}  (C_2^{PF})^{2-q} \left(\|\nabla\varphi_\e^k(u)\|_{L^2(\Omega)}^2+A(\e)\right)
\\ + \frac{1}{\sigma-1}\|f\|_{L^m(\Omega)} S_2^{\frac{2^*}{m'}} \left(\|\nabla \varphi_\e^k(u)\|_{L^2(\Omega)}^2+A(\e)\right)^{\frac{2^*}{2m'}} \,.
\end{multline*}
Therefore, the condition $\frac{\gamma}{\sigma -1}\frac{2^q}{\sigma^q}  (C_2^{PF})^{2-q}<\lm\frac 4{\sigma^2}$ implies a uniform estimate of $\nabla\varphi_\e^k(u)$ in $L^2(\Omega;\R^N)$.
We then infer the estimate
\begin{equation*}
\int_\Omega (1+|u|)^{\sigma-2}|\nabla u|^2 \,dx \le M \,.
\end{equation*}
\end{pf}

\subsection{Renormalized solutions with measure data}\label{ele1}

We recall here the definition and a few properties of \emph{Marcinkiewicz} spaces we are going to employ when dealing with the measure setting. \\
Let $0<\zeta<\infty$. Then, the Marcinkiewicz space $M^{\zeta}(\Omega)$ is defined as the set of measurable functions $u:\Omega\to \mathbb{R}$ such that
\[
[\,u\,]_\zeta=
\sup_{k>0}\left\{k^\zeta |\{x\in \Omega\>:\>|u(x)|>k \}|\right\}^{\frac{1}{\zeta}}<\infty\,.
\]
Furthermore, the following continuous embeddings hold
\[
L^\zeta(\Omega)\hookrightarrow M^{\zeta}(\Omega)\hookrightarrow L^{\zeta-\omega}(\Omega)
\]
for every $\omega>0$ such that $\zeta-\omega>1$. More precisely,
\begin{equation}\label{inm-mar}
\|f\|_{L^{\zeta-\omega}(\Omega)}\le \left(\frac\zeta\omega\right)^{\frac1{\zeta-\omega}}|\Omega|^{\frac\omega{\zeta(\zeta-\omega)}}[\,f\,]_\zeta
\end{equation}
holds for all $f\in M^\zeta(\Omega)$. We point out that the constant in the embedding depends on $\zeta$, $\omega$ and $|\Omega|$, and it blows up just when $\omega$ tends to $0$.

\begin{Lemma}\label{marc}
Let $\Omega\subset\R^N$ be a bounded open set. Let $1<p<N$ and $0<r<p$.
Consider $u\>:\>\Omega\to\R$ a measurable and a.e. finite function satisfying
\begin{equation*}
  T_\ell (u)\in W_0^{1, p}(\Omega)\,,\qquad\hbox{for all }\ell>0\,.
\end{equation*}
Assume that there exists $M>0$ such that
\begin{equation*}
\int_\Omega |\nabla T_\ell(u)|^p\le \ell^r M\,,\qquad\hbox{for all }\ell>0\,.
\end{equation*}

Then
\begin{gather}
 \label{lem3} [\,u\,]_{\frac{N}{N-p}(p-r)}\le C_1(N, p, r) M^{1/(p-r)}\,; \\
 \label{lem4} [\,|\nabla u|\,]_{\frac {N}{N-r}(p-r)}\le C_2(N, p, r) M^{1/(p-r)}\,.
\end{gather}
\end{Lemma}

\begin{pf} Applying Sobolev's inequality (and denoting by $S_p$ the Sobolev constant), we obtain
\begin{align*}
|\{|u|>\ell\}| & \le |\{|T_\ell(u)|\ge\ell\}|\le \int_\Omega\frac{|T_\ell(u)|^{p^*}}{\ell^{p^*}}\le\frac{S_p^{p^*}\left(\int_\Omega |\nabla T_\ell(u)|^p\right)^{p^*/p}}{\ell^{p^*}}
\\ & \le S_p^{p^*} M^{\frac{N}{N-p}} \ell^{-\frac{N}{N-p}(p-r)}\,,
\end{align*}
from where \eqref{lem3} follows.

To see \eqref{lem4}, perform the following manipulations:
\begin{align*}
|\{|\nabla u|>j\}| & \le |\{|\nabla T_\ell (u)|>j\}|+ |\{|u|>\ell\}|
\\ & \le \int_\Omega \frac{|\nabla T_\ell(u)|^p}{j^p}+S_p^{p^*} M^{\frac{N}{N-p}} \ell^{-\frac{N}{N-p}(p-r)}
\\ & \le \frac{\ell^r M}{j^p}+S_p^{p^*} M^{\frac{N}{N-p}} \ell^{-\frac{N}{N-p}(p-r)}\,.
\end{align*}
Since the minimum is obtained for
\[
\ell^*=\Big(\frac{N}{N-p}\frac{p-r}{r}S_p^{p^*}\Big)^{\frac{N-p}{p(N-r)}} j^{\frac{N-p}{N-r}}M^{\frac1{N-r}}\,,
\]
we deduce that
\begin{equation*}
 |\{|\nabla u|>j\}|\le C(N, p, r)  M^{\frac{N}{N-r}}j^{-\frac{N(p-r)}{N-r}},
\end{equation*}
wherewith \eqref{lem4} holds.
\end{pf}

\begin{Remark}\label{beatle1}
For further references, it is convenient to explicit the above constant $C_2(N, p, r)$. It is easy to check that
\[C_2(N, p, r)=C(N, p, r)^{\frac{N-r}{N(p-r)}}\]
and
\[C(N,p,r)=\left[\left(\frac{N(p-r)}{r(N-p)}\right)^{\frac{r(N-p)}{p(N-r)}}+\left(\frac{N(p-r)}{r(N-p)}\right)^{-\frac{N(p-r)}{p(N-r)}}\right]S_p^{\frac{Nr}{N-r}}\,.\]
We point out that
\begin{equation*}
\lim_{r\to0}C(N,p,r)=1\,.
\end{equation*}
\end{Remark}

\begin{Theorem}\label{teo3}
Assume \eqref{Hy1}, \eqref{Hy2}, \eqref{Hy2.5}.
Let $\mu \in \M(\Omega)$ and $\frac{N(q-1)-q}{N-2} < \alpha <q-1$.
If $\|\mu\|_{\M(\Omega)}$ is small enough, then every renormalized solution to problem \eqref{problem_measure} in the sense of Definition \ref{measure} satisfies
\begin{equation}\label{en_tk}
\int_\Omega |\nabla T_j(u)|^2 \,dx \le Mj \,,\quad\hbox{for all }j>0\,,
\end{equation}
and
\begin{equation*}
   \int_\Omega g(u) |\nabla u|^q \,dx \le C_0\,,
 \end{equation*}
where $M$ and $C_0$ are positive constants which only depends on $N$, $q$, $\gamma$, $\alpha$, $\lm$ and $\|\mu\|_{\M(\Omega)}$.
Moreover, for each $k>0$, the following estimate holds
\begin{equation*}
\int_{\{|u|>k\}} g(u) |\nabla u|^q \,dx \le C_k\,,
\end{equation*}
where $C_k$ is a constant which only depends on the above parameters of the problem and it satisfies
\begin{equation*}
\lim_{k\to \infty} C_k=0\,.
\end{equation*}
\end{Theorem}

\begin{pf}
Most of the proof consists of estimating the gradient term in $L^1(\Omega)$. \\
The case we are considering does not states any $\frac{\si}{2}$--class as in the previous results (see, however, Remark \ref{beatle2} below). We thus want to ``recreate" an analogous tool.\\
We choose
\begin{equation*}
  \theta=\frac{2qN-sq-2\alpha N+2\alpha s-N s}{s(N-q)}\,,
\end{equation*}
with  $s$ such that
\[
q<s<\frac{2N(q-\alpha)}{N+q-2\alpha},
\]
 in order to have $\theta>0$. Note that this condition is not restrictive since $q>\al+1>2\al$.\\
%
We now analyze the connection among all these parameters. Observe that $0<1-\frac{2\alpha}q$ holds because of the restriction $\alpha<\frac{q}{2}$, and $q<s$ implies
\[\theta <1-\frac{2\alpha}q\,.\]
On the other hand, it follows from $\displaystyle \frac{N(q-1)-q}{N-2}<\alpha$ that
\begin{equation}\label{alpha-theta}
  \displaystyle 0< \theta<\frac{N(2-s)}{s(N-2)}\,.
\end{equation}

Let $0 < j,k$ and let $\e>0$ satisfy $\e < k$. We start by taking $\psi(u)=T_j((\e+|G_k(u)|)^{\theta}-\e^{\theta})$ as test function in problem \eqref{problem_measure}.
Notice that $\psi(u)$ vanishes on the set $\{|u|\le k\}$.
Then, thanks to \eqref{Hy2},
\begin{equation}\label{eq10m}
\int_\Omega [A(x)\cdot\nabla u] \cdot \nabla \psi (u) \,dx \le  j\int_{\{|u|>k\}} g(u)|\nabla u|^q \,dx+j\|\mu\|_{\mathcal M_b(\Omega)} \,.
\end{equation}
On the left hand side we get
\begin{align*}
\int_\Omega [A(x)\cdot\nabla u] \cdot \nabla \psi(u) \,dx &=\int_\Omega [A(x)\cdot\nabla u] \cdot \nabla T_j((\e+|G_k(u)|)^{\theta}-\e^{\theta}) \,dx
\\ & =\theta \int_{\{ (\e+|G_k(u)|)^{\theta}-\e^\theta<j\}} (\e+|G_k(u)|)^{\theta-1} [A(x)\cdot\nabla u] \cdot \nabla G_k(u) \,dx
\\ & \ge \theta \int_{\{ (\e+|G_k(u)|)^{\theta}<j\}} (\e+|G_k(u)|)^{\theta-1} [A(x)\cdot\nabla u] \cdot \nabla G_k(u) \,dx
\\ &\ge \lm \theta \int_{\{\e+|G_k(u)|<j^{\frac{1}{\theta}}\}} (\e+|G_k(u)|)^{\theta-1}|\nabla G_k(u)|^2 \,dx
\\ & = \lm\theta \int_{\{ (\e+|G_k(u)|)^{\frac{\theta+1}{2}}<j^{\frac{\theta+1}{2\theta}} \}} \big[ (\e+|G_k(u)|)^{\frac{\theta-1}{2}}|\nabla G_k(u)|\big]^{2} \,dx
\\ & = \lm\theta \frac{4}{(\theta+1)^2} \int_\Omega |\nabla T_{j^{\frac{\theta+1}{2\theta}}}(\e+|G_k(u)|)^{\frac{\theta+1}{2}} |^{2} \,dx  \,,
\end{align*}
due to \eqref{Hy1}.\\
Thus, invoking \eqref{Hy2.5} too, \eqref{eq10m} becomes
\begin{align*}
\lm \frac{4\theta}{(\theta+1)^2} & \int_\Omega |\nabla T_{j^{\frac{\theta+1}{2\theta}}}(\e+|G_k(u)|)^{\frac{\theta+1}{2}} |^{2} \,dx \le  j\int_{\{|u|>k\}} g(u)  |\nabla u|^q \,dx + j\|\mu\|_{\mathcal M_b(\Omega)}
\\ & \le j \left\{\gamma \int_{\{|u|>k \}} \frac{1}{|u|^\alpha}|\nabla u|^q\,dx +\|\mu\|_{\M(\Omega)}\right\} =j (I+\|\mu\|_{\M(\Omega)})\,,
\end{align*}
where $\displaystyle I=\gamma \int_{\{|u|>k \}} \frac{1}{|u|^\alpha}|\nabla u|^q\,dx$. Moreover, writing $\ell=j^{\frac{\theta+1}{2\theta}}$ and $r=\frac{2\theta}{\theta+1}$ (i.e., $\theta=\frac{r}{2-r}$), we get
\begin{equation}\label{ecu-2m}
\lm \frac{4\theta}{(\theta+1)^2} \int_\Omega |\nabla T_\ell(\e+|G_k(u)|)^{\frac{\theta+1}{2}} |^{2} \,dx \le \ell^r (I +\|\mu\|_{\M(\Omega)})\,.
\end{equation}
We note that it follows from $\displaystyle \theta<\frac{N(2-s)}{s(N-2)}$ (see \eqref{alpha-theta}) that
$\displaystyle s<\frac{N(2-r)}{N-r}$.

We go on by performing some simple computations on the gradient term $I$.
\begin{align*}
I & = \gamma\,\int_{\{ |u|>k \}}\frac{1}{|u|^\alpha}|\nabla u|^q\, dx \le \gamma\int_{\{|u|>k\}} (\e+|G_k(u)|)^{-\alpha} |\nabla G_k(u)|^q \,dx
\\ & = \gamma\int_{\{|u|>k\}} (\e+|G_k(u)|)^{-\alpha}\frac{(\e+|G_k(u)|)^{\frac{\theta-1}{2}q}}{(\e+|G_k(u)|)^{\frac{\theta-1}{2}q}} |\nabla G_k(u)|^q \,dx
\\ & =  \gamma\frac{2^q}{(\theta+1)^q}\int_{\Omega} (\e+|G_k(u)|)^{-\alpha-\frac{\theta-1}{2} q}|\nabla (\e+ |G_k(u)|)^{\frac{\theta+1}{2}}|^q \,dx\,,
\end{align*}
where we have used that $0<\alpha$ and that $\e+|G_k(u)|\le |u|$ holds; we remark that no singularity appears since we are integrating on the set $\{|u|>k\}$. Then, applying H\"older's inequality with $\pare{\frac sq, \frac s{s-q}}$, we deduce
\begin{equation}\label{ecu-3m}
I \le  \gamma \frac{2^q}{(\theta+1)^q} \left(\int_{\Omega} (\e+|G_k(u)|)^{[-\alpha-\frac{\theta-1}{2}q]\frac{s}{s-q}} \,dx \right)^{\frac{s-q}{s}} \left( \int_{\Omega} |\nabla (\e+|G_k(u)|)^{\frac{\theta +1}{2}}|^s \,dx \right)^{\frac{q}{s}} \,.
\end{equation}
\\[2mm]
The next step is to estimate $I$ in terms of the function
\[\varphi_\e^k(u)=(\e+|G_k(u)|)^{\frac{\theta+1}2}\,.\]
To this end, we will apply Sobolev's inequality taking into account Lemma \ref{sob}. Indeed, the definition of $\theta$ implies that the power of $(\epsilon+|G_k(u)|)$ in the first integrand in \eqref{ecu-3m} changes to
\begin{equation}\label{expm}
\left[\frac{1-\theta}{2}q-\alpha\right]\frac{s}{s-q} = \frac{\theta+1}{2}\, s^* \,.
\end{equation}
Therefore, estimate \eqref{ecu-3m} becomes
\begin{align}\label{eq11}
I & \le \gamma \frac{2^q}{(\theta+1)^q} \left(\int_{\Omega} (\e+|G_k(u)|)^{\frac{\theta+1}{2}s^*} \,dx \right)^{\frac{s-q}{s}} \|\nabla \varphi_\e^k(u)\|^{q}_{L^s(\Omega)}
\\ \nonumber & \le \gamma \frac{2^{q}}{(\theta+1)^q} S_s^{s^*\frac{s-q}{s}} \bigg( \|\nabla \varphi_\e^k(u)\|_{L^s(\Omega)}^s+A(\e)\bigg)^{s^*\frac{s-q}{s^2}} \| \nabla \varphi_\e^k(u) \|^{q}_{L^s(\Omega)}\,.
\end{align}
Going back to inequality \eqref{ecu-2m}, we deduce
\begin{multline*}
 \int_\Omega|\nabla T_{\ell}( \varphi_\e^k(u))|^2\,dx\\
\le \ell^r\left[ \gamma \, C_1(r,s) \bigg( \|\nabla \varphi_\e^k(u)\|_{L^s(\Omega)}^s+A(\e)\bigg)^{s^*\frac{s-q}{s^2}} \| \nabla \varphi_\e^k(u) \|_{L^s(\Omega)}^q+C_2(r,s)\|\mu\|_{\M(\Omega)}\right]\,,
\end{multline*}
being
\begin{equation}\label{cucd}
C_1(r,s)=\frac{(\theta+1)^{2-q}}{2^{2-q}}\frac{S_s^{s^*\frac{s-q}{s}} }{\lm\theta} \q\t{and}\q C_2(r,s)=\frac{(\theta+1)^{2}}{4\theta\lm}\,.
\end{equation}
Note that $C_i(r,s)$, $i=1,\,2$, continuously depend on $r$ and $s$ (besides depending on $N$, $\lambda$ and $q$).
Using Lemma \ref{marc} it yields
\begin{multline*}
\Big[|\nabla \varphi_\e^k(u)|\Big]_{\frac{N(2-r)}{N-r}}\\
 \le c_0(r,s)\left[ \gamma \, {C}_1(r,s) \bigg( \|\nabla \varphi_\e^k(u)\|_{L^s(\Omega)}^s+A(\e)\bigg)^{s^*\frac{s-q}{s^2}} \|\nabla  \varphi_\e^k(u) \|_{L^s(\Omega)}^q+{C}_2(r,s)\|\mu\|_{\M(\Omega)} \right]^{\frac{1}{2-r}}\,,
\end{multline*}
for  some $c_0(r,s)$ continuously depending on $r$ and $s$, besides $N$.

Now recall we have taken $s<\frac{N(2-r)}{N-r}$, so that for each $(r,s)$ there exists a positive constant $C_0(r,s)$ continuously depending on $r$ and $s$, jointly with $N$ and $|\Omega|$,  such that
\[\|\nabla \varphi_\e^k(u)\|_{L^s(\Omega)} \le C_0(r,s)\Big[|\nabla \varphi_\e^k(u)|\Big]_{\frac{N(2-r)}{N-r}}\,.\]
Indeed, by \eqref{inm-mar}, we have
\[C_0(r,s)=
\pare{\frac{N(2-r)}{N(2-r)-s(N-r)}}^\frac{1}{s}\,|\Omega|^\frac{N(2-r)-s(N-r)}{sN(2-r)}\,.
\]
Note that $C_0(r,s)$ only blows up when $s\to\frac{N(2-r)}{N-r}$, which is impossible once $\alpha$ is fixed. Hence,
\begin{align*}
\|\nabla \varphi_\e^k(u)\|_{L^s(\Omega)} & \le C_0(r,s) \left[ | \nabla  \varphi_\e^k(u)| \right]_{\frac{N(2-r)}{N-r}}
\\ & \le c_0(r,s)\,C_0(r,s)\left[ \gamma \, {C}_1(r,s) \bigg( \|\nabla \varphi_\e^k(u)\|_{L^s(\Omega)}^s+A(\e)\bigg)^{s^*\frac{s-q}{s^2}} \|\nabla  \varphi_\e^k(u) \|_{L^s(\Omega)}^q+ {C}_2(r,s)\|\mu\|_{\M(\Omega)}\right]^{\frac{1}{2-r}}
\\ & \le  \gamma^{\frac1{2-r}} C_3(r,s)\bigg( \|\nabla \varphi_\e^k(u)\|_{L^s(\Omega)}^s+A(\e)\bigg)^{s^*\frac{s-q}{s^2}\frac1{2-r}} \|\nabla  \varphi_\e^k(u) \|_{L^s(\Omega)}^{\frac q{2-r}} +
C_4(r,s)\|\mu\|_{\M(\Omega)}^\frac{1}{2-r}
\\ & \le  \gamma^{\frac1{2-r}} C_3(r,s)\|\nabla  \varphi_\e^k(u) \|_{L^s(\Omega)}^{(s^*\frac{s-q}{s}+q)\frac{1}{2-r}} +C_4(r,s)\|\mu\|_{\M(\Omega)}^\frac{1}{2-r}+B(\e)\,,
\end{align*}
where
\[
B(\e)=\gamma^{\frac1{2-r}} C_3(r,s)\left[\bigg( \|\nabla \varphi_\e^k(u)\|_{L^s(\Omega)}^s+A(\e)\bigg)^{s^*\frac{s-q}{s^2}\frac1{2-r}}- \|\nabla \varphi_\e^k(u)\|_{L^s(\Omega)}^{s^*\frac{s-q}{s}\frac1{2-r}}\right]\|\nabla  \varphi_\e^k(u) \|_{L^s(\Omega)}^{\frac q{2-r}}
\]
satisfies $\lim_{\e\to0}B(\e)=0$.
Now, denoting $Y_{k,\e}=\|\nabla \varphi_\e^k(u)\|_{L^s(\Omega)}$, we have
\begin{equation*}
Y_{k,\e} - \gamma^{\frac1{2-r}} C_3(r,s) Y_{k,\e}^{(s^*\frac{s-q}{s}+q)\frac{1}{2-r}}\le C_4(r,s)\|\mu\|_{\M(\Omega)}^\frac{1}{2-r}+B(\e)\,.
\end{equation*}
We explicitly  note that the power of $Y_{k,\e}$ does not depend on either $r$ or $s$. Indeed, it is straightforward
\begin{equation*}
  \left(s^*\frac{s-q}{s}+q\right)\frac{1}{2-r}=\frac{s(N-q)}{N-s}\frac{1}{2-r}
\end{equation*}
and our definitions of $\theta$ and $r$ yield
\begin{equation*}
  \frac{r}{2-r}=\theta=\frac{2qN-sq-2\alpha N+2\alpha s-Ns}{s(N-q)}\q\t{and}\q r=\frac{2qN-sq-2\alpha N+2\alpha s-Ns}{(N-s)(q-2\al)}\,,
\end{equation*}
so that
\begin{equation*}
  \frac{s(N-q)}{N-s}\frac{1}{2-r}=\frac{2qN-sq-2\alpha N+2\alpha s-Ns}{r(N-s)}=q-\alpha\,.
\end{equation*}

We define the family of functions ($r>0$ and $q<s<\frac{2N(q-\alpha)}{N+q-2\alpha}$)
\begin{equation*}
F_{r,s}(y)=y-\gamma^{\frac1{2-r}} C_3(r,s)y^{q-\alpha}\,, \quad y>0 \,,
\end{equation*}
each one satisfying the same properties of that considered in the previous theorems and having a maximum $M_{r,s}^*$ at the point $y_{r,s}^*$.
We remark that we are not able to take limits when $r\to0\Leftrightarrow \theta\to0\Leftrightarrow s\to\frac{2N(q-\alpha)}{N+q-2\alpha}$ since in this case the constants $C_3$ and $C_4$ blow up.
Choose $\mu$ such that
\begin{equation*}
\|\mu\|_{\M(\Omega)}^\frac{1}{2-r}< \frac{M_{r,s}^*}{C_4(r,s)}\,,
\end{equation*}
for some $r>0$ and $q<s<\frac{2N(q-\alpha)}{N+q-2\alpha}$. From now on, we fix such parameters $r$ and $s$. Since $\lim_{\e\to0}B(\e)=0$ (note that $B(\e)$ also depends on $r$ and $s$), it follows that there exists $\e_0 >0$ such that
\begin{equation*}
 M_\e:=C_4(r,s)\|\mu\|_{\M(\Omega)}^\frac{1}{2-r}+B(\e) < M_{r,s}^*\,,
\end{equation*}
for all $0<\e<\e_0$.

Observe that the equation $F_{r,s}(y)=M_\e$ has two roots:
\begin{equation*}
(Y_{\e}^-)_{r,s} \;\mbox{ and }\; (Y_{\e}^+)_{r,s}\,, \quad\mbox{ with }\quad (Y_{\e}^-)_{r,s} < y_{r,s}^* < (Y_{\e}^+)_{r,s} \,,
\end{equation*}
and the continuity of $F_{r,s}$ leads to the continuity of the function $\e\mapsto (Y_{\e}^-)_{r,s}$.

Since the function $k\mapsto Y_{k,\e}$ is also continuous and goes to $0$ when $k \to \infty$, it follows from
$F(Y_{k,\e}) < M^*$ that $Y_{k,\e}\le (Y_{\e}^-)_{r,s}$ for all $k >\e$. As a consequence,
\begin{equation*}
\|\nabla \varphi_\e^k (u)\|_{L^s(\Omega)} \le (Y_{\e}^-)_{r,s} \,,
\end{equation*}
and so
\begin{equation*}
  \left(\frac{\theta+1}{2}\right)^s\int_{\{|u|>k\}}\frac{|\nabla u|^s}{|u|^{\frac{(1-\theta)s}{2}}}\, dx\le\left(\frac{\theta+1}{2}\right)^s\int_{\Omega}\frac{|\nabla G_k(u)|^s}{(\e+|G_k(u)|)^{\frac{(1-\theta)s}{2}}}\, dx\le (Y_{\e}^-)_{r,s}
\end{equation*}
for all $k>\e$ such that $\e<1$. Letting $\e\to0$, we obtain
\begin{equation}\label{lim_k}
  \int_{\{|u|>k\}}|\nabla |u|^{\frac{\theta+1}2}|^s\, dx\le (Y^-)_{r,s} \,,
\end{equation}
for all $k>0$. Here $(Y^-)_{r,s}$ stands for the smaller root of equation $F_{r,s}(y)=C_4(r,s)\|\mu\|_{\M(\Omega)}^\frac{1}{2-r}$. It is then straightforward that
\begin{equation}\label{ck}
  \lim_{k\to\infty}\int_{\{|u|>k\}}|\nabla |u|^{\frac{\theta+1}2}|^s\, dx=0\,.
\end{equation}

Taking into account \eqref{Hy2.5} and  \eqref{eq11}, it yields
\begin{align*}
\int_{\{|u|>k\}}g(u)|\nabla u|^q\, dx & \le \gamma\,\int_{\{ |u|>k \}}\frac{1}{|u|^\alpha}|\nabla u|^q\, dx
\\ & \le \gamma \frac{2^q}{(\theta+1)^q} S_s^{s^*\frac{s-q}{s}} \bigg(\|\nabla \varphi_\e^k(u)\|_{L^s(\Omega)}^s+A(\e)\bigg)^{s^*\frac{s-q}{s^2}} \|\nabla \varphi_\e^k(u) \|_{L^s(\Omega)}^{q}\,,
\end{align*}
and so, letting $\eps\to0$,
\begin{equation*}
  \int_{\{|u|>k\}}g(u)|\nabla u|^q\, dx  \le
  \gamma \frac{2^q}{(\theta+1)^q} S_s^{s^*\frac{s-q}{s}} \bigg(\int_{\{|u|>k\}}|\nabla |u|^{\frac{\theta+1}2}|^s\, dx\bigg)^{s^*\frac{s-q}{s^2}+\frac qs} =:C_k\,.
\end{equation*}
This is the key estimate we are looking for. Now it is enough to choose $C_0=\lim_{k\to0}C_k$ (on account of the estimate \eqref{lim_k}) and to realise that $\lim_{k\to\infty}C_k=0$ (see \eqref{ck}).

It just remains to check that \eqref{en_tk} holds. We take $T_j(u)$ as test function in problem \eqref{problem_measure}. It follows that
\[
\lm\int_\Omega|\nabla T_j(u)|^2\, dx\le j\gamma\,\int_{\Omega}g(u)|\nabla u|^q\, dx+j\|\mu\|_{\M(\Omega)}\le j\big(C_0 +\|\mu\|_{\M(\Omega)}\big)
\]
and we are done.
\end{pf}

\begin{Remark}\label{beatle2}
In contrast to what happens in Proposition \ref{teo1} and Proposition \ref{teo2}, in Theorem \ref{teo3} we do not provide any regularity condition on the solution. It is worth finding the regularity that results in our problem with measure datum. We point out that it is inadvisable to use \eqref{lim_k} because the values of $\theta$ and $s$ do not necessarily supply optimal regularity, besides they are not fully determined.

In problems with measure data, the regularity one obtains is
\[\int_\Omega\frac{|\nabla u|^2}{(1+|u|)^{1+\rho}}\, dx\le \frac M\rho\qquad \forall \rho>0\,,\]
here $M$ is the same constant stated in \eqref{en_tk}.
This inequality is easily deduced by taking \[S(u)=\left(1-\frac1{(1+|u|)^{\rho}}\right)\sg (u)\] as test function (in the sense of Remark \ref{test}).
\end{Remark}

We now turn to analyze the limit case $\alpha=q-1$.
\begin{Proposition}\label{teo4}
Let $\mu \in \M(\Omega)$ and let $ \alpha=q-1$. Assume \eqref{Hy1}, \eqref{Hy2}, \eqref{Hy2.5} and that $u$ is a renormalized solution to problem \eqref{problem_measure} in the sense of Definition \ref{renor-mu}.

If there exists $q<s<2$ satisfying
\[
\gamma<\lambda \frac{(2-s)^2}{Ns^{q-1}c_0(s)(C_s^{PF})^{s-q}|\Omega|^{\frac{2-s}{N}}}\,,
\]
where
\begin{equation}\label{bbeatle2}
c_0(s)=\left[\left(\frac{Ns}{(N-2)(2-s)}\right)^{\frac{(N-2)(2-s)}{2(N-2+s)}}+\left(\frac{Ns}{(N-2)(2-s)}\right)^{-\frac{Ns}{2(N-2+s)}}\right]S_2^{\frac{N(2-s)}{N-2+s}}\,,
\end{equation}
then every such solution $u$  satisfies the following estimate:
\begin{equation*}
\int_\Omega |\nabla T_k(u)|^2 \,dx \le M k\,,
\end{equation*}
and
\begin{equation*}
\int_{\{|u|>k\}} g(u)|\nabla u|^q \,dx \le C_k\,,
\end{equation*}
for every $k>0$,
where $M$ and $C_k$ are positive constants which only depend on $N$, $q$, $s$, $\lm$, $\Omega$, $\|\mu\|_{\M(\Omega)}$ and $\gamma$, and $\lim_{k\to\infty}C_k=0$.
\end{Proposition}
\begin{pf}
Since we follow a similar argument that that of the previous proof, we just sketch the proof. Take $s$ such that $q<s<2$ and if we define $r=2-s$ and $\theta=(2-s)/s$ (i.e. $r=2\theta/(\theta +1)$), then  $0<\theta<(s-q)/q$.
Fix $j>0$ and $0<\e<k$, and take again the test function $\psi(u)=T_j((\e+|G_k(u)|)^{\theta}-\e^{\theta})$ in problem \eqref{problem}.
Arguing as in the previous proof we also obtain \eqref{ecu-2m} and \eqref{ecu-3m}. Nevertheless, we now have
\begin{equation*}
\left[\frac{1-\theta}{2}q-q+1\right]\frac{s}{s-q} = \frac{\theta+1}{2}\, s
\end{equation*}
instead of \eqref{expm}, and so \eqref{ecu-3m} becomes
\begin{equation*}
 \int_{\{|u|>k\}}g(u)|\nabla u|^q dx \le \gamma \frac{2^q}{(\theta+1)^q} \left(\int_{\Omega} (\e+|G_k(u)|)^{\frac{\theta+1}{2}s} \,dx \right)^{\frac{s-q}{s}} \|\nabla \varphi_\e^k(u)\|^{q}_{L^s(\Omega)}\,.
\end{equation*}
Applying the Poincar\'e--Friedrichs inequality (recall Lemma \ref{sob} and Remark \ref{remark2.5}) we deduce that
\begin{equation}\label{ecu-4}
  \int_{\{|u|>k\}}g(u)|\nabla u|^q dx \le \gamma \frac{2^q}{(\theta+1)^q}(C_s^{PF})^{s-q}\left(\|\nabla \varphi_\e^k(u)\|^{s}_{L^s(\Omega)} +A(\e)\right)^{\frac{s-q}s} \|\nabla\varphi_\e^k(u)\|^{q}_{L^s(\Omega)}
\end{equation}
and then \eqref{ecu-2m} leads to
\begin{multline*}
\int_\Omega|\nabla T_{\ell}( \varphi_\e^k(u))|^2\,dx\\
\le \ell^r\left[ \gamma \frac{(\theta+1)^{2-q}}{\lm\theta 2^{2-q}}(C_s^{PF})^{s-q} \bigg( \|\nabla \varphi_\e^k(u)\|_{L^s(\Omega)}^s+A(\e)\bigg)^{s^*\frac{s-q}{s^2}} \| \nabla \varphi_\e^k(u) \|_{L^s(\Omega)}^q
+\frac{(\theta+1)^2}{4\lm\theta}\|\mu\|_{\M(\Omega)} \right]\,.
\end{multline*}
Therefore, recalling that $2-r=s$, Lemma \ref{marc} gives
\begin{multline*}
\Big[|\nabla \varphi_\e^k(u)|\Big]_{\frac{Ns}{N-2+s}}^s\\
 \le  c_0(s) \pare{\gamma \frac{(\theta+1)^{2-q}}{\lm\theta 2^{2-q}}(C_s^{PF})^{s-q} \bigg( \|\nabla \varphi_\e^k(u)\|_{L^s(\Omega)}^s+A(\e)\bigg)^{\frac{s-q}s} \|\nabla  \varphi_\e^k(u) \|_{L^s(\Omega)}^q+\frac{(\theta+1)^2}{4\lm\theta}\|\mu\|_{\M(\Omega)}}\,.
\end{multline*}
Taking on account Remark \ref{beatle1}, we deduce that $c_0(s)$ is given by \eqref{bbeatle2}.
Now observe that $s<\frac{Ns}{N-2+s}$ and so, having in mind \eqref{inm-mar}, there exists a constant $C_0(s)>0$ such that
\[\|\nabla \varphi_\e^k(u)\|_{L^s(\Omega)} \le C_0(s)\Big[|\nabla \varphi_\e^k(u)|\Big]_{\frac{Ns}{N-2+s}}\]
and $C_0(s)$ tends to $+\infty$ as $s\to2$; indeed,
\[C_0(s)=\left(\frac{N}{2-s}\right)^{\frac 1s}|\Omega|^{\frac{2-s}{Ns}}\,.\]
 Hence,
\begin{multline*}
\|\nabla \varphi_\e^k(u)\|_{L^s(\Omega)}^s \le c_0(s) C_0(s)^s \left[ | \nabla  \varphi_\e^k(u)| \right]_{\frac{Ns}{N-2+s}}^s\\
 \le c_0(s) C_0(s)^s\gamma \frac{(\theta+1)^{2-q}}{\lm\theta 2^{2-q}}(C_s^{PF})^{s-q} \bigg( \|\nabla \varphi_\e^k(u)\|_{L^s(\Omega)}^s+A(\e)\bigg)^{\frac{s-q}s} \|\nabla  \varphi_\e^k(u) \|_{L^s(\Omega)}^q\\
 +c_0(s) C_0(s)^s\frac{(\theta+1)^2}{4\lm\theta}\|\mu\|_{\M(\Omega)}\,.
\end{multline*}

Thus, recalling that $\theta=(2-s)/s$, we have obtained an estimate for
\[
\varphi(u)=(1+|u|)^{\frac{\theta+1}2}-1=(1+|u|)^{\frac{1}s}-1
\]
in $W_0^{1,s}(\Omega)$ if
\[
\gamma<\lm\frac{ (2-s)}{s^{q-1}(C_s^{PF})^{s-q}c_0(s)C_0(s)^{s}}=\lambda \frac{(2-s)^2}{Ns^{q-1}c_0(s)(C_s^{PF})^{s-q}|\Omega|^{\frac{2-s}{N}}}\,.
\]
Going back to \eqref{ecu-4}, letting $\e$ go to $0$ and denoting $\varphi^k(u)=(1+|G_k(u)|)^\frac{\theta+1}{2}$, we obtain
\begin{equation*}
  \int_{\{|u|>k\}}g(u)|\nabla u|^q dx \le \gamma \frac{2^q}{(\theta+1)^q}(C_s^{PF})^{s-q} \|\nabla\varphi^k(u)\|^{s}_{L^s(\Omega)}
  \le\gamma \frac{2^{q-s}}{(\theta+1)^{q-s}}(C_s^{PF})^{s-q}\int_{\{|u|>k\}}|\nabla |u|^{\frac{\theta+1}{2}}|^sdx\,,
\end{equation*}
wherewith $\int_{\{|u|>k\}}g(u)|\nabla u|^q dx \le C_k$ for certain $C_k$ such that $\lim_{k\to\infty}C_k=0$.

Finally, since the gradient term is bounded in $L^1(\Omega)$, it follows that the remaining estimate holds.
\end{pf}
\begin{Remark}\label{remark3}\rm
It is worth remarking what happens when $s\to2$ (i.e. $\theta\to 0$). Observe that it is not possible to choose any $\theta\in \pare{0,\frac{2-q}q}$, so that the above proof does not apply. Furthermore, since $\lim_{s\to2}c_0(s)=1$, it follows that
\[\lim_{s\to2}\frac{ (2-s)^2}{Ns^{q-1}c_0(s)(C_s^{PF})^{s-q}|\Omega|^{\frac{2-s}{N}}}=0\,.\]
Thus, no estimate is obtained for the equation
\[-\Delta u=\gamma\frac{|\nabla u|^2}{|u|}+f(x)\]
when $f\in L^1(\Omega)$ and $\gamma>0$.
This is in total agreement with \cite[Proposition 5.1]{PS}.
\end{Remark}

\subsection{The sublinear case with measure data.}
When $\alpha>q-1$, our problem lies in the sublinear setting.   Then we expect existence of a solution for each datum that is a finite Radon measure. To our knowledge, the range $q-1<\alpha\le \frac q2$ is not covered in previous papers, so that it will next be studied. We remark that the above proof can be extended to $\alpha$ satisfying $q-1<\alpha< \frac q2$ by choosing $\frac q{q-\alpha}<s<2$. Nevertheless, it does not work for $\alpha=\frac q2$. Hence, we will use very different test functions in the proof of the following result, which does not apply Lemma \ref{marc}.

\begin{Proposition}\label{teo5}
Let $\mu \in \M(\Omega)$. Assume \eqref{Hy1}, \eqref{Hy2}, \eqref{Hy2.5} and that $u$ is a renormalized solution to problem \eqref{problem_measure} in the sense of Definition \ref{measure}.

If $q-1<\alpha\le \frac q2$, then every such solution $u$  satisfies the following estimates:
\begin{equation*}
\int_\Omega |\nabla T_k(u)|^2 \,dx \le M k\,,
\end{equation*}
and
\begin{equation*}
\int_{\{|u|>k\}} g(u)|\nabla u|^q \,dx \le C_k\,,
\end{equation*}
for every $k>0$,
where $M$ and $C_k$ are positive constants only depending on the parameters of our problem, and $\lim_{k\to\infty}C_k=0$.
\end{Proposition}
\begin{pf}
We take
\begin{equation*}
  \varphi^k(u)=1-\frac1{(1+|G_k(u)|)^\theta}
\end{equation*}
as test function in \eqref{problem_measure}; here $k>1$ and $\theta$ is a positive parameter to be chosen. Then
\begin{align}\label{ecu-6}
  \lambda\theta\int_{\{|u|>k\}}\frac{|\nabla u|^2}{(1+|G_k(u)|)^{1+\theta}}\, dx & \le \int_{\{|u|>k\}}g(u)|\nabla u|^q\, dx +\|\mu\|_{\M(\Omega)}\\
\nonumber &  \le \ga \int_{\{|u|>k\}}\frac{|\nabla u|^q}{(1+|G_k(u)|)^{\alpha}}\, dx+\|\mu\|_{\M(\Omega)}\,.
\end{align}
In order to estimate the right hand side, we apply H\"older's inequality with exponents $\pare{\frac2q,\frac2{2-q}}$, getting
\begin{align*}
  I & =\int_{\{|u|>k\}}\frac{|\nabla u|^q}{(1+|G_k(u)|)^{\alpha}}\, dx=
  \int_{\{|u|>k\}}\frac{|\nabla u|^q}{[(1+|G_k(u)|)^{1+\theta}]^{\frac q2}}(1+|G_k(u)|)^{(1+\theta)\frac q2-\alpha}\, dx \\
&  \le \left(\int_{\{|u|>k\}}\frac{|\nabla u|^2}{(1+|G_k(u)|)^{1+\theta}}\, dx\right)^{\frac q2}
  \left(\int_{\{|u|>k\}}(1+|G_k(u)|)^{(1+\theta)\frac q{2-q}-\alpha\frac 2{2-q}}\, dx\right)^{\frac{2-q}2}\,.
\end{align*}
Now $\theta$ is chosen to satisfy $0<\theta<1+\alpha-q$, so that $(1+\theta)\frac q{2-q}-\alpha\frac 2{2-q}<1-\theta$ wherewith $\beta:=\frac{(1-\theta)(2-q)}{(1+\theta) q-2\alpha}>1$. H\"older's inequality, now with exponents $(\beta,\beta')$, and the Poincar\'e--Friedrichs inequality lead to
\begin{multline*}
  I\le \left(\int_{\Omega}\frac{|\nabla G_k(u)|^2}{(1+|G_k(u)|)^{1+\theta}}\, dx\right)^{\frac q2}
  \left(\int_{\Omega}\left[(1+|G_k(u)|)^{\frac{1-\theta}2}\right]^2\, dx\right)^{\frac{2-q}{2\beta}}|\Omega|^{\frac{2-q}{2\beta'}}\\
  \le |\Omega|^{\frac{2-q}{2\beta'}} (C_2^{PF})^{\frac{2-q}\beta} \left(\int_{\Omega}\frac{|\nabla G_k(u)|^2}{(1+|G_k(u)|)^{1+\theta}}\, dx\right)^{\frac q2}
  \left(\int_{\Omega}\left|\nabla(1+|G_k(u)|)^{\frac{1-\theta}2}\right|^2\, dx+ A(1)\right)^{\frac{2-q}{2\beta}}\\
  \le |\Omega|^{\frac{2-q}{2\beta'}} \left(\frac{1-\theta}2\right)^{\frac{2-q}\beta} (C_2^{PF})^{\frac{2-q}\beta} \left(\int_{\Omega}\frac{|\nabla G_k(u)|^2}{(1+|G_k(u)|)^{1+\theta}}\, dx\right)^{\frac q2}
  \left(\int_{\Omega} \frac{|\nabla G_k(u)|^2}{(1+|G_k(u)|)^{1+\theta}}\, dx+ A(1)\right)^{\frac{2-q}{2\beta}}\,.
\end{multline*}
Going back to \eqref{ecu-6} we obtain
\begin{multline*}
  \lambda\theta\int_{\{|u|>k\}}\frac{|\nabla u|^2}{(1+|G_k(u)|)^{1+\theta}}\, dx\le \|\mu\|_{\M(\Omega)}\\
  +\ga |\Omega|^{\frac{2-q}{2\beta'}} \left(\frac{1-\theta}2\right)^{\frac{2-q}\beta} (C_2^{PF})^{\frac{2-q}\beta} \left(\int_{\Omega}\frac{|\nabla G_k(u)|^2}{(1+|G_k(u)|)^{1+\theta}}\, dx\right)^{\frac q2}
  \left(\int_{\Omega} \frac{|\nabla G_k(u)|^2}{(1+|G_k(u)|)^{1+\theta}}\, dx+ A(1)\right)^{\frac{2-q}{2\beta}}
\end{multline*}
and it follows from $\frac q2+\frac{2-q}{2\beta}<1$ that there exists $M_1>0$ satisfying
\begin{equation*}
  \int_{\Omega}\frac{|\nabla G_k(u)|^2}{(1+|G_k(u)|)^{1+\theta}}\, dx\le M_1 \quad\hbox{ for all }k>1\,,
\end{equation*}
where $M_1$ only depends on $\lambda$, $q$, $\ga$, $\alpha$, $\Omega$ and $\|\mu\|_{\M(\Omega)}$. As a consequence of the above procedure, we also find $M_2>0$, depending on the same parameters, such that
\begin{equation*}
  \int_{\Omega}\frac{|\nabla G_k(u)|^q}{(1+|G_k(u)|)^{\alpha}}\, dx\le M_2 \quad\hbox{ for all }k>1\,.
\end{equation*}
A further estimate can be obtained observing that
\begin{align*}
  \int_{\{|u|>k\}}g(u)|\nabla u|^q\, dx
&  \le\ga  \left(\int_{\{|u|>k\}}\frac{|\nabla u|^2}{(1+|G_k(u)|)^{1+\theta}}\, dx\right)^{\frac q2}
  \left(\int_{\{|u|>k\}}(1+|G_k(u)|)^{(1+\theta)\frac q{2-q}-\alpha\frac2{2-q}}\, dx\right)^{\frac{2-q}2}\\
&  \le \ga \left(\int_{\Omega}\frac{|\nabla G_k(u)|^2}{(1+|G_k(u)|)^{1+\theta}}\, dx\right)^{\frac q2}
  \left(\int_{\Omega}\left[(1+|G_k(u)|)^{\frac{1-\theta}2}\right]^2\, dx\right)^{\frac{2-q}{2\beta}}|\{|u|>k\}|^{\frac{2-q}{2\beta'}}\\
&  \le M_3|\{|u|>k\}|^{\frac{2-q}{2\beta'}}=C_k\,,
\end{align*}
that holds, at least, for every $k>1$.

Taking $T_k(u)$, for some $k>1$ fixed, as test function in \eqref{problem_measure}, we derive
\begin{align*}
  \lambda\int_\Omega |\nabla T_k(u)|^2\, dx & \le
  \ga\int_{\{|u|<k\}}|\nabla u|^q|u|^{1-\alpha}\, dx+ k\int_{\{|u|\ge k\}}g(u) |\nabla u|^q\, dx+ k\|\mu\|_{\M(\Omega)}\\
&  \le \ga k^{1-\alpha}\int_{\Omega}|\nabla T_k(u)|^q\, dx+kC_k+ k\|\mu\|_{\M(\Omega)}\\
&  \le \ga k^{1-\frac q2}\int_{\Omega}|\nabla T_k(u)|^q\, dx+kC_k + k\|\mu\|_{\M(\Omega)}\,.
\end{align*}
Then Young's inequality implies an estimate of $T_k(u)$ in $H_0^1(\Omega)$ for every $k>1$ (and so for every $k>0$). We finally deduce an estimate of the gradient term in $L^1(\Omega)$. Indeed, fix $k>1$, denote $\bar g_k=\sup_{|s|\le k}|g(s)|$ and split the gradient term as follows
\begin{align*}
 \int_\Omega g(u)|\nabla u|^q\, dx &= \int_{\{|u|\le k\}} g(u)|\nabla u|^q\, dx+\int_{\{|u|>k\}} g(u)|\nabla u|^q\, dx \\
&  \le \bar g_k\int_\Omega |\nabla T_k(u)|^q\ dx+C_k\,.
\end{align*}

Once the gradient term is estimated in $L^1(\Omega)$, the remaining estimate is easy.
\end{pf}

\section{Compactness and convergence results}
Let us consider the approximating problems
\begin{equation}\label{appr}
\begin{cases}
\begin{array}{ll}
-\dive[A(x)\cdot\N u_n] = H(x,u_n,\N u_n) +f_n(x)& \mbox{ in } \Omega\,,\\[2mm]
u_n=0 & \mbox{ on } \partial\Omega\,,
\end{array}
\end{cases}
\end{equation}
with $f_n=T_n(f)$.
Proposition \ref{prop1} implies that there exists at least a solution $u_n\in L^\infty(\Omega)\cap\ensp$ such that
\begin{equation}\label{solappr}
\integrale [A(x)\cdot\nabla u_n]\cdot\N \vp\,dx=\integrale H(x,u_n,\N u_n) \vp\,dx+\integrale  T_n(f(x))\varphi\, dx \qq\forall\vp\in L^\infty(\Omega)\cap\ensp .
\end{equation}
We also handle measure data in Subsection \ref{ss4} but considering different approximating problems for \eqref{problem_measure}.

This Section is devoted to check that, up to subsequences, $\{u_n\}_n$ converges to a solution to problem \eqref{problem}.

\subsection{The case of solutions with finite energy}

\begin{Proposition}\label{convdebH2}
Let $f \in L^m(\Omega)$ with $\frac{2N}{N+2}\le m < \frac{N}{2}$, $ \alpha=\frac{N(q-1)-mq}{N-2m}$, $\sigma=\frac{(N-2)m}{N-2m}$ and $\{u_n\}_n$ be a sequence of solutions of \eqref{appr}. Assume also \eqref{Hy1}, \eqref{Hy2} and \eqref{Hy2.5}. Then
\begin{equation}\label{Hsigma}
\{|u_n|^{\frac{\si}{2}}\}_n\q\t{is uniformly bounded in}\q \ensp\,,
\end{equation}
\begin{equation}\label{HH}
\{u_n\}_n\q\t{is uniformly bounded in}\q \ensp \,,
\end{equation}
and
\begin{equation}\label{L1est}
\{H(x,u_n,\nabla u_n)u_n \}_n\q\t{is uniformly bounded in}\q L^1(\Omega)\,.
\end{equation}
Furthermore, up to subsequences, there exists a function $u$ such that
\begin{equation}\label{convdebH}
u_n\rightharpoonup u\q\t{in }\ensp\,,
\end{equation}
and
\begin{equation}\label{convun}
u_n\to u \qq \t{a.e. in }\Omega\,.
\end{equation}
\end{Proposition}
\begin{pf}
We apply Proposition \ref{teo1} to \eqref{appr} and deduce \eqref{Hsigma}.

Taking $\vp=u_n$ in \eqref{solappr} and recalling \eqref{Hy1}, \eqref{Hy2} and \eqref{Hy2.5}, we get
\[
\lm\integrale |\N u_n|^2\,dx\le\gamma \integrale |\N u_n|^q |u_n|^{1-\al}\,dx+\integrale|f||u_n|\,dx.
\]
We apply H\"older's inequality with indices $\left(\frac{2}{q},\frac{2}{2-q}\right)$ and $(m,m')$, respectively, on the integrals on the right hand side obtaining
\begin{align*}
\lm\|u_n\|_{H_0^1(\Omega)}^2 & \le \gamma\|u_n\|_{H_0^1(\Omega)}^q
\left(\integrale |u_n|^{\frac{2(1-\al)}{2-q}}\,dx\right)^{\frac{2-q}{2}}
+\|f\|_{L^{m}(\Omega)}\|u_n\|_{L^{m'}(\Omega)}\\
&\le \gamma\|u_n\|_{H_0^1(\Omega)}^q
\left(\integrale |u_n|^{\frac{2(1-\al)}{2-q}}\,dx\right)^{\frac{2-q}{2}}
+|\Omega|^{\frac 1{m'}-\frac1{2^*}}\|f\|_{L^{m}(\Omega)}\|u_n\|_{L^{2^*}(\Omega)}\\
&\le \gamma\|u_n\|_{H_0^1(\Omega)}^q
\left(\integrale |u_n|^{\frac{2(1-\al)}{2-q}}\,dx\right)^{\frac{2-q}{2}}
+S_2|\Omega|^{\frac 1{m'}-\frac1{2^*}}\|f\|_{L^{m}(\Omega)}\|u_n\|_{H_0^1(\Omega)}
\end{align*}
thanks, also, to Lebesgue spaces inclusion (indeed $m'\le 2^*$ by assumptions) and to Sobolev's embedding.
Then, twice applications of Young's inequality with $\left(\frac{2}{q},\frac{2}{2-q}\right)$ and $(2,2)$ yield to
\begin{equation}\label{boundH}
\lm\left(\frac{2-q-\eps}{2}\right)\|u_n\|_{H_0^1(\Omega)}^2 \le
\frac{2-q}{2}\,\frac{\gamma^{\frac{2}{2-q}}}{\lm^{\frac{q}{2-q}}}\integrale |u_n|^{\frac{2(1-\al)}{2-q}}\,dx
+\frac{S_2^2}{2\eps\lm}|\Omega|^{\frac2{m'}-\frac2{2^*}}\|f\|^2_{L^{m}(\Omega)}\,.
\end{equation}
We now take advantage of the power regularity in \eqref{Hsigma}, namely: $\{u_n\}_n$ is bounded in $L^{\frac{2^*\sigma}2}(\Omega)$. Observe that
$$
1-\alpha=\frac{(N-m)(2-q)}{N-2m}\le \frac{Nm(2-q)}{2(N-2m)}\,,
$$
owed to $m\ge\frac{2N}{N+2}$. Hence,
$$
\frac{2(1-\alpha)}{2-q}\le \frac{Nm}{N-2m}=\frac{2^*\sigma}2\,,
$$
so that
 the right hand side of \eqref{boundH} is uniformly bounded in $n$ and this means that \eqref{HH} holds. In particular we deduce \eqref{convdebH} and \eqref{convun} too.

 As far as the $L^1$--bound \eqref{L1est} is concerned, it is also a consequence of the inequality
 $$
 \gamma \integrale |\N u_n|^q |u_n|^{1-\al}\,dx\le \gamma\|u_n\|_{H_0^1(\Omega)}^q
\left(\integrale |u_n|^{\frac{2(1-\al)}{2-q}}\,dx\right)^{\frac{2-q}{2}}
 $$
 which we already know being bounded.
\end{pf}

\begin{Proposition}\label{propequiintdge2}
Let $f \in L^m(\Omega)$ with $\frac{2N}{N+2}\le m < \frac{N}{2}$, $ \alpha=\frac{N(q-1)-mq}{N-2m}$ and $\{u_n\}_n$ be a sequence of solutions of \eqref{appr}. Assume also \eqref{Hy1}, \eqref{Hy2} and \eqref{Hy2.5}. Then,
\begin{equation}\label{equi}
H(x,u_n,\N u_n)\q\t{ is equi--integrable in}\q L^1(\Omega)\,.
\end{equation}
Moreover,
up to subsequences, we have
\begin{equation}
\N u_n\to \N u\qq \t{a.e. in }\Omega.\label{convNun}
\end{equation}
In particular
\begin{equation}
H(x,u_n,\N u_n)\to H(x,u,\N u)\qq \t{in }L^1(\Omega).\label{rhs}
\end{equation}
 Furthermore,
\begin{equation}\label{L1est1}
H(x,u,\nabla u)u \in L^1(\Omega)
\end{equation}
and
\begin{equation*}
H(x,u,\nabla u) \in L^{2/q}(\Omega).
\end{equation*}
\end{Proposition}
\begin{pf}
We begin by showing that the sequence $\{H(x,u_n,\nabla u_n)\}_n$ is uniformly bounded in $L^\frac{2}{q}(\Omega)$. In fact,
\begin{equation*}
\int_\Omega |H(x,u_n,\nabla u_n)|^\frac{2}{q}\,dx \le \bar{g_1}^\frac{2}{q} \int_{\{|u_n|<1\}} |\nabla u_n|^2\,dx +\gamma^\frac{2}{q}\int_{\{|u_n|\ge1\}} |\nabla u_n|^2 \,dx \le (\bar{g_1}^\frac{2}{q}+\gamma^\frac{2}{q})\; C\,,
\end{equation*}
where $\ds \bar{g_1}=\max_{\{|s|\le 1\}} |g(s)|<\infty$, being $g(\cdot)$ a continuous function and for some positive constant $C$  depending on $n$ (thanks to \eqref{HH}).
Therefore, \eqref{equi} follows.

As far as the proof of \eqref{convNun} is concerned, we want to apply \cite[Theorem $2.1$ and Remark $2.2$]{BM}. To this aim, we need \eqref{convdebH},  \eqref{convun} as well as the $L^1$--estimate of $\{H(x, u_n, \nabla u_n)\}_n$.\\

Having \eqref{equi}, \eqref{convun} and \eqref{convNun}, we are allowed to apply  Vitali's Theorem  and conclude with \eqref{rhs}.
Finally \eqref{L1est1} follows from Fatou's Lemma and the a.e. convergences \eqref{convun} and \eqref{convNun}.
\end{pf}

\begin{Theorem}
Let $f \in L^m(\Omega)$ with $\frac{2N}{N+2} \le m < \frac{N}{2}$ and $ \alpha=\frac{N(q-1)-mq}{N-2m}$. Assume \eqref{Hy1}, \eqref{Hy2} and \eqref{Hy2.5}. Then, there exists at least a solution $u\in \ensp$ of \eqref{problem} in the sense of Definition \ref{fin-ener} such that $H(x, u, \nabla u)\in L^{2/q}(\Omega)$, $H(x, u, \nabla u)u\in L^1(\Omega)$ and
\begin{equation}\label{H}
\integrale |u|^{\si-2}|\N u|^2\,dx<M,
\end{equation}
that is, $|u|^\frac{\si}{2}\in \ensp$.
\end{Theorem}

\begin{pf}
We can take the limit in $n\to\infty$ in the approximating formulation \eqref{solappr} thanks to \eqref{convdebH}--\eqref{rhs}, recovering \eqref{finenform}.
The regularity \eqref{H} follows from \eqref{Hsigma}.
\end{pf}

\begin{Remark}\label{remark4}\rm
Having in mind Remark \ref{remark2}, we have a similar a priori estimate when $$\frac{N(q-1)-mq}{N-2m}<\alpha<q-1\,.$$ Thus, we may follow the proofs of Propositions \ref{convdebH2} and \ref{propequiintdge2} with this new exponent $\alpha$. We point out that we only need to check that
$$\frac{2(1-\alpha)}{2-q}\le \frac{2^*\sigma}2$$
which obviously holds with a bigger $\al$. Therefore, the above existence result applies as well.

The limit case $\alpha=q-1$ also holds taking into account the a priori estimate stated in Proposition \ref{prop2}.
\end{Remark}

\subsection{The case of renormalized solutions}

\begin{Proposition}\label{propr}
Let $f \in L^m(\Omega)$ with $1< m < \frac{2N}{N+2}$, $ \alpha=\frac{N(q-1)-mq}{N-2m}$ and $\{u_n\}_n$ be a sequence of solutions of \eqref{appr}. Assume also \eqref{Hy1}, \eqref{Hy2} and \eqref{Hy2.5}. Then, up to subsequences, there exists a function $u$ such that
\begin{equation}\label{convun2}
u_n\to u \qq \t{a.e. in }\Omega.
\end{equation}
\end{Proposition}
\begin{pf}
We claim that the uniform bound
\begin{equation}\label{Hsigma2}
\int_\Omega (1+|u_n|)^{\sigma-2}|\nabla u_n|^2 \,dx \le M
\end{equation}
holds. Indeed, Proposition \ref{teo2} applies with the same test function evaluated in $u_n$.\\
Now, set $1<r<2$ to be determined. Then, the above inequality allows us to estimate
\begin{equation*}
\int_\Omega |\nabla u_n |^r \,dx \le \left(\int_\Omega (1+|u_n|)^{\sigma-2}|\nabla u_n|^2\,dx \right)^{\frac{r}{2}} \left(\int_\Omega (1+|u_n|)^{\frac{r(2-\sigma)}{2-r}} \,dx \right)^{\frac{2-r}{2}}\,.
\end{equation*}

Requiring $\frac{r(2-\sigma)}{2-r}=r^*$ (that is $\frac{r(2-\sigma)}{2-r}=\frac{\si}{2}2^*$), we obtain $r=\frac{N\si}{N+\si-2}>1$ since $1<\si<2$. Note  that $r=m^*=\frac{N(q-1-\al)}{1-\al}$ which, for $\al=0$, becomes the exponent of the gradient regularity in \cite{GMP}. Since $\{u_n\}$ is bounded in $W_0^{1,r}(\Omega)$, an appeal to the compact embedding allows us to conclude \eqref{convun2}.
\end{pf}

\medskip

\begin{Proposition}\label{propr1}
Let $f \in L^m(\Omega)$ with $1< m < \frac{2N}{N+2}$, $ \alpha=\frac{N(q-1)-mq}{N-2m}$ and $\{u_n\}_n$ be a sequence of solutions of \eqref{appr}. Assume also \eqref{Hy1}, \eqref{Hy2} and \eqref{Hy2.5}. Then,
\begin{equation}\label{equi2.0}
H(x,u_n,\N u_n)\q\t{ is bounded in}\q L^m(\Omega)\,,
\end{equation}
and
\begin{equation}\label{equi2}
H(x,u_n,\N u_n)\q\t{ is equi--integrable in}\q L^1(\Omega)\,.
\end{equation}

Furthermore, up to subsequences, we have
\begin{equation}
\N u_n\to \N u\qq \t{a.e. in }\Omega\,,\label{convNun2}
\end{equation}
\begin{equation}
H(x,u_n,\N u_n)\to H(x,u,\N u)\qq \t{in }L^1(\Omega)\,,\label{rhs2}
\end{equation}
and, for all $j>0$.
\begin{equation}\label{convT}
T_j(u_n)\to T_j(u)\q\t{strongly in } \ensp \,.
\end{equation}
\end{Proposition}
\begin{pf}
Let us begin with the proof of \eqref{equi2.0}. Again, due to the assumption \eqref{Hy2} on $H(x, t,\xi)$ and to the regularity of $f$, we focus only on the gradient term.
Observe that, for some $\ga_0>\ga$, it holds that
\begin{align}
\int_\Omega \left|g(u_n)|\nabla u_n|^q\right|^b\,dx&\le
\gamma_0 \int_\Omega \frac{|\nabla u_n|^{qb}(1+|u_n|)^{\frac{bq(\sigma-2)}{2}}}{(1+|u_n|)^{\alpha b}(1+|u_n|)^{\frac{bq(\sigma-2)}{2}}}\,dx\nonumber\\
&\le \gamma_0\left(\int_\Omega |\nabla u_n|^2(1+|u_n|)^{\sigma-2}\,dx\right)^{\frac{qb}{2}} \left(\int_\Omega (1+|u_n|)^{\frac{bq(2-\sigma)}{2-bq}-\frac{2\alpha b}{2-bq}} \,dx\right)^\frac{bq-2}{2}\label{regrhs}
\end{align}
thanks to H\"older's inequality with $\left(\frac{2}{qb},\frac{2}{2-bq}\right)$. We impose
\begin{equation*}
\frac{bq(2-\si)}{2-bq}-\frac{2\al b}{2-bq}=2^*\frac{\si}{2}
\end{equation*}
and by \eqref{Hsigma2}, the integral \eqref{regrhs} is bounded.
Now, thanks also to the definitions of $\si=\si(q,\al)$ and $m=m(q,\al)$, we deduce
\[
b=\frac{N\si}{q(N-2+\si)-\al(N-2)}=\frac{N(q-1-\al)}{q-2\al}=m>1.
\]
Once we have obtained \eqref{regrhs}, then \eqref{equi2} follows by observing that
\[
\int_E |g(u_n)||\N u_n|^q\,dx\le |E|^{\frac{1}{m'}}\left(\integrale \left|g(u_n)|\N u_n|^q\right|^m\,dx\right)^\frac{1}{m}
\]
for every $E\subset \Omega$. \\

If, in particular, we take $E=\Omega$, then we have proved that the right hand side of \eqref{appr} is uniformly bounded in $L^1(\Omega)$ and this fact yields to \eqref{convNun2} thanks to \cite{BG} (see also \cite[Theorem $2.1$]{P}). Note that the limit function $u$ satisfies $|\N u|\in L^r(\Omega)$ with the same $r$ as in Proposition \ref{propr}.\\

Having \eqref{equi2}, \eqref{convun2} and \eqref{convNun2}, we are allowed to apply  Vitali's Theorem  and conclude with \eqref{rhs2}.\\

The uniform boundedness in \eqref{Hsigma2} implies that $T_j(u_n)$ is uniformly bounded in $\ensp$. We deduce the compactness of $T_j(u_n)$ in $\ensp$ from the compactness of the right hand side in $L^1(\Omega)$ (see  \cite{M} or \cite{LP}).
\end{pf}

\medskip

\begin{Theorem}\label{teoexsi<2}
Let $f \in L^m(\Omega)$ with $1 < m <\frac{2N}{N+2}$ and let $ \alpha=\frac{N(q-1)-mq}{N-2m}$.  Assume \eqref{Hy1}, \eqref{Hy2} and \eqref{Hy2.5} as well. Then, there exists at least a solution $u$ of \eqref{problem} in the sense of Definition \ref{renor} such that $H(x, u, \nabla u)\in L^m(\Omega)$ and
\begin{equation}\label{H2}
\integrale (1+|u|)^{\si-2}|\N u|^2\,dx<M,
\end{equation}
that is, $ (1+|u|)^{\frac{\si}{2}-1}u\in \ensp$.
\end{Theorem}

\begin{pf}
Consider in \eqref{solappr} a test function of the kind $S(u_n)\vp$, where $\vp\in H^{1}(\Omega)\cap L^\infty(\Omega)$ and $S:\R\to\R$ is a  Lipschitz function having compact support, say $\t{supp}(S(u_n))\subseteq [-j,j]$, and such that $S(u)\vp\in \ensp$. Then
\begin{equation*}
\integrale [A(x)\cdot\N u_n]\cdot\N (S(u_n)\vp)\,dx=\integrale H(x,u_n,\N u_n)S(u_n)\vp\,dx +\int_\Omega T_n(f) S(u_n)\vp\, dx.
\end{equation*}
Due to the support assumption on $S(u_n)$, the above equation only takes into account  $T_j(u_n)$, and so we rewrite the approximating formulation as
\begin{equation*}
\begin{array}{c}
\ds\integrale [A(x)\cdot\N T_j(u_n)]\cdot\N\vp\, S(u_n)\,dx+\integrale [A(x)\cdot\N T_j(u_n)]\cdot \N T_j(u_n) S'(u_n)\vp\,dx\\[4mm]\ds
=\integrale H(x,u_n,\N u_n)S(u_n)\vp\,dx+\int_\Omega T_n(f) S(u_n)\vp\, dx.
\end{array}
\end{equation*}
The convergence of the right hand side follows from \eqref{rhs2} and \eqref{convun2}. Furthermore
\[
\integrale [A(x)\cdot\N T_j(u_n)]\cdot\N\vp\, S(u_n)\,dx\to \integrale[A(x)\cdot\N T_j(u)]\cdot\N\vp\, S(u)\,dx
=\integrale [A(x)\cdot\N u]\cdot\N \vp\, S(u)\, dx
\]
and
\begin{align*}
\integrale [A(x)\cdot\N T_j(u_n)]\cdot\N T_j(u_n)  S'(u_n)\vp\,dx\to &\integrale [A(x)\cdot\N T_j(u)]\cdot \N T_j(u) S'(u)\vp\,dx
\\ & =\integrale [A(x)\cdot\N u]\cdot\N u S'(u)\vp \, dx
\end{align*}
thanks to \eqref{convun2} and \eqref{convT}.

We point out that  \eqref{equi2.0} and Fatou's lemma imply that $H(x, u, \nabla u)\in L^m(\Omega)$ holds.

The regularity \eqref{H2} directly follows from Proposition \ref{teo2} applied on \eqref{appr}.
\end{pf}

\begin{Remark}\label{remark5}\rm
As in Remark \ref{remark4}, we may consider exponents satisfying
$$\frac{N(q-1)-mq}{N-2m}<\alpha<q-1\,.$$
Indeed, it is enough to have in mind Remark \ref{al-beta} and follow the proofs of Propositions \ref{propr} and \ref{propr1} as well as Theorem \ref{teoexsi<2} with this new exponent $\alpha$. We point out that now we have to check that
\begin{equation*}
\frac{bq(2-\si)}{2-bq}-\frac{2\al b}{2-bq}\le 2^*\frac{\si}{2}
\end{equation*}
which obviously holds with a bigger exponent $b$. Therefore, the above existence result applies as well.

The limit case $\alpha=q-1$ also holds taking into account the a priori estimate stated in Proposition \ref{prop3}.
\end{Remark}

\subsection{The limit case}\label{ss3}

We have already analyzed the situation when $q>\frac{N+\alpha(N-2)}{N-1}$ with data $f\in L^m(\Omega)$ ($m>1$). It remains to study the limit case $q=\frac{N+\alpha(N-2)}{N-1}$, where existence of a renormalized solution with $L^1$--data should be expected. Nevertheless, this is not so as a variant of \cite[Example 4.1]{GMP} shows.

\begin{Example} Let $q=\frac{N+\alpha(N-2)}{N-1}$, and consider a nonnegative $f\in L^1(\Omega)$ and a continuous function $g\>:\>\R\to(0,+\infty)$ satisfying $g(s)=s^{-\alpha}$ for $s>k_0>0$.

Assume that there exists a renormalized solution $u$ to problem
\begin{equation*}
\left\{
\begin{array}{ll}
-\Delta u = g(u)\,|\nabla u |^q + f(x) & \mbox{ in } \Omega\,,\\[2mm]
u=0 & \mbox{ on } \partial\Omega\,,
\end{array}\right.
\end{equation*}
which is obviously nonnegative. Then $g(u)\,|\nabla u |^q\in L^1(\Omega)$, so that $g(k+G_k(u))\,|\nabla G_k(u) |^q\in L^1(\Omega)$ for all $k>0$. Fixing $k>k_0$, we deduce that
\[
\Big|\nabla \Big((k+G_k(u))^{1-\frac\alpha q}-k^{1-\frac\alpha q}\Big)\Big|^q\in L^1(\Omega)\,,
\]
that is,
\[
(k+G_k(u))^{1-\frac\alpha q}-k^{1-\frac\alpha q}\in W_0^{1,q}(\Omega)\,.
\]
Hence, the Sobolev embedding implies
$
(k+G_k(u))^{1-\frac\alpha q}\in L^{\frac{Nq}{N-q}}(\Omega)
$
and consequently it follows from $0\le u \le k+G_k(u)$ that
$u\in L^{\frac{N(q-\alpha)}{N-q}}(\Omega)$, where $q=\frac{N+\alpha(N-2)}{N-1}$. Observing that
\[
\frac{N(q-\alpha)}{N-q}=\frac{N}{N-2}\,,
\]
 it yields $u\in L^{\frac{N}{N-2}}(\Omega)$. To get a contradiction, we just need to compare with the unique renormalized solution of
\begin{equation*}
\left\{
\begin{array}{ll}
-\Delta v =  f(x) & \mbox{ in } \Omega\,,\\[2mm]
v=0 & \mbox{ on } \partial\Omega\,,
\end{array}\right.
\end{equation*}
which satisfies $0\le v\le u$ and so $v \in L^{\frac{N}{N-2}}(\Omega)$, but this summability does not hold for a general $L^1$-data.
\end{Example}

We may expect existence of solution to problem \eqref{problem} when we take $q=\frac{N+\alpha(N-2)}{N-1}$ and the datum belongs to the Orlicz space $L^1((\log L)^N)$. However, since we are focus in the setting of Lebesgue spaces, we must assume data $f\in L^m(\Omega)$ (with $m>1$) to deal with this limit case. Observe that it is enough to consider $1<m<\frac{N(q-1)}{q}$ due to embeddings in Lebesgue spaces. In this situation we have existence for a problem with exponent $\alpha_0=\frac{N(q-1)-mq}{N-2m}$. Owed to Remark \ref{remark5}, then we obtain an existence result for
\[
\alpha=\frac{N(q-1)-q}{N-2}>\frac{N(q-1)-mq}{N-2m}=\alpha_0\,.
\]
Therefore, we have proved the following result.

\begin{Theorem}
Let $f \in L^m(\Omega)$ with $m>1$ and let $\sigma=\frac{m(N-2)}{N-2m}$. Take $ \alpha=\frac{N(q-1)-q}{N-2}$ and  assume \eqref{Hy1}, \eqref{Hy2} and \eqref{Hy2.5}. Then, there exists at least a solution $u$ of \eqref{problem} in the sense of Definition \ref{renor} such that $ (1+|u|)^{\frac{\si}{2}-1}u\in \ensp$.

\end{Theorem}

\subsection{The case of measure data}\label{ss4}

We now discuss the case with measure data. Since we reason, as we have done before, through approximation techniques, we make some comments on the approximating problem we are going to consider.\\
Given $\mu\in \M(\Omega)$, we choose a sequence $\{\mu_n\}_n$ in $L^\infty(\Omega)$ which approximates $\mu$ as in \cite[Section 3]{DMOP} and satisfies
\begin{equation*}
\|\mu_n\|_{L^1(\Omega)}\le \|\mu\|_{\M(\Omega)}.
\end{equation*}
Now consider the following approximating problems of \eqref{problem_measure}:
\begin{equation}\label{appr_measure}
\left\{
\begin{array}{ll}
-\Div[A(x)\cdot\nabla u_n] = H(x, u_n, \nabla u_n) + \mu_n & \mbox{ in } \Omega\,,\\[2mm]
u_n=0 & \mbox{ on } \partial\Omega\,.
\end{array}\right.
\end{equation}
We already know that there exists solutions $u_n\in H_0^1(\Omega)\cap L^\infty(\Omega)$ to problem \eqref{appr_measure}.
We recall that the definition of the sequence $\{\mu_n\}_n$ in \cite[Section 3]{DMOP} is made in such a way that the following result holds.

\begin{Proposition}\label{conv_approx_mea}
Using the same notation as above, consider a Lipschitz--continuous function  $S\>:\>\R\to\R$ such that $S'$ has compact support and denote by $S(+\infty)$ and $S(-\infty)$
the limits of $S(t)$ at $+\infty$ and $-\infty$, respectively. Take  $\varphi\in W^{1,r}(\Omega)\cap L^\infty(\Omega)$, with $r>N$, such that $S(u)\varphi\in H_0^{1}(\Omega)$.

If, for some function $u$,
\begin{align*}
&u_n(x)\to u(x)\qquad \hbox{a.e. in } \Omega\\
&\nabla u_n(x)\to \nabla u(x)\qquad \hbox{a.e. in } \Omega\\
&T_k(u_n) \rightharpoonup T_k(u) \qquad \hbox{weakly in }H_0^1(\Omega)\quad\hbox{for all } k>0\\
&u_n\to u \qquad \hbox{strongly in }W_0^{1,s}(\Omega)\quad\hbox{for all }  1\le s<\frac {N}{N-1},
\end{align*}
then
\begin{equation*}
  \lim_{n\to\infty} \int_\Omega S(u_n)\varphi \mu_n\, dx=\int_\Omega S(u)\varphi\, d\mu_0+S(+\infty)\int_\Omega \varphi\, d\mu_s^+-S(-\infty)\int_\Omega \varphi\, d\mu_s^-\,.
\end{equation*}

\end{Proposition}

\begin{Proposition}\label{equiintL1}
Let $\mu \in \M(\Omega)$ have a norm small enough. Let $\frac{N(q-1)-q}{N-2} < \alpha <q-1$ and $\{u_n\}_n$ be a sequence of solutions of \eqref{appr_measure}. Assume also \eqref{Hy1}, \eqref{Hy2} and \eqref{Hy2.5}. Then, the a.e. convergences \eqref{convun2} and \eqref{convNun2}, the equi--integrability \eqref{equi2} and the strong convergences \eqref{rhs2}, \eqref{convT} and
\begin{equation}\label{W2s}
u_n \to u \quad\text{ in } W^{1,s}_0(\Omega)
\end{equation}
for all $1\le s < \frac{N}{N-1}$.
\end{Proposition}

\begin{pf}
Theorem \ref{teo3} implies that
\begin{equation}\label{esti}
\int_\Omega |\nabla T_k(u_n)|^2 \,dx \le Mk \,, \text{ for all } k>0\,,
\end{equation}
and then, using Lemma \ref{marc} we get
\begin{equation*}
|\{ |\nabla u_n| > k \}| \le C \frac{M^{\frac{N}{N-1}}}{k^{\frac{N}{N-1}}}\,, \text{ for all } k>0\,.
\end{equation*}
Hence, the sequence $\{u_n\}$ is bounded in $W^{1,s}_0(\Omega)$ for all $1\le s < \frac{N}{N-1}$ and there exist $u \in W^{1,s}_0(\Omega)$ and a subsequence (not relabelled) such that
\begin{align}
& u_n \rightharpoonup u \quad\text{ weakly in } W^{1,s}_0(\Omega)\,,\label{W1s}\\
& u_n \to u \quad\text{ in } L^s(\Omega)\,,\nonumber\\
& u_n \to u \quad\text{ a.e. in } \Omega\,.\label{conv_punt_med}
\end{align}
Moreover, condition \eqref{esti} also implies that
\begin{equation}\label{trunc_1}
\nabla T_k(u_n) \rightharpoonup \nabla T_k(u) \quad\text{weakly in } L^2(\Omega;\R^N)\,.
\end{equation}
To prove the equi-integrability of the right hand side we use that
\begin{equation}\label{gu}
\int_{\{|u|>k\}} g(u_n)|\nabla u_n|^q \,dx \le C_k \quad \hbox{for all }k>0\,,
\end{equation}
with $\lim_{k \to \infty} C_k=0$ (see Theorem \ref{teo3}). Thus, given $\e>0$ we may find $k_0>0$ such that
\begin{equation*}
\int_{\{|u|>k\}} g(u_n)|\nabla u_n|^q \,dx \le \frac{\e}{2} \quad\text{ for all }\; k\ge k_0 \;\text{ and for all }\; n \in \mathbb{N}\,.
\end{equation*}
Let $E \subset \Omega$ and let $k\ge k_0$ be fixed. Denoting $\bar g_k=\sup_{|s|\le k}|g(s)|$, the following inequalities hold:
\begin{align*}
\int_E g(u_n)|\nabla u_n|^q \,dx &= \int_{E\cap\{|u_n|\le k\}} g(u_n)|\nabla u_n|^q \,dx + \int_{E \cap \{|u_n|>k\}}g(u_n)|\nabla u_n|^q \,dx\\
& \le \bar{g}_k\int_E |\nabla T_k(u_n)| \,dx +\int_{\{|u|>k\}} g(u_n)|\nabla G_k(u_n)|^q\,dx \\
& \le \bar{g}_k \left(\int_\Omega |\nabla T_k(u_n)|^2\,dx \right)^{\frac{q}{2}}|E|^{1-\frac{q}{2}} + \frac{\e}{2}\\
& \le \bar{g}_k \left(Mk\right)^{\frac{q}{2}} |E|^{1-\frac{q}{2}} +\frac{\e}{2}\,,
\end{align*}
which goes to 0 when $|E|$ is small and so \eqref{equi2} is proved.

On the other hand, taking $E=\Omega$ and applying \cite[Lemma 1]{BM} we deduce \eqref{W2s}.
As a consequence we get
\begin{align*}
& \nabla u_n \to \nabla u \quad\text{ a.e. in } \Omega \,,\\
& g(u_n)|\nabla u_n|^q \to g(u)|\nabla u|^q \quad \text{ in } L^1(\Omega)
\,.
\end{align*}

This last convergence implies \eqref{rhs2}. Finally, appealing to the proof of  \cite[Theorem $3.4$, Step 5]{DMOP}, we deduce \eqref{convT}.

\end{pf}

\begin{Theorem}
Assume $\frac{N(q-1)-q}{N-2}<\al<q-1$, \eqref{Hy1}, \eqref{Hy2} and \eqref{Hy2.5}. If $\mu \in\M(\Omega)$ has a norm small enough, then there exists at least a renormalized solution to problem \eqref{problem_measure} in the sense of Definition \ref{measure}.
\end{Theorem}

\begin{pf}
We take advantage of the results contained in \cite[Theorem $2.33$]{DMOP}. \\
Proposition \ref{equiintL1} provides us with $u\in\mathcal{T}^{1,2}_0(\Omega)$ (by \eqref{trunc_1}), $u\in W_0^{1,s}(\Omega)$ for all $1\le s<\frac N{N-1}$ (by \eqref{W1s}) and $H(x,u,\nabla u)\in L^1(\Omega)$ (by \eqref{rhs2}).

Now, consider a Lipschitz--continuous function  $S\>:\>\R\to\R$ such that $S'$ has compact support and a $\varphi\in W^{1,r}(\Omega)\cap L^\infty(\Omega)$, with $r>N$, such that $S(u)\varphi\in H_0^{1}(\Omega)$. Since $S'$ has compact support, it follows that there exists $j>0$ such that $S'(u_n)\nabla u_n=S'(u_n)\nabla T_j(u_n)$. As a consequence, we get
\[\int_\Omega |\nabla S(u_n)|^2 dx=\int_\Omega S'(u_n)^2|\nabla T_j(u_n)|^2 dx\le \|S'\|_\infty \int_\Omega |\nabla T_j(u_n)|^2 dx<\infty\]
and so
\begin{equation}\label{truncad_1}
  S'(u_n)\nabla T_j(u_n) \rightharpoonup S'(u)\nabla T_j(u) \qquad\hbox{weakly in } L^2(\Omega; \R^N)\,.
\end{equation}

Taking $S(u_n)\varphi$ as test function in problem \eqref{appr_measure}, we have
\begin{multline}\label{conv_11}
  \int_\Omega \varphi S'(u_n)[A(x)\cdot\nabla u_n]\cdot \nabla u_n \, dx+ \int_\Omega  S(u_n)[A(x)\cdot\nabla u_n]\cdot \nabla \varphi \, dx \\
  =\int_\Omega  S(u_n)\varphi H(x, u_n, \nabla u_n) \, dx +\int_\Omega S(u_n)\varphi \mu_n\, dx\,.
\end{multline}
The first term on the left hand side can be written as
\[\int_\Omega \varphi S'(u_n)[A(x)\cdot\nabla u_n]\cdot \nabla u_n \, dx=\int_\Omega \varphi S'(u_n)[A(x)\cdot\nabla T_j(u_n)]\cdot \nabla T_j(u_n) \, dx\]
Hence, the strong convergence of $A(x)\cdot\nabla T_j(u_n)$ to $A(x)\cdot\nabla T_j(u)$ in $L^2(\Omega; \R^N)$ (due to \eqref{convT}) and the weak convergence of $S'(u_n)\nabla T_j(u_n)$ to $S'(u)\nabla T_j(u)$ (by \eqref{truncad_1}) imply the convergence of this first term to
\[\int_\Omega \varphi S'(u)[A(x)\cdot\nabla T_j(u)]\cdot \nabla T_j(u) \, dx=\int_\Omega \varphi S'(u)[A(x)\cdot\nabla u]\cdot \nabla u \, dx\,.\]
The convergence of the second term follows from \eqref{conv_punt_med}, the boundedness of function $S$ and \eqref{W2s}.

As far as the right hand side of \eqref{conv_11} is concerned, the convergence in the first term yields as a consequence of $S(u_n){\buildrel * \over\rightharpoonup} S(u)$ in $L^\infty(\Omega)$ and \eqref{rhs2}.  We deal with the second term applying Proposition \ref{conv_approx_mea}. Therefore, we can pass to limit in \eqref{conv_11} obtaining
\begin{multline*}
  \int_\Omega \varphi S'(u)[A(x)\cdot\nabla u]\cdot \nabla u \, dx+ \int_\Omega S(u)[A(x)\cdot\nabla u]\cdot \nabla \varphi \, dx \\
  =\int_\Omega S(u)\varphi H(x, u, \nabla u) \, dx +\int_\Omega S(u)\varphi\, d\mu_0+S(+\infty)\int_\Omega \varphi\, d\mu_s^+-S(-\infty)\int_\Omega \varphi\, d\mu_s^-\,.
\end{multline*}
Since we have proved one of the equivalent definitions of renormalized solution stated in \cite{DMOP}, we are done.

\end{pf}

Applying Proposition \ref{teo4} and Proposition \ref{teo5} in Proposition \ref{equiintL1}, instead of Theorem \ref{teo3}, it leads to Theorem \ref{lin1} and Theorem \ref{sub-lin1}, respectively.

\begin{Theorem}\label{lin1}
Let $\mu \in\M(\Omega)$ and $\al=q-1$. Assume \eqref{Hy1}, \eqref{Hy2} and \eqref{Hy2.5}, and $\ga$ is small enough.  Then, there exists at least a renormalized solution to problem \eqref{problem_measure} in the sense of Definition \ref{measure}.
\end{Theorem}

\begin{Theorem}\label{sub-lin1}
Let $\mu \in\M(\Omega)$ and $q-1< \al\le \frac q2$. Assume also \eqref{Hy1}, \eqref{Hy2} and \eqref{Hy2.5}.  Then, there exists at least a renormalized solution to problem \eqref{problem_measure} in the sense of Definition \ref{measure}.
\end{Theorem}

\section{Results on further regularity}

Throughout this paper, we have shown that the features of our problem can be shortened in the parameters $q$ and $\al$ and illustrated by a $(q,\al)$--plane (recall pictures in the introduction). In this way, the linear setting corresponds to the line $\al=q-1$, whereas the superlinear one coincides with a triangle whose sides are that line, $\al=0$ and $q=2$. In this triangle, points may be grouped according to the suitable summability of sources. Thus, between lines $\al=q-1$ and $\al=\frac{N-1}{N-2}q-\frac{N}{N-2}$, we may take measure data. In the remaining triangle the best line that enables data $f\in L^m(\Omega)$ is given by
\begin{equation*}
  \al=\frac{N-m}{N-2m}q-\frac{N}{N-2m}\,,\qquad 1<m<\frac N2\,.
\end{equation*}
The existence result also works if we have a bigger $\al$ (as was emphasized several times) or if a smaller $q$ is considered. On the other hand, fixed $q$ and $\al$, if we have more regular data, then existence is guaranteed by embeddings between Lebesgue spaces. In this case, however, the solution should be more regular as well.

In this Section, we propose an analogous of the bootstrapping results contained in  \cite{GMP}  in the case of being $g$ constant  in \eqref{Hy2}, that is, when $g(s)=\ga$ for all $s\in \R$.
It is worth comparing the results between the cases \cite[$g(s)\equiv \ga$]{GMP} and $g(s)\le \frac{\ga}{|s|^\al}$ below. In the constant case it is proved
\begin{itemize}
\item[R1.]  $u\in L^\infty(\Omega)$  if $f\in L^s(\Omega)$, $s>\frac{N}{2}$  (see  \cite[Theorems $3.8$ \& $4.9$ points $(i)$  \& Theorem $5.8$]{GMP}), while solutions to \eqref{problem} are bounded even  for $f\in L^r(\Omega)$ with $\frac{N}{2}<s<r$;
\item[R2.]  $u\in L^d(\Omega)$, $d<\infty$, if $f\in L^{\frac{N}{2}}(\Omega)$ (see  \cite[Theorems $3.8$ \& $4.9$ points $(ii)$ \& Theorem $5.8$]{GMP}) while, if \eqref{Hy2.5} is in force, then the same holds for greater values of $q$;
\item[R3.]  $|u|^{\frac{t}{2}}\in H_0^1(\Omega)$ if $f\in L^{r}(\Omega)$, $2^*<r<\frac{N}{2}$ (see  \cite[Theorems $3.8$ \& $4.9$ points $(iii)$ \& Theorem $5.8$]{GMP}),  while solutions to \eqref{problem} with $f\in L^{s}(\Omega)$ and (eventually) $2^*<s<r<\frac{N}{2}$ verify $|u|^{\frac{\tau}{2}}\in H_0^1(\Omega)$ for $\tau>t$;
\item[R4.] $(1+|u|)^{\frac{t}{2}-1}u\in H_0^1(\Omega)$ if $f\in L^{r}(\Omega)$, $\frac{N(q-1)}{q}<r\le 2^*$ (see  \cite[Theorems $4.9$ point $(iv)$ \& Theorem $5.8$]{GMP} ),  while solutions to \eqref{problem} with $f\in L^{s}(\Omega)$ (eventually) for $\frac{N(q-(1+\al))}{q-2\al}<s<\frac{N(q-1)}{q}<r\le 2^*$  verify $(1+|u|)^{\frac{\tau}{2}-1}u\in H_0^1(\Omega)$ for $\tau>t$.

\end{itemize}

In our case, both R1 and R2 directly follow from Proposition \ref{prop1}. So we will focus on R3 and R4. Notice that
the assumption \eqref{Hy2.5} (which generalise the constant case) in \eqref{Hy2} allows us to get the same results of the Theorems quoted above for lower data regularity/greater values of $q$ with respect to \cite{GMP}. Indeed the values $s$ in Theorems \ref{reg1} and \ref{reg2} can be taken
\[
\frac{N(q-1-\alpha)}{q-2\al}<s<\frac{N(q-1)}{q}.
\]

\subsection{The case $\frac{(N-2)(q-1-\alpha)}{2-q}\ge 2$}

\begin{Theorem}\label{reg1}
Let $ f\in L^s(\Omega)$ with   $s>m=\frac{N(q-1-\alpha)}{q-2\alpha}$ and $\al\frac{N-2}{N}+\frac{2}{N}+1\le q<2$ (i.e. $m\ge 2_*$ and  $\sigma=\frac{m(N-2)}{N-2m}=\frac{(N-2)(q-1-\alpha)}{2-q}\ge 2$). {Assume also \eqref{Hy1}, \eqref{Hy2} and \eqref{Hy2.5}.} Then, solutions $u$ to \eqref{problem}  in the sense of Definition \ref{fin-ener}, verifying \eqref{H} and $2^*<s<\frac{N}{2}$ satisfy $|u|^\frac{\tau}{2}\in H_0^1(\Omega)$, $u\in L^{2^*\tau}(\Omega)$ and
\begin{equation}\label{Ltau}
\||u|^\frac{\tau}{2}\|_{H_0^1(\Omega)}+\|u\|_{L^{s^{**}}(\Omega)}\le M \qq\t{for }\tau=\frac{s(N-2)}{N-2s}.
\end{equation}
{The constant $M$ depends on  $q,\,m,\,N,\,\lm,\,\ga,\,\al,\,|\Omega|$ and on $\|f\|_{L^s(\Omega)},\,\|u\|_{L^1(\Omega)}$. }
\end{Theorem}

Note that $\tau>\si=\frac{m(N-2)}{N-2m}$ since $s>m$.\\

\begin{pf}
We take $\vp_n(u)=\Phi_n(G_k(u))$ in \eqref{finenform} with
\[
\Phi_n(w)=\int_0^w |T_n(z)|^{\tau-2}\,dz,\qq \tau>\si\,.
\]
Then, the assumptions \eqref{Hy1}, \eqref{Hy2}, \eqref{Hy2.5} lead us to
\begin{multline}\label{disii}
\lm\int_\Omega|\N G_k(u)|^2|T_n(G_k(u))|^{\tau-2}\,dx\\
\le \ga\int_\Omega\frac{|\N G_k(u)|^q|\Phi_n(G_k(u))|}{|u|^\al}\,dx+\int_\Omega |f||\Phi_n(G_k(u))|\,dx\,.
\end{multline}
Observe that
\[
|\Phi_n(w)|
\le |T_n(w)|^{q\left(\frac{\tau}{2}-1\right)}\int_0^{w} |T_n(\xi)|^{(2-q)\left(\frac{\tau}{2}-1\right)}\,d\xi
\le |T_n(w)|^{q\left(\frac{\tau}{2}-1\right)}
\left(
\int_0^{w} |T_n(\xi)|^{\frac{\tau}{2}-1}\,d\xi
\right)^{2-q}
|w|^{q-1}
\]
by H\"older's inequality with indices $\left( \frac{1}{2-q},\frac{1}{q-1} \right)$. This fact allows us to estimate the first integral on the right hand side of \eqref{disii} as
\begin{multline*}
I=\int_\Omega\frac{|\N G_k(u)|^q|\Phi_n(G_k(u))|}{|u|^\al}\,dx\\
\le
\int_\Omega
\Biggl[
\left(|\N G_k(u)|^q|T_n(G_k(u))|^{q\left(\frac{\tau}{2}-1\right)}\right)
\left(
\int_0^{|G_k(u)|} |T_n(\xi)|^{\frac{\tau}{2}-1}\,d\xi
\right)^{2-q}|G_k(u)|^{q-1-\al}\Biggr]\,dx\,.
\end{multline*}
We now apply H\"older's inequality  with three indices $\left( \frac{2}{q},\frac{2^*}{2-q},\frac{N}{2-q} \right)$ in the right hand side above, so we get
\begin{align*}
I&\le
\left(\int_\Omega|\N G_k(u)|^2 |T_n(G_k(u))|^{\tau-2}\,dx
\right)^\frac{q}{2}
\left(\int_\Omega\left(
\int_0^{|G_k(u)|} |T_n(\xi)|^{\frac{\tau}{2}-1}\,d\xi
\right)^{2^*}\,dx
\right)^{\frac{2-q}{2^*}}\\
&\qq\times\left(\int_\Omega |G_k(u)|^{\frac{N(q-1-\alpha)}{2-q}}\,dx  \right)^{\frac{2-q}{N}}\,.
\end{align*}
Then, since
\[
\left(\int_\Omega\left(
\int_0^{G_k(u)} |T_n(\xi)|^{\frac{\tau}{2}-1}\,d\xi
\right)^{2^*}\,dx
\right)^{\frac{2-q}{2^*}}\le S_2^{\frac{2-q}2}
\left(\int_\Omega|\N G_k(u)|^2|T_n(G_k(u))|^{\tau-2}\,dx
\right)^\frac{2-q}{2}
\]
by Sobolev's embedding and $\frac{N(q-1-\al)}{2-q}=2^*\frac{\si}{2}$ (which is due to the definition of $\si$), we rewrite
\begin{equation}\label{disii2}
I\le
S_2^{\frac{2-q}2}\left(\int_\Omega|\N G_k(u)|^2|T_n(G_k(u))|^{\tau-2}\,dx
\right)
\left(\int_\Omega |G_k(u)|^{2^*\frac{\si}{2}}\,dx  \right)^{\frac{2-q}{N}}\,.
\end{equation}
Take $k_0$ such that
\begin{equation*}
S_2^{\frac{2-q}2}\ga \||G_k(u)|^{\frac{\si}{2}}\|_{L^{2^*}(\Omega)}^{\frac{2^*}{N}(2-q)}\le \frac{\lm}{2}\qq\forall k\ge k_0,
\end{equation*}
so that combining \eqref{disii}--\eqref{disii2} we obtain
\[
\frac\lm2\int_\Omega|\N G_k(u)|^2|T_n(G_k(u))|^{\tau-2}\,dx\le  \int_\Omega |f||\Phi_n(G_k(u))|\,dx
\]
for $k\ge k_0$.\\
Defining $\psi_n(w)=\int_0^w |T_n(z)|^{\frac{\tau}{2}-1}\,dz$, we are left with the study of
\[
\frac{\lm}{2}\int_\Omega|\N \psi_n(G_k(u))|^2\,dx\le  \int_\Omega |f||\Phi_n(G_k(u))|\,dx\,,
\]
and our current aim becomes finding the relation between $\Phi_n(\cdot)$ and $\psi_n(\cdot)$, in order to obtain an inequality only involving $\psi_n(\cdot)$.
Using the definition of $\psi_n$ and considering the sets $\{|w|<n \}$ and $\{|w|\ge n \}$, we get the inequality $|\psi_n(w)|\ge \frac2\tau |T_n(w)|^{\frac{\tau}{2}-1}|w|$. And then, since function $|T_n(z)|^{\tau-2} $ is non-decreasing, we deduce
\[
|\Phi_n(w)|\le |T_n(w)|^{\tau-2}|w|\le |T_n(w)|^{\tau-2}|w|\left(\frac{|w|}{|T_n(w)|}\right)^{\frac{\tau-2}{\tau}}=\Big[|T_n(w)|^{\frac{\tau}{2}-1}|w|\Big]^{2\frac{\tau-1}{\tau}}.
\]
Hence,
\[
|\Phi_n(w)|\le
 c  |\psi_n(w)|^{2\frac{\tau-1}{\tau}}
\]
for some constant  $c$, which does not depend on $n$.
This estimate and  H\"older's inequality with $(s,s')$ provide us with
\begin{multline}\label{utauen}
\frac{\lm}{2}\int_\Omega|\N \psi_n(G_k(u))|^2\,dx\le c \int_\Omega |f||\psi_n(G_k(u))|^{2\frac{\tau-1}{\tau}}\,dx\\
\le c \|f\|_{L^s(\Omega)}\left(
 \int_\Omega |\psi_n(G_k(u))|^{2s'\frac{\tau-1}{\tau}}\,dx
\right)^\frac{1}{s'}\,.
\end{multline}
Now, by Sobolev's embedding and since the definition of $\tau$ implies $2s'\frac{\tau-1}{\tau}=2^*$, we obtain
\begin{equation}\label{Phi2star}
\frac{\lm}{2S_2}\left(\int_\Omega| \psi_n(G_k(u))|^{2^*}\,dx\right)^{\frac{2}{2^*}-2\frac{\tau-1}{2^*\tau}}\le c \|f\|_{L^s(\Omega)}\,.
\end{equation}
We point out that
 $\frac{2}{2^*}-2\frac{\tau-1}{2^*\tau}=\frac{2}{2^*\tau}>0$.  Then we let $n\to\infty$ in \eqref{Phi2star} getting $u\in L^{2^*\frac{\tau}{2}}(\Omega)$.\\
We are now allowed to consider the limit $n\to\infty$ in \eqref{utauen}, deducing $|u|^\frac{\tau}{2}\in H_0^1(\Omega)$.
\end{pf}

{\begin{Remark}\label{uL1}\rm
The dependence of $M$ in \eqref{Ltau} on $\|u\|_{L^1(\Omega)}$ follows from the following fact.\\
Let us come back to \eqref{FM}. Recalling that the integral in the  right hand side of \eqref{eq10} is evaluated over $\{|u|>k \}$ and proceeding as in the proof of Proposition \ref{teo1}, we get the same results with  $C_3 \||f|\chi_{\{|u|>k\}}\|_{L^m(\Omega)}$ instead of $C_3 \|f\|_{L^m(\Omega)}$. Thus, we rewrite \eqref{FM} as
\[
F(Y_k) \le C_3 \||f|\chi_{\{|u|>k\}}\|_{L^m(\Omega)}\qq\forall k>0\,.
\]
Observe that H\"older's inequality gives
\[C_3\||f|\chi_{\{|u|>k\}}\|_{L^m(\Omega)}\le C_3 \|f\|_{L^s(\Omega)}\left(\frac{\|u\|_{L^1(\Omega)}}{k}\right)^{\frac{s-m}{sm}}\]
and this value will be less than $M^*$ for $k$ large enough.
This fact implies that the constants $M$ in Proposition \ref{teo1} depend (in this case) on $\|u\|_{L^1(\Omega)}$ too.\\
\end{Remark}}

\subsection{The case $1< \frac{(N-2)(q-1-\alpha)}{2-q} <2$}

\begin{Theorem}\label{reg2}
Let $ f\in L^s(\Omega)$ with  $s>m=\frac{N(q-1-\alpha)}{q-2\al}$ and $1+\frac{1}{N-1}+\al\frac{N-2}{N-1}<q<\al\frac{N-2}{N}+\frac{2}{N}+1$ (i.e. $1<m<2_*$ and $1<\si=\frac{m(N-2)}{N-2m}=\frac{(N-2)(q-1-\alpha)}{2-q}< 2$). Then, solutions $u$ to \eqref{problem}  in the sense of Definition \ref{renor}, verifying \eqref{H2} and
\begin{enumerate}[(i)]
\item $2_*\le s<\frac{N}{2}$ satisfy $|u|^\frac{\tau}{2}\in H_0^1(\Omega)$, $u\in L^{2^*\frac{\tau}{2}}(\Omega)$ and
\[
\||u|^\frac{\tau}{2}\|_{H_0^1(\Omega)}+\|u\|_{L^{2^*\frac{\tau}{2}}(\Omega)}\le M \qq\t{for }\tau=\frac{s(N-2)}{N-2s}\ge2>\si\,;
\]
\item $m<s< 2_*$ satisfy $(1+|u|)^{\frac{\tau}{2}-1}u\in H_0^1(\Omega)$, $|\N u|\in L^{s^*}(\Omega)$ and
\begin{equation}\label{Ltt}
\|(1+|u|)^{\frac{\tau}{2}-1}u\|_{H_0^1(\Omega)}+\||\N u|\|_{L^{s^*}(\Omega)}\le M \qq\text{for }\tau=\frac{s(N-2)}{N-2s}\in (\si,2)\,.
\end{equation}

\end{enumerate}
The constants $M$ above depend on  $q,\,s,\,N,\,\lm,\,\ga,\,\al,\,|\Omega|$ and on $\|f\|_{L^s(\Omega)},\,\|u\|_{L^1(\Omega)}$.
\end{Theorem}
Here, cases $(i)$ and $(ii)$ differ in how the interval $(\si,\infty)$ is split by the parameter $\tau$ as $\tau\ge2(>\si)$ and $\si<\tau< 2$, respectively.

\begin{pf}
 First consider a function $\Phi\in W^{1,\infty}(\R)$ satisfying
\begin{equation}\label{Phi1}
0\le \Phi'(w)\le L(1+|w|)^{\si-2},\qq L>0\,,
\end{equation}
and
\begin{equation}\label{Phi2.0}
|\Phi(w)|\le C(\Phi'(w))^{\frac{q}{2}} \int_0^{|w|}(\Phi'(z))^{\frac{2-q}{2}}\,dz,\qq C>0 \,.
\end{equation}
We remark that H\"older's inequality yields
\begin{equation}\label{Phi2}
|\Phi(w)|\le C(\Phi'(w))^{\frac{q}{2}} \int_0^{|w|}(\Phi'(z))^{\frac{2-q}{2}}\,dz \le C(\Phi'(w))^{\frac{q}{2}} \left[\int_0^{|w|}(\Phi'(z))^{\frac{1}{2}}\,dz\right]^{2-q}|w|^{q-1} \,.
\end{equation}

We take $\Phi(G_k(u))$ as test function. Note that the condition \eqref{Phi1} is needed in order to make the test function admissible.
Recalling \eqref{Hy1}, \eqref{Hy2} and \eqref{Hy2.5}, we have
\begin{equation}\label{back}
  \lm\int_\Omega \Phi'(G_k(u)) |\nabla G_k(u)|^2 \, dx
\le   \int_\Omega g(G_k(u))|\N G_k(u)|^q\Phi(G_k(u))\, dx+\int_\Omega |f| |\Phi(G_k(u))|\, dx\,.
\end{equation}
We use \eqref{Phi2} and H\"older's inequality with indices $\left(\frac2q, \frac{2^*}{2-q}, \frac{N}{2-q}\right)$ to estimate the integral involving the gradient term in the right hand side as
\begin{multline*}
 \int_\Omega g(u)|\N G_k(u)|^q\Phi(G_k(u))\, dx\le\ga \int_\Omega \frac{|\N G_k(u)|^q}{|u|^\alpha}\Phi(G_k(u))\, dx\\
 \le C\ga \int_\Omega |\N G_k(u)|^q(\Phi'(G_k(u)))^{\frac{q}{2}}
\left( \int_0^{|G_k(u)|}(\Phi'(z))^\frac{1}{2}\,dz \right)^{2-q}|G_k(u)|^{q-1-\alpha}\,dx\\
\le C\ga\left( \int_\Omega |\N G_k(u)|^2\Phi'(G_k(u))\,dx\right)^{\frac q2}\left[ \int_\Omega \left(\int_0^{|G_k(u)|}  (\Phi'(z))^\frac{1}{2}\,dz \right)^{2^*}\,dx\right]^{\frac{2-q}{2^*}}\\
\hskip3cm\times\left(
\int_\Omega |G_k(u)|^{\frac{N(q-1-\alpha)}{2-q}}\,dx
\right)^{\frac{2-q}{N}}\\
\le C\ga S_2^{2-q}\left( \int_\Omega |\N G_k(u)|^2\Phi'(G_k(u))\,dx\right)
\left(
\int_\Omega |G_k(u)|^{\frac\sigma2 2^*}\,dx
\right)^{\frac{2-q}{N}}\,,
\end{multline*}
thanks to Sobolev's embedding too. Now we choose $k_0$ such that
\begin{equation*}
\ga CS_2^{2-q}\||G_k(u)|^{\frac\sigma2}\|_{L^{2^*}(\Omega)}^{\frac{2^*(2-q)}{N}}\le \frac{\lm}{2}\qq\forall k\ge k_0\,.
\end{equation*}
\\
Going back to \eqref{back}, we have found that
\begin{equation}\label{fin}
\frac{\lm}{2}\int_\Omega |\N G_k(u)|^2\Phi'(G_k(u))\,dx\le \int_\Omega |f| |\Phi(G_k(u))|\, dx\,.
\end{equation}
If $\tau\ge2$, we argue as in Theorem \ref{reg1}. Instead, if $\tau<2$, we set
\[
\Phi(w)=\int_0^w (1+|z|)^{\si-2}(1+|T_n(z)|)^{\tau-\si}\,dz\,.
\]
It is straightforward that $\Phi$ satisfies  \eqref{Phi1}. We are showing that \eqref{Phi2.0} holds as well. To this end, we study the limits at $0$ and at $+\infty$:
\[\lim_{w\to0}\frac{|\Phi(w)|}{(\Phi'(w))^{\frac{q}{2}} \int_0^{|w|}(\Phi'(z))^{\frac{2-q}{2}}\,dz}=\frac{\tau(2-q)+2(q-1)}{2(\tau-1)}\,.\]
\[\lim_{w\to+\infty}\frac{|\Phi(w)|}{(\Phi'(w))^{\frac{q}{2}} \int_0^{|w|}(\Phi'(z))^{\frac{2-q}{2}}\,dz}=\frac{\sigma(2-q)+2q(q-1)}{2(\sigma-1)}\,.\]
Hence, \eqref{Phi2.0} follows.
We also consider
\begin{equation*}
\psi(w)=\int_0^{w}(1+|z|)^{\frac{\sigma-2}{2}}(1+|T_n(z)|)^{\frac{\tau-\si}{2}}\,dz
\end{equation*}
with $\tau=\frac{s(N-2)}{N-2s}>\si=\frac{m(N-2)}{N-2m}$  (i.e. $\tau 2^*=(2\tau-1)s'$ as in Theorem \ref{reg1}).

We are now considering the function
\[w\mapsto \frac{|\Phi(w)|}{|\Psi(w)|^{2\frac{\tau-1}{\tau}}}\]
to check that it is bounded.
It is easy to see that this quotient defines an even function which is increasing  in $\{0\le w\le n\}$ and decreasing in $\{w>n\}$. Since
\begin{equation}\label{final1}
\lim_{w\to\infty}\frac{\int_0^w(1+z)^{\tau-2}dz}{\left[\int_0^w(1+z)^{\frac{\tau-2}2}dz\right]^{2\frac{\tau-1}{\tau}}}
=\left(\frac{\tau}2\right)^{2\frac{\tau-1}{\tau}}\frac 1{\tau-1}\,,
\end{equation}
it follows that
\[|\Phi(w)|\le c|\Psi(w)|^{2\frac{\tau-1}{\tau}}\]
for certain constant $c$ not depending on $n$. Observe that the choice of the exponent involving $\tau$ in \eqref{final1} is justified to argue as in Theorem \ref{reg1}.
Indeed, an analogous inequality as \eqref{utauen} can be recovered reasoning in a similar way.

We just remark that the regularity $(1+|u|)^{\frac{\tau}{2}-1}u\in H_0^1(\Omega)$ in point $(ii)$ implies that
\[
\int_{\Omega}|\N u|^{s^*}\,dx\le
\left(\int_{\Omega}|\N u|^2(1+|u|)^{\tau -2}\,dx\right)^\frac{s^*}{2}
\left(\int_{\Omega}(1+|u|)^{2^*\frac{\tau}{2} }\,dx\right)^\frac{2-s^*}{2}\,,
\]
so \eqref{Ltt} follows.
\end{pf}

\begin{Remark}\label{uL12}\rm
An analogous of Remark \ref{uL1} holds in this case.
\end{Remark}

\subsection{The case $\frac{(N-2)(q-1-\alpha)}{2-q}<1$}

In this Subsection, we consider the case $m=1=\si$.

\begin{Theorem}\label{reg3}
Let $ f\in L^s(\Omega)$ with  $s>1$ and  $ \frac{N(q-1)-q}{N-2}<\alpha<q-1$ (i.e. $\al+1<q<1+\frac{1}{N-1}+\al\frac{N-2}{N-1}$). Assume also \eqref{Hy1}, \eqref{Hy2}, \eqref{Hy2.5} and that $u$ is a solution $u$ to \eqref{problem_measure}  in the sense of Definition  \ref{renor}. Then
\begin{enumerate}[(i)]
\item $2_*\le s<\frac{N}{2}$ satisfy $|u|^\frac{\tau}{2}\in H_0^1(\Omega)$, $u\in L^{2^*\frac{\tau}{2}}(\Omega)$ and
\[
\||u|^\frac{\tau}{2}\|_{H_0^1(\Omega)}+\|u\|_{L^{2^*\frac{\tau}{2}}(\Omega)}\le M \qq\t{for }\tau=\frac{s(N-2)}{N-2s}\ge 2\,;
\]
\item {  $1<s< 2_*$, then $(1+|u|)^{\frac{\tau}{2}-1}u\in H_0^1(\Omega)$, $|\N u|\in L^{s^*}(\Omega)$ and
\begin{equation*}
\|(1+|u|)^{\frac{\tau}{2}-1}u\|_{H_0^1(\Omega)}+\||\N u|\|_{L^{s^*}(\Omega)}\le M \qq\t{for }\tau=\frac{s(N-2)}{N-2s}\in (1,2)\,.
\end{equation*}
}
\end{enumerate}
The constants $M$ above depend on  $q,\,s,\,N,\,\lm,\,\ga,\,\al,\,|\Omega|$ and on $\|f\|_{L^s(\Omega)},\,\|u\|_{L^1(\Omega)}$.
\end{Theorem}
\begin{pf}
The proof follows the same argument of Theorem \ref{reg2}, although some changes in the case $1<\tau<2$ must be done. Indeed, we first choose $0<\rho<1$ and consider the function given by
\[\Phi(w)=\int_0^w\frac{(1+|T_n(z)|)^{\tau-1+\rho}}{(1+|z|)^{1+\rho}}dz\,,\]
which is a bounded function and so can be taken as test function (in the sense of Remark \ref{test}). It can be checked as in Theorem \ref{reg2} that it also satisfies condition \eqref{Phi2.0}. Arguing as above,  we arrive at the inequality
\begin{multline*}
 \int_\Omega g(u)|\N G_k(u)|^q\Phi(G_k(u))\, dx\le\\
C\ga\left( \int_\Omega |\N G_k(u)|^2\Phi'(G_k(u))\,dx\right)
\left(
\int_\Omega |G_k(u)|^{\frac{N(q-1-\alpha)}{2-q}}\,dx
\right)^{\frac{2-q}{N}}.
\end{multline*}
Next, we use the fact that
\[
\frac{N(q-1-\alpha)}{2-q}< \frac{2^*}2
\]
and take into account Remark \ref{beatle2}.
So, an inequality similar to \eqref{fin} may be deduced.

On the other hand, we set
\[\Psi(w)=\int_0^w\frac{(1+|T_n(z)|)^{\frac{\tau-1+\rho}2}}{(1+|z|)^{\frac{1+\rho}2}}dz\,.\]
Contrary to what happens in the above theorem, now the function
\[w\mapsto \frac{|\Phi(w)|}{|\Psi(w)|^{2\frac{\tau-1}{\tau}}}\]
is increasing in the whole interval $[0,+\infty)$. Nevertheless, it is not difficult to find a bound:
\[\frac{\frac{\tau-1+\rho}{\rho(\tau-1)}(1+n)^{\tau-1}}{  \left(\frac 2\tau\right)^{2\frac{\tau-1}{\tau}}\left[(1+n)^{\frac\tau2}-1\right]^{2\frac{\tau-1}{\tau}}}\le \left(\frac{\tau}{2}\right)^{2\frac{\tau-1}{\tau}}\frac{\tau-1+\rho}{\rho(\tau-1)}\left[\frac{2^{\frac\tau2}}{2^{\frac\tau2}-1}\right]^{2\frac{\tau-1}{\tau}}\]
and it does not depend on $n$. Hence, the analogous of \eqref{Phi2star} follows because
$2s'\frac{\tau-1}{\tau}=2^*$.
%
\end{pf}

\section*{Acknowledgements}
The third author has been partially supported  by the Spanish Ministerio de Ciencia, Innovaci\'on y Universidades and
FEDER, under project PGC2018--094775--B--I00.

  \makeatletter
    \providecommand\@dotsep{5}
  \makeatother

\end{document}

%% file: style.tex
\usepackage{amsmath,bm}
\usepackage{amsfonts}
\usepackage{amssymb}
\usepackage{graphicx}
\usepackage{xcolor}
\usepackage[utf8]{inputenc}
\usepackage{enumerate}

\usepackage{hyperref}
\hypersetup{colorlinks=true,linkcolor=red!70!black,citecolor=orange,urlcolor=blue}

\usepackage{chngcntr}

\numberwithin{equation}{section}

\newenvironment{pf}{\noindent{ Proof}.\enspace}{\rule{2mm}{2mm}\medskip}

\newtheorem{Theorem}{Theorem}[section]

\newtheorem{Definition}[Theorem]{Definition}
\newtheorem{Lemma}[Theorem]{Lemma}
\newtheorem{Proposition}[Theorem]{Proposition}

\newtheorem{Remark}[Theorem]{Remark}
\newtheorem{Example}[Theorem]{Example}

\usepackage{tikz-cd}
\usepackage{caption}

\usepackage{float}
\usepackage{xcolor}
\usepackage{pgfplots}
\usetikzlibrary{tikzmark,fit,shapes.geometric}
\usepackage{tikz-cd}
\usetikzlibrary{shapes.geometric,calc,shapes,matrix,arrows,fit,decorations.text}

\usepackage{tkz-graph}
\usetikzlibrary{patterns}
\usepackage{pgfplotstable}
\usetikzlibrary{intersections}

\allowdisplaybreaks

\usepackage{xcolor}

%% file: macro.tex
\newcommand{\R}{\mathbb{R}}
\newcommand{\N}{\mathbb{N}}
\newcommand{\e}{\varepsilon}
\newcommand{\sg}{{\rm \; sign}}
\newcommand{\Div}{\hbox{\rm div\,}}

\long\def\salta#1{\relax}
\def\eps{\varepsilon}
\def\N{\nabla}
\def\ga{\gamma}
\def\al{\alpha}
\def\vp{\varphi}

\def\dive{\t{div }}
\newcommand{\lm}{\lambda}
\newcommand{\si}{\sigma}
\newcommand{\D}{\Delta}
\def\q{\quad}
\def\qq{\qquad}
\newcommand{\integrale}{\int_\Omega}

\definecolor{ncs}{rgb}{0.0, 0.53, 0.74}

\def\t{\text}
\newcommand{\ds}{\displaystyle}
\def\ensp{H_0^1(\Omega)}
\def\M{\mathcal{M}_b}

\newcommand{\pare}[1]{\left( #1 \right)}
\newcommand{\bra}[1]{\left[ #1 \right]}

%% file: todo.tex
\usepackage{xargs}                      
\usepackage{xkeyval}
\usepackage[colorinlistoftodos,prependcaption,textsize=small]{todonotes}
\usetikzlibrary{positioning}

\newcounter{todocountert}

\newcounter{todocounterq}

\newcounter{todocountera}

%% file: 2019_10_07.bbl
\begin{thebibliography}{plain}
\bibitem{AFM} A. Alvino, V. Ferone \& A. Mercaldo, {Sharp a priori estimates for a class of nonlinear elliptic equations with lower order terms}, Ann. Mat. Pura Appl. (4) 194 (2015), no. 4, 1169--1201.



\bibitem{B-V} Ph. B\'enilan, L. Boccardo, Th. Gallou\"et, R. Gariepy, M. Pierre \& J.L. V\'azquez, \href{http://www.numdam.org/item/?id=ASNSP_1995_4_22_2_241_0}{An $L^1$--theory of existence and uniqueness of solutions of nonlinear elliptic equations}, Ann. Scuola Norm. Sup. Pisa Cl. Sci. (4) {\bf 22} (1995), no. 2, 241–-273.

    \bibitem{BBM} A. Bensoussan, L. Boccardo \& F. Murat, \href{https://www.sciencedirect.com/science/article/pii/S0294144916303420?via\%3Dihub}{On a nonlinear partial differential equation having natural growth terms and unbounded solution} Ann. Inst. Henri Poincaré, Anal. Non Linégreene 5, No. 4, (1988) 347--364.

\bibitem{BG} L. Boccardo \& T. Gallou\"et, \href{https://www.tandfonline.com/doi/abs/10.1080/03605309208820857?journalCode=lpde20}{Nonlinear Elliptic Equations with Right Hand Side Measures}, Comm. in Par. Diff. Eq., Vol. {\bf 17} (1992), pp. 189--258.

    \bibitem{BGO} L. Boccardo, T. Gallou\"et \& L. Orsina, \href{https://www.sciencedirect.com/science/article/pii/S0294144916301135?via\%3Dihub}{Existence and uniqueness of entropy solutions for nonlinear elliptic equations with measure data}, Ann. Inst. Henri Poincaré, Anal. Non Linégreene 13, No. 5, (1996), 539--551.

\bibitem{BM} L. Boccardo \& F. Murat, \href{https://ac.els-cdn.com/0362546X92900238/1-s2.0-0362546X92900238-main.pdf?_tid=f76725e4-f080-444a-9f17-7086ee790125&acdnat=1529828053_ce1dced2fecff821b3718c76bc5bcdde}{Almost everywhere convergence of the gradients of solutions to elliptic and parabolic equations
}, Nonlin. Anal. TMA, Vol. {\bf 19} (1992), pp. 581--597.

\bibitem{BMP1} L. Boccardo, F. Murat \& J.--P. Puel, \href{https://eudml.org/doc/115517}{Existence de solutions non born\'ees pour certaines \'equations quasi--lin\'eaires} Port. Math. 41, (1982) 507--534.

    \bibitem{BMP2} L. Boccardo, F. Murat \& J.--P. Puel, \href{http://www.numdam.org/item/?id=ASNSP_1984_4_11_2_213_0}{R\'esultats d’existence pour certains probl\`emes elliptiques quasilin\'eaires} Ann. Sc. Norm. Super. Pisa, Cl. Sci., IV. Ser. 11, (1984) 213--235.


\bibitem{BMP3} L. Boccardo, F. Murat \& J.--P. Puel, \href{https://link.springer.com/article/10.1007\%2FBF01766148}{Existence of bounded solutions for non linear elliptic unilateral problems} Ann. Mat. Pura Appl., IV. Ser. 152, (1988) 183--196.

\bibitem{BMP4} L. Boccardo, F. Murat \& J.--P. Puel, \href{https://epubs.siam.org/doi/10.1137/0523016}{$L^\infty$ estimate for some nonlinear elliptic partial differential equations and application to an existence result} SIAM J. Math. Anal. 23, No. 2, (1992) 326--333.


\bibitem{BMa} G. Bottaro \& M.E. Marina, {Problema di Dirichlet per equazioni ellittiche di tipo variazionale su insiemi non limitati}
Boll. Unione Mat. Ital., IV. Ser. 8, (1973) 46--56.


\bibitem{DMOP} G. Dal Maso, F. Murat, L. Orsina \& A. Prignet, \href{http://www.numdam.org/item/ASNSP_1999_4_28_4_741_0}{Renormalized solutions of elliptic equations with general measure data}, Ann. Scuola Norm. Sup. Pisa Cl. Sci. (4) {\bf 28} (1999), no. 4, 741-–808.

 \bibitem{DVP}  T. Del Vecchio \& M.M. Porzio, {Existence results for a class of non coercive Dirichlet problems} Ric. Mat. 44, No. 2, (1995) 421--438.

 \bibitem{FM} V. Ferone \& F. Murat, {Quasilinear problems having quadratic growth in the gradient: An existence result when the source term is small} \'Equations aux d\'eriv\'ees partielles et applications. Articles d\'edi\'es \`a Jacques--Louis Lions. Gauthier--Villars: Paris. (1998) 497--515.

      \bibitem{FM1} V. Ferone \& F. Murat, \href{https://www.sciencedirect.com/science/article/pii/S0362546X99001650?via\%3Dihub}{Nonlinear problems having natural growth in the gradient: an existence result when the source terms are small} Nonlinear Anal., Theory Methods Appl., Ser. A, Theory Methods 42, No. 7, (2000) 1309--1326.

    \bibitem{FM2} V. Ferone \& F. Murat, \href{https://www.sciencedirect.com/science/article/pii/S0022039613004270?via\%3Dihub}{Nonlinear elliptic equations with natural growth in the gradient and source terms in Lorentz spaces} J. Differ. Equations 256, No. 2, (2014) 577--608.


     \bibitem{GMP0} N. Grenon, F. Murat \& A. Porretta, \href{https://www.sciencedirect.com/science/article/pii/S1631073X05004231?via\%3Dihub}{
Existence and a priori estimate for elliptic problems with subquadratic gradient dependent terms}  C. R., Math., Acad. Sci. Paris 342, No. 1, (2005) 23--28.

\bibitem{GMP} N. Grenon, F. Murat \& A. Porretta, \href{http://annaliscienze.sns.it/index.php?page=Article&id=300}{A priori estimates and existence for elliptic equations with gradient dependent terms}, Ann. Sc. Norm. Super. Pisa Cl. Sci., {\bf 13}, Issue 1, (2014), pp. 137--205.

  \bibitem{LP}  Ch. Leone \& A. Porretta,
  \href{https://www.sciencedirect.com/science/article/pii/S0362546X96003239?via\%3Dihub}{Entropy solutions for nonlinear elliptic equations in $L^1$}
Nonlinear Anal., Theory Methods Appl. {\bf 32}, Issue 3, (1998), 325--334.

\bibitem{LL} J. Leray \& J.--L. Lions, \href{http://www.numdam.org/item/?id=BSMF_1965__93__97_0}{Quelques r\'esultats de Visik sur les probl\`emes elliptiques non lin\'egreenes par les m\'ethodes de Minty--Browder}, Bull. Soc. Math. Fr. 93, (1965) 97--107.

\bibitem{salva} S. López-Martínez, A singularity as a break point for the multiplicity of solutions to quasilinear elliptic problems,  to appear in Adv. Nonlin. Anal..

\bibitem{M} F. Murat, \href{http://archive.schools.cimpa.info/anciensite/NotesCours/PDF/2009/Alexandrie_Murat_2.pdf}{Soluciones renormalizadas de EDP elipticas no lineales}, Laboratoire d’Analyse Numérique de l’Université Paris
VI, Technical report R93023, (1993).

\bibitem{P} A. Porretta, \href{http://citeseerx.ist.psu.edu/viewdoc/download?doi=10.1.1.432.2611&rep=rep1&type=pdf}{Some remarks on the regularity of solutions for a class of elliptic equations with measure data}, Houston J. Math., Vol. {\bf 26} (2000), pp. 183-–213.

    \bibitem{P1} A. Porretta, \href{https://eudml.org/doc/127132}{Nonlinear equations with natural growth terms and measure data}, Electron. J. Differ. Equ. 2002, Conf. 09, (2002) 183--202.

\bibitem{PS} A. Porretta \& S. Segura de Le\'on, \href{http://www.sciencedirect.com/science/article/pii/S0021782405001212?via\%3Dihub}
 {Nonlinear elliptic equations having a gradient term with natural growth}, J. Math. Pures Appl. (9) {\bf 85}, No 3, (2006), pp. 465--492.

 \bibitem{S} S. Segura de Le\'on, {}{Existence and uniqueness for $L^1$ data of some elliptic equations with natural growth}, Adv. Differ. Equ. 8, No. 11,  (2003) 1377--1408.


\end{thebibliography}
